\documentclass[final,1p,times,twocolumn]{elsarticle}
\usepackage{natbib}

\usepackage{graphicx}
\usepackage{subfigure}

\usepackage{amsmath}
\usepackage{amssymb}

\usepackage{color}
\usepackage{bm}

\usepackage{booktabs}
\usepackage{multirow}
\usepackage{array}

\newcommand{\norm}[1]{\left\Vert#1\right\Vert}
\newcommand{\pref}[1]{(\ref{#1})}
\newtheorem{theorem}{Theorem}[section]
\newtheorem{remark}[theorem]{Remark}

\begin{document}

\begin{frontmatter}



\title{Implicit Progressive-Iterative Approximation for Curve and Surface Reconstruction}


\author[add_math]{Yusuf Fatihu Hamza}
\author[add_math,add_cad]{Hongwei Lin\corref{cor1}}
\ead{hwlin@zju.edu.cn}
\author[add_cad]{Zihao Li}
\cortext[cor1]{Corresponding author:
  Tel.: +86-571-87951860-304;
  fax: +86-571-87953867
	}
\address[add_math]{School of Mathematics, Zhejiang University, Hangzhou, 310027, China}
\address[add_cad]{State Key Laboratory of CAD\&CG, Zhejiang University, Hangzhou, 310058, China}

%\author[1,2]{Hongwei Lin \corref{cor1}}
%\cortext[cor1]{Corresponding author:
%  Tel.: +86-571-87951860-304;
%  fax: +86-571-87953867}
%  \ead{hwlin@zju.edu.cn}

\begin{abstract}
Implicit curve and surface reconstruction attracts the attention of many
     researchers and gains a wide range of applications,
     due to its ability to describe objects with complicated geometry and topology.
However, extra zero-level sets or spurious sheets arise
    in the reconstruction process makes the reconstruction result challenging to be interpreted and damage the final result.
In this paper, we proposed an implicit curve and surface reconstruction
     method based on the progressive-iterative approximation method,
     named implicit progressive-iterative approximation (I-PIA).
The proposed method elegantly eliminates the spurious sheets naturally
    without requiring any explicit minimization procedure,
    thus reducing the computational cost greatly and providing high-quality reconstruction results.
Numerical examples are provided to demonstrate the efficiency
    and effectiveness of the proposed method.
\end{abstract}

\end{frontmatter}

% main text

%\begin{figure}
%\includegraphics[width=1.0\linewidth]{buddha8k150a}
%\caption{ Bunny by LSPIA. }

%------------------------------------------------------------------------------
% Section: Introduction
%------------------------------------------------------------------------------
\section{Introduction}
\label{sec:introduction}

Implicit representation and parametric representation are two common
    representation techniques in geometric design.
With parametric representations,
    it is difficult to fit a data set with complicated geometry,
    and parametrization is always a challenging problem.
Without requiring any parametrization,
    the implicit function can describe an object with complicated geometry and supply flexible and smooth surface representation.
Thus, implicit surface reconstruction receives great attention
    due to its capability to create an object with complicated topology and geometry.

However, the extra zero-level sets generated in the implicit curve and
    surface reconstruction procedure make the reconstruction results challenging to be interpreted and damage the resulting curve and surface.
To eliminate the extra zero-level sets,
    regularization terms are usually required to be added in the objective functions of minimization problems for implicit curve and surface reconstruction.
For example, Liu~\cite{Liu2017Implicit} incorporated the total variation
    of implicit representation to reduce the appearance of the extra zero-level sets as minimum as possible.
In Refs. ~\cite{Juttler2002Least,Rouhani2011Implicit,Rouhani2015Implicit,Yang2005Fitting},
    tension terms are added in the minimization problem to get rid of extra zero-level sets and avoid a singular system of equations.
Because the forms of the regularization terms are usually complicated,
    their addition to the objective functions seriously affects the efficiency of the implicit curve and surface reconstruction algorithms.

On the other hand, progressive-iterative approximation (PIA) is a series of
    efficient data fitting methods with intuitive geometric meaning.
They have been extensively employed in parametric curve and surface fitting,
    and subdivision curve and surface fitting \cite{Chen2008Subdivision,Cheng2009Loop,Deng2014Progressive,Deng2014Weighted,Lin2010Local,
    Lin1707convergence,Lin2004non-uniform,Lin11extended,Lin2005Totally,
    Lu2010Weighted},
    but never been used in implicit curve and surface reconstruction.
In this paper, we developed a progressive-iterative approximation method for
    implicit curve and surface reconstruction,
    named implicit progressive-iterative approximation (I-PIA).
We proved the convergence of I-PIA,
    and showed that the I-PIA method itself naturally solves a minimization
    problem with a regularization term,
    without any extra computation effort.
Therefore, not only the I-PIA method can eliminate the extra level sets,
    but, more importantly,
    it can improve the reconstruction efficiency of the implicit curves and surfaces.
Lots of numerical examples illustrated in this paper show that,
    by I-PIA, the implicit curve and surface reconstruction time is improved one to three orders of magnitude,
    compared with the state-of-the-art implicit curve and surface reconstruction methods.
In conclusion, the main contributions of this paper include:
\begin{itemize}
    \item I-PIA naturally solves a minimization problem with a
        regularization term in its iteration procedure,
        without any extra computation effort.
    \item No extra zero-level sets exist in the iteration of I-PIA,
        so the reconstruction result by I-PIA is very clear.
    \item The implicit curve and surface reconstruction time by I-PIA
        is improved at least one to three orders of magnitude,
        compared with the state-of-the-art methods.
\end{itemize}

The structure of this paper is as follows.
In Section \ref{sec:Related work}, we reviewed some related work.
Preliminary definitions and statements of the problems are given in Section
 \ref{sec:Fitting}.
In Section  \ref{sec:I-PIA}, we presented the I-PIA
 method for implicit curve and surface reconstruction.
Experimental results and discussions are given in Section \ref{sec:Numerical Examples}.
Finally, we conclude the paper in Section \ref {sec:conclude}.

%------------------------------------------------------------------------------
% Subsection: Related work
%------------------------------------------------------------------------------
\subsection{Related work}
\label{sec:Related work}

 In this section, some work related to implicit curve and surface
    reconstruction and progressive iterative approximation (PIA) will be briefly reviewed.

  \textbf{Implicit curve and surface reconstruction:} Most implicit surface reconstruction algorithms blend local implicit primitives to represent surfaces based on the idea developed by Blinn \cite{Blinn82generalization}. Muraki \cite{Muraki1991Volumetric} developed the Blobby model to fit a very complicated data set by blending implicit primitives. Hoppe et al. \cite{Hoppe1992Surface} proposed a surface reconstruction algorithm based on the locally defined signed distance function. Curless and Levoy \cite{Curless96volumetric} used the volumetric representation consisting of a cumulative weighted signed distance function. Carr et al. \cite{Carr2001Reconstruction} proposed a fast method for fitting and evaluating Radial basis functions
(RBFs) to model a large data set.

So far, many schemata had been considered to reduce the amount of storage cost in implicit      reconstruction.
	Morse \cite{Morse2005Interpolating} proposed using compactly supported radial basis functions to reduce the computational cost and memory requirement, which allows reconstructing surface from large data-sets that are impractical  from the previous method like \cite{Turk1998Variational}. Furthermore, Refs.~\cite{Kojekine2003Software,ohtake2005multi} consider compactly supported RBF to reduce the computational cost and improve the efficiency of the reconstruction process. Pan et al. \cite{Pan2016Compact} incorporated a low-rank tensor approximation technique and reduced the storage requirement efficiently. %Rouhani et al. \cite{Rouhani2015Implicit} proposed an efficient version of implicit B-spline (IBS) referred to truncated IBS, that only finds the active coefficients, which significantly reduced the computational cost and storage space.

Ohtake et al. \cite{Ohtake2003Multi-level} proposed a multi-level partition of unity (MPU) representation to reconstruct surface models from a huge set of points. Wang et al. \cite{Wang2011Parallel} presented a surface reconstruction algorithm based on the implicit PHT-spline, which reconstructs a surface from a large point cloud efficiently.
The Poisson surface reconstruction proposed by Kazhdan et al. \cite{Kazhdan2006Poisson} expresses the surface reconstruction as a Poisson problem and approximates the indicator function of the surface.
Moreover, screened Poisson surface reconstruction was presented to avoid the over-smoothing by incorporating the positional constraints in the optimization problem \cite{Kazhdan2013Screened}.

\textbf{Progressive iterative approximation (PIA):} PIA is widely used for its ability to fit data efficiently. PIA method elegantly generates a sequence of curves/surfaces by refining the control points of blending curves/surfaces, and the data points are interpolated by the limit of the sequence.
Both the numbers of control points and data points are required to be equal in the classical PIA.
With the recent advancement of the big data era,
    it is infeasible to fit large-scale data points by PIA.
To overcome this drawback, Lin and Zhang \cite{Lin11extended} developed an extended PIA (EPIA) method, which allows the number of data points to be higher than the number of control points. Progressive and iterative approximation for least square fitting (LSPIA) \cite{Deng2014Progressive} is another elegant PIA method, which allows the number of data points to be higher than the number of control points, and its limit is the least square fitting result to a given data set.

Initially, Qi et al. \cite{Qi1975numeric} and de Boor \cite{deBoor} discovered the property of PIA for uniform cubic B-spline curve.
Later, Lin et al. \cite{Lin2004non-uniform} showed that non-uniform B-spline curve and surface have PIA property.
Moreover, PIA property holds for curves and surfaces with normalized totally positive (NTP) basis \cite{Lin2005Totally}.
Also Rational B-spline curve and surface posses this property \cite{Shi2006terative}.
%Delgado and pe?a \cite{Delgado2007Progressive} examine the convergence rate of different NTP basis and prove that normalized B-basis has the fastest convergence rate.
Lu \cite{Lu2010Weighted} proposed a weighted PIA technique to increase the convergence rate of the PIA method.
For more details, a comprehensive overview of PIA is provided in \cite{Lin2018Survey}.

 As stated above, some PIA methods have been developed for both parametric
     curve and surface fitting,
     and subdivision curve and surface fitting,
     but never been used in implicit curve and surface reconstruction.

%------------------------------------------------------------------------------
% Section: Definitions and preliminaries
%------------------------------------------------------------------------------		
\section{Definitions and preliminaries }
\label{sec:Fitting}

In this  section,
    statement of the problems, definition of implicit B-spline curve and surface are given.
    Specifically, the following problems will be handled in this paper.

 \textbf{Implicit curve reconstruction problem:}
 Given a collection of unorganized data  points in the two-dimensional space,
\begin{equation}\label{eq:Datapnt_2d}
\{\bm{p}_i=(x_i,y_i), i=1,2,...,n\},
\end{equation}
with a set of associated oriented unit normals $\{\bm{n}_i, i=1,2,...,n\}$,
 find a function $f(x,y)$ so that the zero level sets of $f(x,y)$,
    i.e., $f(x,y) = 0$, fit the unorganized point set~\pref{eq:Datapnt_2d}.

 \textbf{Implicit surface reconstruction problem:}
 Given a collection of unorganized data points in the three-dimensional space,
\begin{equation}\label{eq:Datapnt_3d}
\{p_i=(x_i,y_i,z_i), i=1,2,...,n\},
\end{equation}
with a set of associated oriented unit normals $\{\bm{n}_i, i=1,2,...,n\}$,
 find a function $f(x,y,z)$ so that the zero level sets of $f(x,y,z)$,
    i.e., $f(x,y,z) = 0$, fit the unorganized point set~\pref{eq:Datapnt_3d}.

 Let $f(x,y)$ be  a bivariate tensor-product B-spline
function of bi-degree $(d_1,d_2)$  defined over some domain $\Omega$ \cite{Liu2017Implicit,Yang2005Fitting}:
\begin{equation}\label{eq:implicit_curve}
f(x,y) = \sum_{i=1}^{N_u}\sum_{j=1}^{N_v}C_{ij}B_i(x)B_j(y),
\end{equation}
where $C_{i,j}$ are the control coefficients and $B_i(x), B_j(y)$ are the B-spline basis functions  with some given knot sequences.
We consider the bi-cubic B-spline functions,
    i.e. $d_1=d_2=3$ with uniform  knot sequences.
The implicit B-spline curve reconstructed by fitting
    the data point set~\pref{eq:Datapnt_2d} is given by
\begin{equation}\label{eq:zeroset_2d}
z_f=\{(x,y)\in\Omega:f(x,y)=0\}.
\end{equation}

Similarly define the trivariate tensor product B-spline function
    $f(x,y,z)$ on some domain $\Omega$ as \cite{Liu2017Implicit,Yang2005Fitting}:
\begin{equation}
f(x,y,z) = \sum_{i=1}^{N_u}\sum_{j=1}^{N_v}\sum_{k=1}^{N_w}C_{ijk}B_i(x)B_j(y)B_k(z),
\end{equation}
where $C_{ijk}$ are the control coefficients,
    and $B_i(x),$ $B_j(y),$ $B_k(z)$ are the B-spline
 basis functions.
In our implementation, the tri-cubic B-spline function is employed.
Analogously, the implicit B-spline surface reconstructed by fitting
    the data point set~\pref{eq:Datapnt_3d} is,
\begin{equation}\label{eq:zeroset_3d}
z_f=\{(x,y,z)\in\Omega:f(x,y,z)=0\}.
\end{equation}

 For simplicity,
    we only present the method for implicit curve reconstruction in the following,
    and that for the implicit surface reconstruction is similar.
 Actually, the implicit curve is reconstructed by minimizing the sum of the
    squared algebraic distances~\cite{Juttler2002Least}, i.e.,
\begin{equation}\label{eq:minprob}
 \text{min} \ E(C_{11},C_{12},\cdots,C_{N_uN_v})=\sum_{i=1}^{n}f^2(p_i),
\end{equation}
    where $C_{ij},i =0,1,\cdots,N_u, j=0,1,\cdots,N_v$ are the coefficients of the tensor product B-spline function.
 However, in the implicit curve and surface reconstruction problem,
    the minimization problem~\pref{eq:minprob} is usually underdetermined,
    i.e., the number of unknowns is larger than that of the data points,
    thus leading to extra zero level sets in the reconstruction results.

 To make the undetermined problem~\pref{eq:minprob} determined,
    and eliminate the extra zero level sets,
    global regularization terms are added in the minimization problem~\pref{eq:minprob} \cite{Juttler2002Least,Rouhani2011Implicit,Rouhani2015Implicit,Yang2005Fitting}.
 However, the addition of the regularization terms seriously affects
    the efficiency of implicit curve and surface reconstruction.
 In this paper, we developed the I-PIA method for implicit curve and surface
    reconstruction,
    which solves the minimization problem with regularization terms naturally,
    and improves the efficiency of implicit curve and surface reconstruction significantly,
    while eliminating the extra zero level sets.

 \vspace{0.5cm}

\section{Implicit progressive iterative approximation}
\label{sec:I-PIA}

 In this section, we will develop the I-PIA iteration method for the
    implicit curve and surface reconstruction,
    and prove its convergence.

\subsection{I-PIA for implicit curve reconstruction}
\label{subsec:I-PIA_Curve}

Given the unorganized data point set
    $\{p_i=(x_i,y_i)\}_{i=1}^{n}$~\pref{eq:Datapnt_2d},
		with a set of associated oriented unit normals $\{\bm{n}_i, i=1,2,...,n\}$,
    we want to reconstruct an implicit curve from the data point set.
To avoid the trivial solution,
    we need to add extra offset points $\{p_l = (x_l,y_l)\}_{l=n+1}^{2n}$ to the data point set~\pref{eq:Datapnt_2d}.
 The offset points are generated along the normal vector $\bm{n}$ at small distance $\sigma$ \cite{Carr2001Reconstruction,Fasshauer2007Meshfree,Rouhani2015Implicit}, i.e.,
\begin{equation*}
p_l=p_i+\sigma \bm{n}_i, \quad l=n+i,\quad i=1,2,...,n.
\end{equation*}
Let $\epsilon$ be the value of the implicit function at the offset points, i.e.,
 $$f(p_l)=\epsilon, \quad l=n+1,n+2,...,2n.$$

Define the initial implicit B-spline function as follows:
\begin{equation}\label{eq:Impcurve_0}
f^{(0)}(x,y)=\sum_{i=1}^{N_u}\sum_{j=1}^{N_v}C_{i,j}^{(0)}B_i(x)B_j(y),
\end{equation}
 where the initial control coefficients are taken as $C_{ij}^{(0)} = 0$.

%\begin{figure*}
%\begin{center}
%$
%{\it \Gamma}=diag\left( \frac{1}{{\sum_{k=1}^{m}B_1(x_k)B_1(y_k)}},...,
%  \frac{1}{\sum_{k=1}^{m}B_1(x_k)B_{N_v}(y_k)},...,\frac{1}{\sum_{k=1}^{m}B_{N_u}(x_k)B_{N_v}(y_k)} \right )$,
%
%\caption{Equation for ${\it \Gamma}$.}
%\label{Gamma:CG1}
%\end{center}
%\end{figure*}
%\begin{figure*}
%\begin{center}
%
%\begin{align*}
%B &=& \left[\begin{array}{ccccccccccccc}
%B_{1}(x_1) B_{1}(y_1)& & B_{1}(x_1)B_{2}(y_1) & & ...& & B_{1}(x_1)B_{N_v}(y_1) & & ...& & B_{N_u}(x_1)B_{N_v}(y_1)\\
%B_{1}(x_2) B_{1}(y_2)& & B_{1}(x_2)B_{2}(y_2) & & ...& & B_{1}(x_2)B_{M_v}(y_2) & & ...& & B_{N_u}(x_2)B_{N_v}(y_2)\\
%... & & ...& &  ... & &  ...& &  ...& &  ... \\
%B_{1}(x_m) B_{1}(y_m)& & B_{1}(x_m)B_{2}(y_m) & & ...& & B_{1}(x_m)B_{N_v}(y_m) & & ...& & B_{N_u}(x_m)B_{N_v}(y_m)\\
% \end{array} \right]_{m\times N_uN_v}.
%\end{align*}
%
%  \caption{Iterative matrix B.}
%\label{B:CG1}
%
%\end{center}
%\end{figure*}

Let $\delta^{(0)}_k, k=1,2,\cdots,2n$ be the difference vectors for the
    data points,
    which can be calculated as,
\begin{eqnarray}\nonumber
\delta_{k}^{(0)}&=&0-f^{(0)}(x_k,y_k), \quad k=1,2, ... ,n,\\\nonumber
\delta_{l}^{(0)}&=&\epsilon-f^{(0)}(x_l,y_l), \quad l=n+1,n+2, ... ,2n.
\end{eqnarray}
Moreover, let $\Delta^{(0)}_{ij},i=1,2, ... ,N_u, j=1,2, ... ,N_v$ be the
    difference vectors for the control coefficients, defined as,
\begin{equation*}
\Delta_{ij}^{(0)} = \mu \sum_{k=1}^{2n}B_i(x_k)B_j(y_k)\delta_{k}^{(0)},
\end{equation*}
    where, $\mu$ is a weight, and its selection method is explained in Remark~\ref{rem:weight_curve}.
The new coefficients are obtained by:
\begin{eqnarray}\nonumber
C_{ij}^{(1)}&=&C_{ij}^{(0)}+\Delta_{ij}^{(0)},\\\nonumber
\end{eqnarray}
and the new implicit B-spline curve,	
\begin{eqnarray*}
f^{(1)}(x,y)&=&\sum_{i=1}^{N_u}\sum_{j=1}^{N_v}C_{ij}^{(1)}B_i(x)B_j(y).\\
\end{eqnarray*}

Likewise, suppose we have obtained the $\alpha$-th implicit curves $f^{(\alpha)}(x,y)$ after the $\alpha$-th
iteration, and let,
\begin{eqnarray}\nonumber
\delta_{k}^{(\alpha)}&=&0-f^{(\alpha)}(x_k,y_k), \quad k=1,2,...,n,\\\nonumber
\delta_{l}^{(\alpha)}&=&\epsilon-f^{(\alpha)}(x_l,y_l), \quad l=n+1,n+2,...,2n,\end{eqnarray}
\begin{eqnarray}\nonumber
\Delta_{ij}^{(\alpha)}&=&
    \mu \sum_{k=1}^{2n}B_i(x_k)B_j(y_k)\delta_{k}^{(\alpha)},\\\nonumber
\end{eqnarray}
\begin{equation}\label{eq:CCefficient}
C_{ij}^{(\alpha+1)} = C_{ij}^{(\alpha)}+\Delta_{ij}^{(\alpha)},
 %{C}_{i,j}^{(\alpha+1)}&=&{C}_{i,j}^{(\alpha)}+{\Delta}_{i,j}^{(\alpha)},
\end{equation}
\begin{equation}\label{eq:Impcurve_alpha}
 f^{(\alpha+1)}(x,y) = \sum_{i=1}^{N_u}\sum_{j=1}^{N_v}C_{ij}^{(\alpha+1)}B_i(x)B_j(y).\\\
\end{equation}
%\begin{eqnarray}\label{eq:Eq23}
%\nonumber f^{(\alpha+1)}({p}_l)=\sum_{i=1}^{M}\sum_{j=1}^{M}{C}_{i,j}^{(\alpha+1)}B_i(x_l)B_j(y_l)\\\
% l=n+1,n+2,...,N.
%\end{eqnarray}
From the above iterative procedure,
    we generate a series of implicit
    B-spline functions $\{f^{(\alpha)}(x,y), \alpha=0,1,2,... \}$.

 Now,
    arrange the control coefficients
    $C_{ij},i=1,2,\cdots,N_u, j=1,2,\cdots,N_v$
    to a vector in the lexicographic order, i.e.,
\begin{equation*}
    \bm{C}^{(\alpha)} = [C^{(\alpha)}_{11},C^{(\alpha)}_{12},...,C^{(\alpha)}_{1N_v},
    ...,C^{(\alpha)}_{N_u,N_v}]^{\rm T},
\end{equation*}
and let,
\begin{equation*}
    \bm{b} = [b_1,b_2,\cdots,b_{2n}]^{\rm T}
           = [\underbrace{0,0,\cdots,0}_{n},
              \underbrace{\epsilon,\epsilon,\cdots,\epsilon}_{n}]^{\rm T}.
\end{equation*}
According to (\ref{eq:CCefficient}), we have
\begin{eqnarray}\nonumber
C_{ij}^{(\alpha+1)}&=&C_{ij}^{(\alpha)}+
    \mu \sum_{k=1}^{2n}B_i(x_k)B_j(y_k)\delta_{k}^{(\alpha)},\\\nonumber
        &=&C_{ij}^{(\alpha)}+
    \mu \sum_{k=1}^{2n}B_i(x_k)B_j(y_k) \left[b_k-\sum_{i=1}^{N_u}\sum_{j=1}^{N_v}C_{ij}^{(\alpha)}B_i(x_k)B_j(y_k)\right]. \\ \nonumber
\end{eqnarray}
Then, the I-PIA~\pref{eq:CCefficient} for implicit curve reconstruction
    can be represented in the matrix form,
\begin{eqnarray}\nonumber
\bm{C}^{(\alpha+1)}&=&\bm{C}^{(\alpha)}+\mu B^{\rm T}\left(\bm{b}-B\bm{C}^{(\alpha)}\right),\\\label{eq:Iterfmt}
&=& \left(I-\mu B^{\rm T}B\right)\bm{C}^{(\alpha)} + \mu B^{\rm T}\bm{b}, \quad \alpha=0,1,2,...
\end{eqnarray}
 where, $B$ is the collocation matrix of the basis functions (arranged in the lexicographical order),
 $$\{B_1(x)B_1(y), B_1(x)B_2(y), \cdots, B_1(x)B_{N_v}(y),\cdots,B_{N_u}(x)B_1(y), \cdots,B_{N_u}(x)B_{N_v}(y)\}$$
 on the point sequence $\{(x_k,y_k), k=1,2,\cdots,2n\}$,
    i.e.,
 \begin{equation} \label{eq:b_2d_collocation}
    B = \begin{bmatrix}
            B_1(x_1)B_1(y_1) & B_1(x_1)B_2(y_1) & \cdots & B_{N_u}(x_1)B_{N_v}(y_1) \\
            B_1(x_2)B_1(y_2) & B_1(x_2)B_2(y_2) & \cdots & B_{N_u}(x_2)B_{N_v}(y_2) \\
            \cdots   & \cdots   & \cdots & \cdots \\
            B_1(x_{2n})B_1(y_{2n}) & B_1(x_{2n})B_2(y_{2n}) & \cdots & B_{N_u}(x_{2n})B_{N_v}(y_{2n}) \\
          \end{bmatrix}.
 \end{equation}

%where, $B$ is the semi-Kronecker product of the matrices $B_x$ and $B_y$,
%    define as:
% \begin{equation} \label{eq:semi_kron}\textcolor[rgb]{0.00,0.07,1.00}{
%    B = B_x \bar{\otimes} B_y = \begin{bmatrix}
%            B_1(x_1) B_yR_1       & B_2(x_1) B_yR_1     & \cdots   & B_{N_u}(x_1) B_yR_1\\
%            B_1(x_2) B_yR_2       & B_2(x_2) B_yR_2     & \cdots   & B_{N_u}(x_2) B_yR_2\\
%            \cdots                & \cdots              & \cdots   & \cdots \\
%            B_1(x_{2n}) B_yR_{2n} & B_2(x_{2n}) B_yR_{2n}& \cdots  & B_{N_u}(x_{2n}) B_yR_{2n}\\
%          \end{bmatrix}.}
% \end{equation}
% \begin{equation} \label{eq:bx_by}
%    B_x = \begin{bmatrix}
%            B_1(x_1) & B_2(x_1) & \cdots & B_{N_u}(x_1) \\
%            B_1(x_2) & B_2(x_2) & \cdots & B_{N_u}(x_2) \\
%            \cdots   & \cdots   & \cdots & \cdots \\
%            B_1(x_{2n}) & B_2(x_{2n}) & \cdots & B_{N_u}(x_{2n}) \\
%          \end{bmatrix},\quad
%    B_y = \begin{bmatrix}
%            B_1(y_1) & B_2(y_1) & \cdots & B_{N_v}(y_1) \\
%            B_1(y_2) & B_2(y_2) & \cdots & B_{N_v}(y_2) \\
%            \cdots   & \cdots   & \cdots & \cdots \\
%            B_1(y_{2n}) & B_2(y_{2n}) & \cdots & B_{N_v}(y_{2n}) \\
%          \end{bmatrix}
%					\textcolor[rgb]{0.00,0.07,1.00}{
%					=\begin{bmatrix}
%            B_yR_1 \\
%            B_yR_2 \\
%            \cdots    \\
%            B_yR_{2n}\\}
%          \end{bmatrix}.
% \end{equation}

 \begin{remark} \label{rem:weight_curve}
    For the convergence of the I-PIA iterative method~\pref{eq:Iterfmt},
    the weight $\mu$~\pref{eq:Iterfmt} should satisfy
    $0 < \mu < \frac{2}{\lambda_{max}(B^T B)}$,
    where $\lambda_{max}(B^T B)$ is the largest eigenvalue of
    $B^T B$~\pref{eq:Iterfmt}.
    Moreover, for the fast convergence of I-PIA,
    a practical selection of $\mu$ in our implementation is (refer to~\cite{Deng2014Progressive}),
    \begin{equation*}
        \mu = \frac{2}{C},\ \text{where},\
        C = \norm{B^T B}_{\infty} = \max_{ij} \sum_k B_i(x_k) B_j(y_k).
    \end{equation*}
 \end{remark}

 The convergence analysis of the I-PIA for the implicit
    curve reconstruction~\pref{eq:Iterfmt} will be presented in Section~\ref{subsec:convergence}.

 %\textcolor[rgb]{0.00,0.07,1.00}{I have revised here.}
 %\vspace{0.5cm}

%------------------------------------------------------------------%surface
%------------------------------------------------------------------

\subsection{I-PIA for implicit surfaces reconstruction}
\label{subsec:I-PIA_Surface}

 The I-PIA iterative method for the implicit curve reconstruction can be
    easily extended to implicit surface reconstruction.
 In the following, the details of I-PIA for implicit surface reconstruction
    will be presented.

  Given the unorganized data point set $\{p_i=(x_i,y_i,z_i)\}_{i=1}^{n}$     \pref{eq:Datapnt_3d}, with a set of associated oriented unit normals
	    $\{\bm{n}_i, i=1,2,...,n\}$, we want to reconstruct an implicit
			surface from the data point set.
 To avoid the trivial solution,
    extra offset point $\{p_l=(x_l,y_l,z_l)\}_{l=n+1}^{2n}$ should be added to the data point set~\pref{eq:Datapnt_3d}.
 The offset points are generated along the normal vector $\bm{n}$ at small
    distance $\sigma$  \cite{Carr2001Reconstruction,Fasshauer2007Meshfree,Rouhani2015Implicit},
    i.e.,
	\begin{equation*}
        p_l=p_i+\sigma \bm{n}_i, \quad l=n+i,\quad i=1,2,...,n.
    \end{equation*}
 Moreover, let $\epsilon$ be the value of the implicit function at the offset points, i.e.,
  $$f(p_l)=\epsilon, \quad l=n+1,n+2,...,2n.$$

   Define the initial B-spline implicit function as follows:
\begin{equation*}\label{eq:Imp-surf_0}
f^{(0)}(x,y,z) = \sum_{i=1}^{N_u}\sum_{j=1}^{N_v}\sum_{k=1}^{N_w}
                    C_{ijk}^{(0)}B_i(x)B_j(y)B_k(z),
\end{equation*}
    where $C_{ijk}^{(0)}=0, i=0,1,\cdots,N_u, j=0,1,\cdots,N_v, k=0,1,\cdots,N_w$.
Let $\delta^{(0)}_r, r=1,2,...,2n$ be the difference vectors for the data points, calculated by,
\begin{eqnarray*}
\delta_{r}^{(0)}&=&0-f^{(0)}(x_r,y_r,z_r), \quad r=1,2,\cdots,n,\\
\delta_{l}^{(0)}&=&\epsilon-f^{(0)}(x_l,y_l,z_l), \quad l=n+1,n+2,\cdots,2n.
\end{eqnarray*}
 Then, the difference vectors for the control coefficients
    $\Delta_{ijk}^{(0)},i=1,2, ... ,N_u, j=1,2, ... ,N_v, k=1,2, ... ,N_w,$
    can be constructed as,
\begin{equation*}
\Delta_{ijk}^{(0)} = \mu
    \sum_{r=1}^{2n}B_i(x_r)B_j(y_r)B_k(z_r)\delta_{r}^{(0)},
\end{equation*}
 where, $\mu$ is a weight,
    and its selection method is explained in Remark~\ref{rem:weight_surface}.
 Similar as the curve case,
    we can get the new coefficients,
\begin{equation*}
    C_{ijk}^{(1)} = C_{ijk}^{(0)}+\Delta_{ijk}^{(0)},
\end{equation*}
and the new implicit B-spline surface,
\begin{equation}\label{eq:Imp-surf}
    f^{(1)}(x,y,z)=\sum_{i=1}^{N_u}\sum_{j=1}^{N_v}\sum_{k=1}^{N_w}C_{ijk}^{(1)}B_i(x)B_j(y)B_k(z).
\end{equation}

Likewise, suppose we have obtained the $\alpha$-th implicit B-spline surface
    $f^{(\alpha)}(x,y,z)$ after the $\alpha$-th iteration,
    and let,
\begin{eqnarray*}
\delta_{r}^{(\alpha)}&=&0-f^{(\alpha)}(x_r,y_r,z_r), \quad r=1,2,\cdots,n,\\
\delta_{l}^{(\alpha)}&=&\epsilon-f^{(\alpha)}(x_r,y_r,z_r), \quad l=n+1,n+2,\cdots,2n,
\end{eqnarray*}
\begin{equation*}
\Delta_{ijk}^{(\alpha)} =
            \mu \sum_{r=1}^{2n}B_i(x_r)B_j(y_r)B_k(z_r)\delta_{r}^{(\alpha)},
\end{equation*}
\begin{equation} \label{eq:i-pia_surface}
    C_{ijk}^{(\alpha+1)} = C_{ijk}^{(\alpha)}+\Delta_{ijk}^{(\alpha)},
\end{equation}
\begin{equation}\label{eq:Imp-surf_alpha}
f^{(\alpha+1)}(x,y,z) = \sum_{i=1}^{N_u}\sum_{j=1}^{N_v}\sum_{k=1}^{N_w}
        C_{ijk}^{(\alpha+1)}B_i(x)B_j(y)B_k(z).%k=1,2,...,n
 %l=n+1,n+2,...,N.
\end{equation}
 In this way, a series of implicit B-spline surface $\{f^{(\alpha)}(x,y,z), \alpha=0,1,2,... \}$ is generated.

 Similar as the curve case in Section~\ref{subsec:I-PIA_Curve},
    arranging the control coefficients
    $$C_{ijk},i=1,2,\cdots,N_u, j=1,2,\cdots,N_v, k=1,2,\cdots,N_w,$$
    to a vector in the lexicographic order, i.e.,
\begin{equation*}
    \bm{C}^{(\alpha)} = [C^{(\alpha)}_{111},C^{(\alpha)}_{112},...,C^{(\alpha)}_{1,1,N_w},
    ...,C^{(\alpha)}_{N_u,N_v,N_w}]^{\rm T},
\end{equation*}
and letting,
\begin{equation*}
    \bm{b} = [b_1,b_2,\cdots,b_{2n}]^{\rm T}
           = [\underbrace{0,0,\cdots,0}_{n},
              \underbrace{\epsilon,\epsilon,\cdots,\epsilon}_{n}]^{\rm T}.
\end{equation*}
 The I-PIA~\pref{eq:i-pia_surface} for implicit surface reconstruction
    can be represented in the matrix form,
\begin{eqnarray}\nonumber
\bm{C}^{(\alpha+1)}&=&\bm{C}^{(\alpha)}+\mu B^{\rm T}\left(\bm{b}-B\bm{C}^{(\alpha)}\right),\\\label{eq:patch_it_fmt}
&=& \left(I-\mu B^{\rm T}B\right)\bm{C}^{(\alpha)} + \mu B^{\rm T}\bm{b}, \quad \alpha=0,1,2,...
\end{eqnarray}
    where, $B$ is the collocation matrix of the basis functions
    (arranged in the lexicographical order),
 \begin{equation*}
 \begin{split}
 \{B_1(x)B_1(y)B_1(z), B_1(x)B_1(y)B_2(z), \cdots,& B_1(x)B_1(y)B_{N_v}(z),\cdots,\\
 & B_{N_u}(x)B_{N_v}(y)B_1(z), \cdots,B_{N_u}(x)B_{N_v}(y)B_{N_w}(z)\}
 \end{split}
 \end{equation*}
 on the point sequence $\{(x_k,y_k,z_k), k=1,2,\cdots,2n\}$,
    i.e.,
 \begin{equation} \label{eq:b_3d_collocation}
    B = \begin{bmatrix}
            B_1(x_1)B_1(y_1)B_1(z_1) & B_1(x_1)B_1(y_1)B_{N_w}(z_1) & \cdots & B_{N_u}(x_1)B_{N_v}(y_1)B_{N_w}(z_1) \\
            B_1(x_2)B_1(y_2)B_1(z_2) & B_1(x_2)B_1(y_2)B_2(z_2) & \cdots & B_{N_u}(x_2)B_{N_v}(y_2)B_{N_w}(z_2) \\
            \cdots   & \cdots   & \cdots & \cdots \\
            B_1(x_{2n})B_1(y_{2n})B_1(z_{2n}) & B_1(x_{2n})B_1(y_{2n})B_2(z_{2n}) & \cdots & B_{N_u}(x_{2n})B_{N_v}(y_{2n})B_{N_w}(z_{2n}) \\
          \end{bmatrix}.
 \end{equation}

 \begin{remark} \label{rem:weight_surface}
    Similar as Remark~\ref{rem:weight_curve},
    for the convergence of the I-PIA iterative method~\pref{eq:patch_it_fmt},
    the weight $\mu$~\pref{eq:patch_it_fmt} should satisfy
    $0 < \mu < \frac{2}{\lambda_{max}(B^T B)}$,
    where $\lambda_{max}(B^T B)$ is the largest eigenvalue of
    $B^T B$~\pref{eq:patch_it_fmt}.
    Moreover,
    a practical selection of $\mu$ for the fast convergence of I-PIA is (refer to~\cite{Deng2014Progressive}),
    \begin{equation*}
        \mu = \frac{2}{C},\ \text{where},\
        C = \norm{B^T B}_{\infty} = \max_{ijk} \sum_l B_i(x_l) B_j(y_l) B_k(z_l).
    \end{equation*}
 \end{remark}

 The convergence analysis of the I-PIA for the implicit
     surface reconstruction~\pref{eq:patch_it_fmt} will be presented in Section~\ref{subsec:convergence}.

%-------------------------------------------------------------------------------
% Section:
%-------------------------------------------------------------------------------
 \subsection{Convergence analysis}
 \label{subsec:convergence}

 As shown above, the I-PIA iteration method for implicit
    curve reconstruction~\pref{eq:Iterfmt} and that for implicit surface reconstruction~\pref{eq:patch_it_fmt} can be represented in a unified form, i.e.,
 \begin{equation} \label{eq:unify_it_fmt}
    \bm{C}^{(\alpha+1)} =
    \left(I-\mu B^{\rm T}B\right)\bm{C}^{(\alpha)} + \mu B^{\rm T}\bm{b}, \quad \alpha=0,1,2,...
 \end{equation}
 The matrices $\mu B^T B$ in Eqs.~\pref{eq:Iterfmt}
    and~\pref{eq:patch_it_fmt} hold some common properties:
 \begin{enumerate}
    \item[Property 1:] The matrix $B^T B$ is positive semi-definite.
        So, its eigenvalues are nonnegative real numbers.
    \item[Property 2:] The matrix $B^T B$ is singular.
        In the implicit B-spline curve and surface reconstruction,
            the number of data points (suppose it is $n$) is less than that of the control coefficients of implicit curve and surface (suppose $m$), i.e, $m > n$.
        Because the order of the matrix $B$ is $n \times m$,
            its rank is at most $n$.
        So the rank of the $m \times m$ matrix $B^T B$ is at most $n$
            ($n < m$).
        It means that the matrix $B^T B$ is singular.
    \item[Property 3:] The eigenvalues of $\mu B^T B$ satisfy
        $0 \leq \lambda(\mu B^T B) < 2$.
        This is because the choice of the weight $\mu$,
        as well as Properties 1 and 2.
 \end{enumerate}

 \begin{remark} \label{rem:singular_mtx}
    According to Property 2, the $m \times m$ matrix $B^T B$ is singular.
    Suppose the dimension of its zero eigenspace is $m_0$.
    Then, the rank of the matrices $B^T B$ and $\mu B^T B$ is,
    $$ rank(\mu B^T B) = rank(B^T B) = m - m_0. $$
 \end{remark}

 On the other hand, the implicit curve and surface reconstruction problem
    can be formulated as solving the following least-squares fitting problem,
    \begin{equation} \label{eq:lsq}
        B^T B \bm{X} = B^T \bm{b},
    \end{equation}
    where $\bm{X}$ is an unknown vector, $B$ and $\bm{b}$ are the same as in Eq.~\pref{eq:unify_it_fmt}.
 As pointed out by Property 2, the coefficient matrix $B^T B$ is singular.
 So, if the solution of the least-square fitting system~\pref{eq:lsq} exists,
    it has infinite solutions,
    which usually leads to extra zero-level sets.
 Therefore, to eliminate the extra zero-level sets and get desirable results,
    it is required to solve the constrained minimization problem,
    \begin{equation} \label{eq:cons_mini}
        \begin{split}
        & \min_{\bm{X}} \norm{\bm{X}}_E \\
        s.t. \quad & B^T B \bm{X} = B^T \bm{b},\\
        \end{split}
    \end{equation}
    where, $\norm{\cdot}_E$ is the Euclidean norm,
    or its variants.
 In the following, we will show that,
    the I-PIA iterative method~\pref{eq:unify_it_fmt} converges to the solution of the constrained minimization problem~\pref{eq:cons_mini}.

 \begin{theorem} \label{thm:convergence}
    When the initial values $C^{(0)}$ equal $0$,
    the I-PIA iterative method~\pref{eq:unify_it_fmt} converges to the
    solution of the constrained minimization problem~\pref{eq:cons_mini},
    i.e., $(B^T B)^+ B^T \bm{b}$,
    where $(B^T B)^+$ is the Moore-Penrose (M-P) pseudo-inverse of the matrix $B^T B$.
 \end{theorem}

 \textbf{Proof:}
 By Remark~\ref{rem:singular_mtx}, $rank(\mu B^T B) = rank(B^T B) = m-m_0$.
 Because the matrix $B^T B$ is both a normal matrix and a positive
     semi-definite matrix,
     its eigen decomposition and singular value decomposition are the same,
     \begin{equation}\label{eq:svd}
        B^T B = V diag(\lambda_1,\lambda_2,...,\lambda_{m-m_0},\underbrace{0,0,...,0}_{m_0})V^T,
    \end{equation}
    where $V$ is an orthogonal matrix,
    and $\lambda_i>0, i=1,2,\cdots,m-m_0$ are both the eigen values and singular values of the matrix $B^T B$.
 Then, the M-P inverse of $B^T B$ is,
 \begin{equation*}
        (B^T B)^+ = V diag \left(\frac{1}{\lambda_1},\frac{1}{\lambda_2},...,\frac{1}{\lambda_{m-m_0}},
            \underbrace{0,0,...,0}_{m_0} \right)V^T.
 \end{equation*}
 Therefore, we have,
 \begin{equation*}
        (B^T B)^+ (B^T B) = V diag(\underbrace{1,1,\cdots,1}_{m-m_0},\underbrace{0,0,...,0}_{m_0})V^T.
    \end{equation*}

 Due to Eq.~\pref{eq:svd}, it holds,
 \begin{equation*}
        \mu B^T B = V diag(\mu \lambda_1, \mu \lambda_2,...,\mu \lambda_{m-m_0},\underbrace{0,0,...,0}_{m_0})V^T,
 \end{equation*}
 where $\mu \lambda_i, i=1,2,\cdots, m-m_0$ are the eigen values of
    the matrix $\mu B^T B$.
 Based on Property 3, they satisfy $0 \leq \mu \lambda_i <2, i=1,2,\cdots, m-m_0$. Therefore,
 \begin{equation}\label{eq:orth-mat}
    \begin{split}
    \lim_{\alpha \rightarrow \infty}\left(I-\mu B^{\rm T}B\right)^{\alpha}
    & = \lim_{\alpha \rightarrow \infty} V
        diag((1-\mu \lambda_1)^{\alpha}, (1-\mu \lambda_2)^{\alpha}, \cdots, (1-\mu \lambda_{m-m_0})^{\alpha}, \underbrace{1,1,\cdots,1}_{m_0}) V^T \\
    & = V diag(\underbrace{0,...,0}_{m-m_0},\underbrace{1,...,1}_{m_0})V^T \\
    & = I - V diag(\underbrace{1,...,1}_{m-m_0},\underbrace{0,...,0}_{m_0})V^T \\
    & = I-VV^T (B^T B)^+ (B^T B)VV^T\\
    & = I - (B^T B)^+ (B^T B).
    \end{split}
 \end{equation}

 Note that the linear system $B^T B \bm{X} = B^T \bm{b}$ has solutions,
    if and only if~\cite{james1978generalised},
    \begin{equation} \label{eq:solution_exist_cond}
        (B^T B) (B^T B)^+ (B^T \bm{b}) = B^T \bm{b}.
    \end{equation}
 Then, subtracting $(B^T B)^+ B^T \bm{b}$ from both sides
    of Eq.~\pref{eq:unify_it_fmt}, leads to,
    \begin{equation*}
        \begin{split}
            C^{(\alpha+1)} - (B^T B)^+ B^T \bm{b} & =
    (I - \mu B^T B) C^{(\alpha)} + \mu B^T \bm{b} - (B^T B)^+ B^T \bm{b}\\
        & = (I - \mu B^T B)C^{(\alpha)} + \mu (B^T B) (B^T B)^+(B^T \bm{b})
            - (B^T B)^+ B^T \bm{b} \\
        &=(I - \mu B^T B) C^{(\alpha)}-(I - \mu B^T B) (B^TB)^+B^T \bm{b} \\
        &=(I - \mu B^T B) (C^{(\alpha)} - (B^TB)^+ B^T \bm{b}) \\
        &= (I - \mu B^T B)^{\alpha+1} (C^{(0)} - (B^TB)^+B^T\bm{b}).
        \end{split}
    \end{equation*}
 So, together with Eqs.~\pref{eq:orth-mat}
    and~\pref{eq:solution_exist_cond}, we have,
 \begin{equation*}
    \begin{split}
        C^{(\infty)}- (B^T B)^+ B^T \bm{b} & = \lim_{{\alpha} \rightarrow \infty} (I-\mu B^T B)^{(\alpha+1)} (C^{(0)}- (B^T B)^+ B^T \bm{b})\\
    & = (I-(B^T B)^+ (B^T B)) (C^{(0)}-(B^T B)^+B^T \bm{b}) \\
    & = (I-(B^T B)^+(B^T B)) C^{(0)}.
    \end{split}
 \end{equation*}
 Therefore,
    \begin{equation*}
        C^{(\infty)} = (B^T B)^+ B^T \bm{b} + (I-(B^T B)^+(B^T B)) C^{(0)},
    \end{equation*}
 which are the solutions of the singular linear system $B^T B X = B^T \bm{b}$ (note that $C^{(0)}$ can take an arbitrary value).
 Among them, $(B^T B)^+ B^T \bm{b}$ is the one with minimum Euclidean norm.
 So, when the initial value $C^{(0)} = 0$,
    the I-PIA iterative format converges to
    $C^{\infty} = (B^T B)^+ B^T \bm{b}$,
    the solution of the singular linear system $B^T B X = B^T \bm{b}$
    with the minimum Euclidean norm.
 It is the solution of the constrained minimization
    problem~\pref{eq:cons_mini}.
 $\Box$

\begin{figure*}[!htb]
\centering
 \subfigure[Iteration 1.]
 {
    \includegraphics[width=0.22\linewidth]{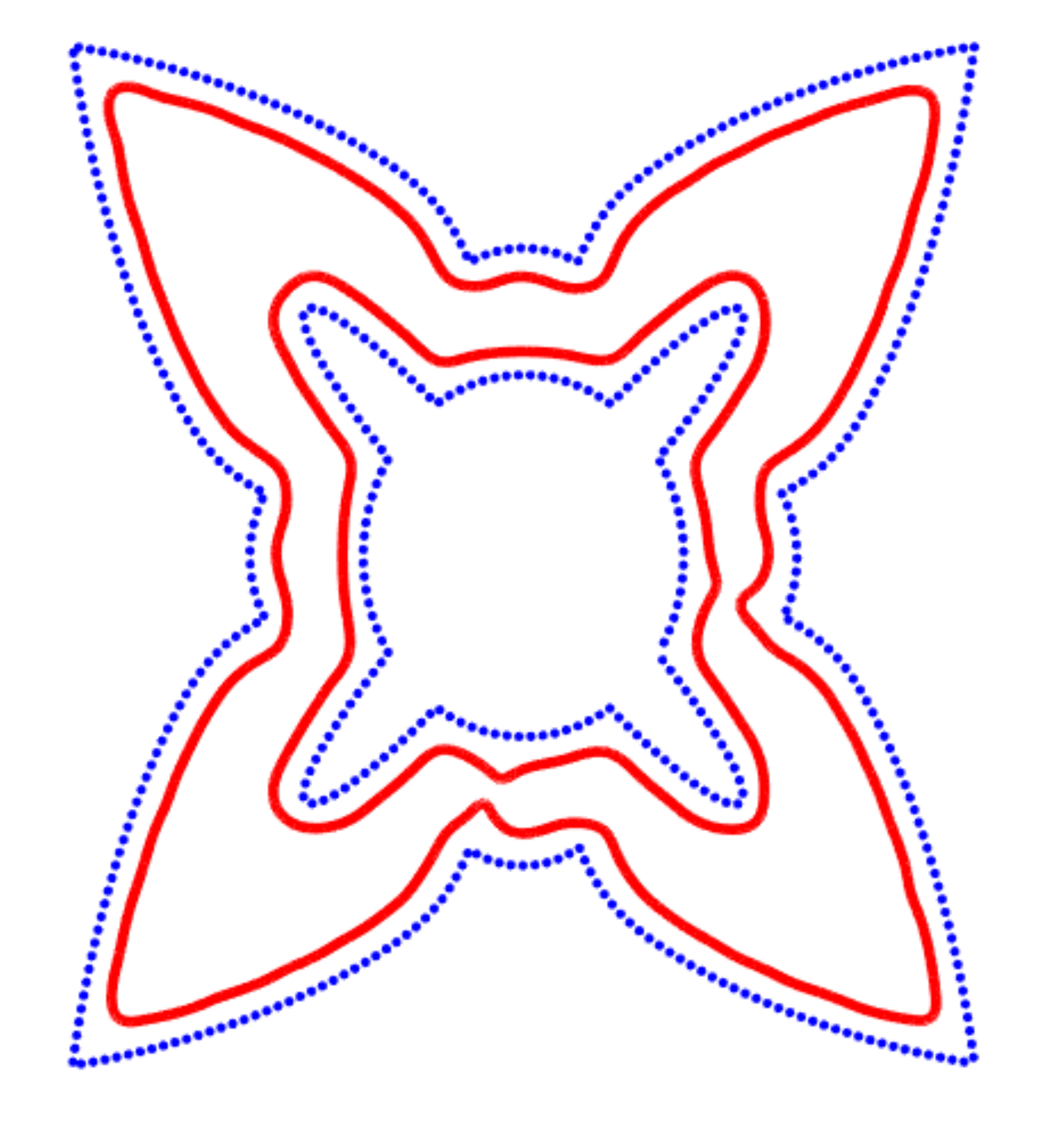}
 }
 \subfigure[Iteration 5.]
 {
    \includegraphics[width=0.22\linewidth]{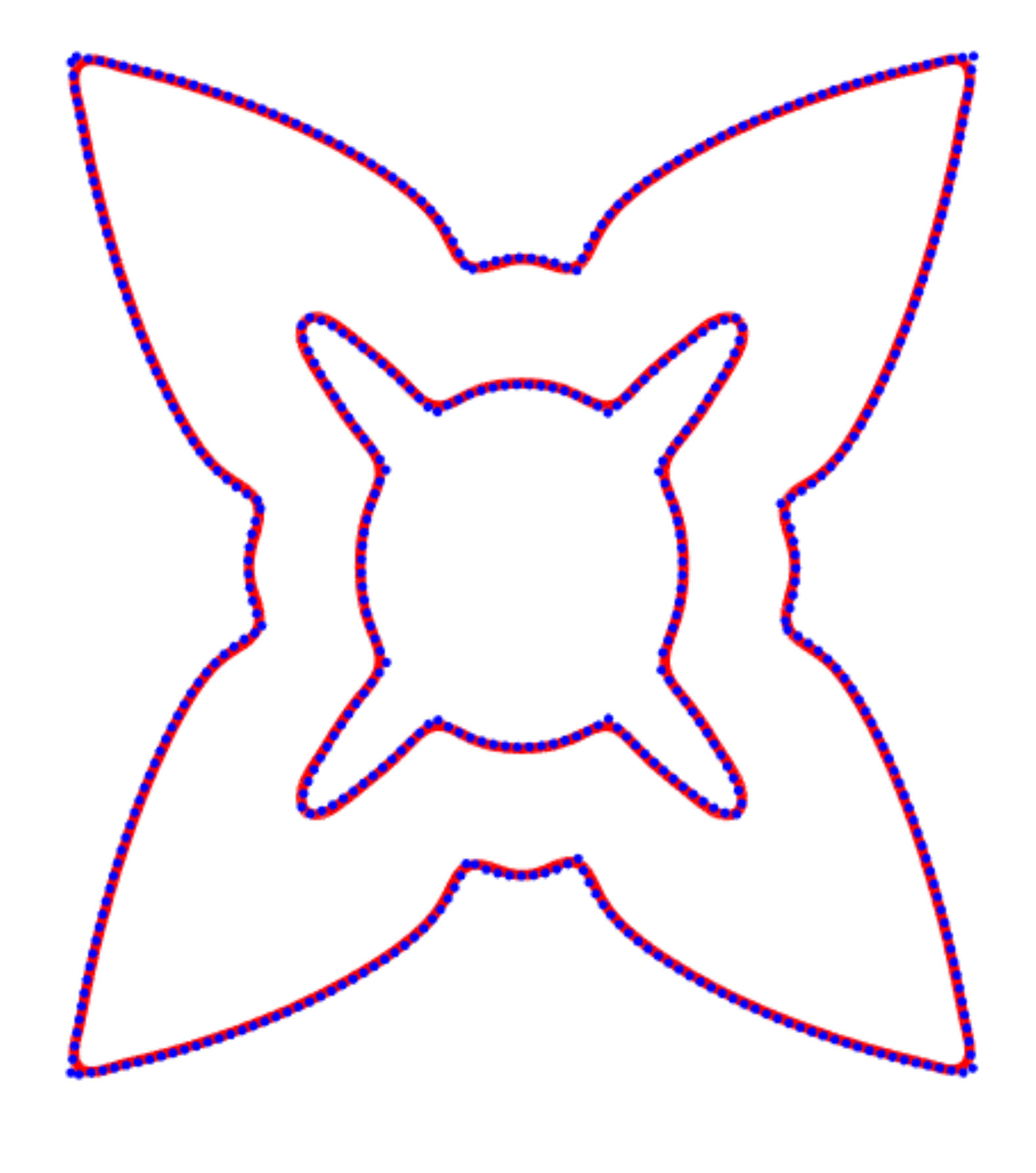}
 }
 \subfigure[Iteration 10.]
 {
    \includegraphics[width=0.22\linewidth]{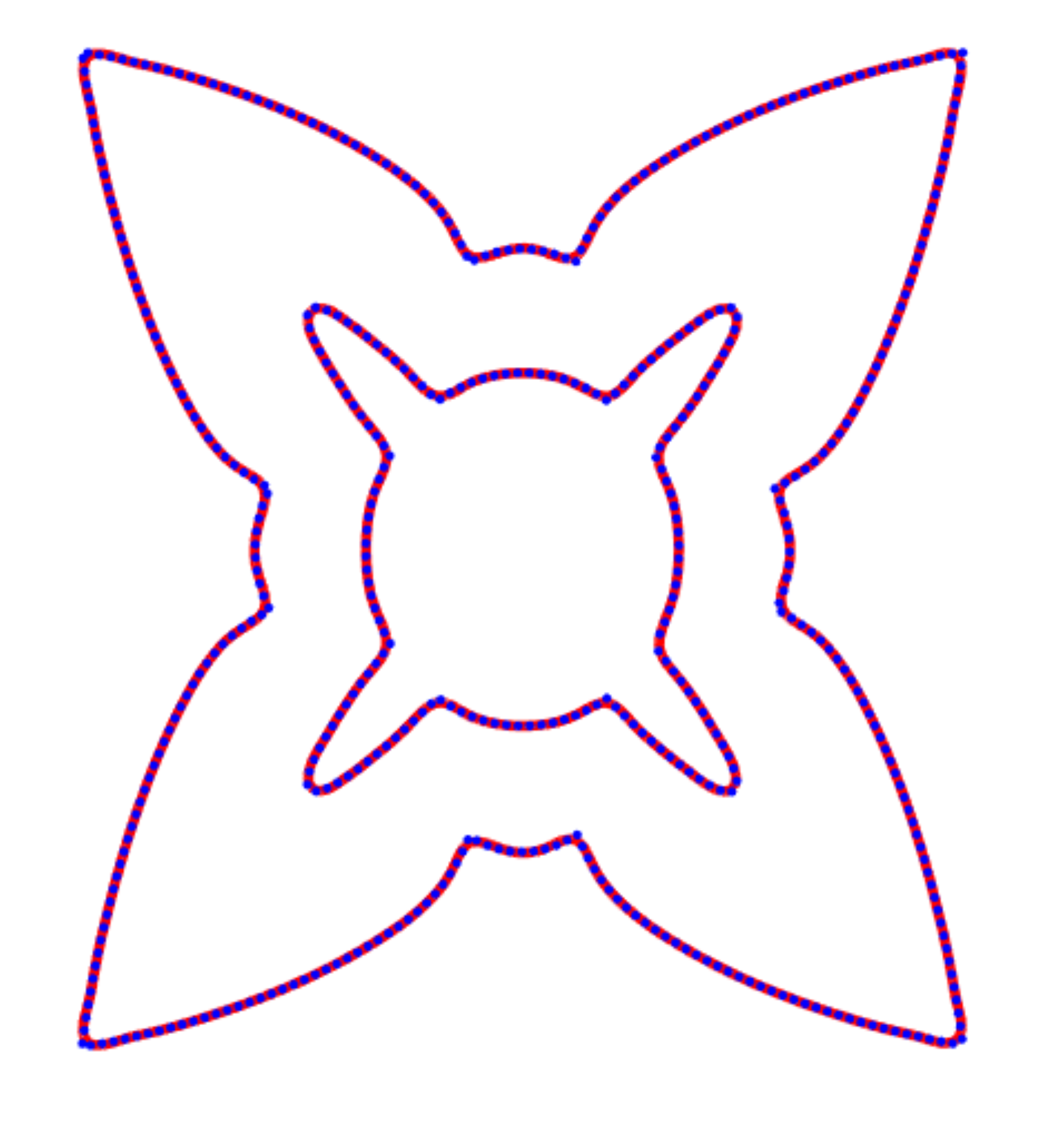}
 }
 \subfigure[Iteration 15.]
 {
    \includegraphics[width=0.22\linewidth]{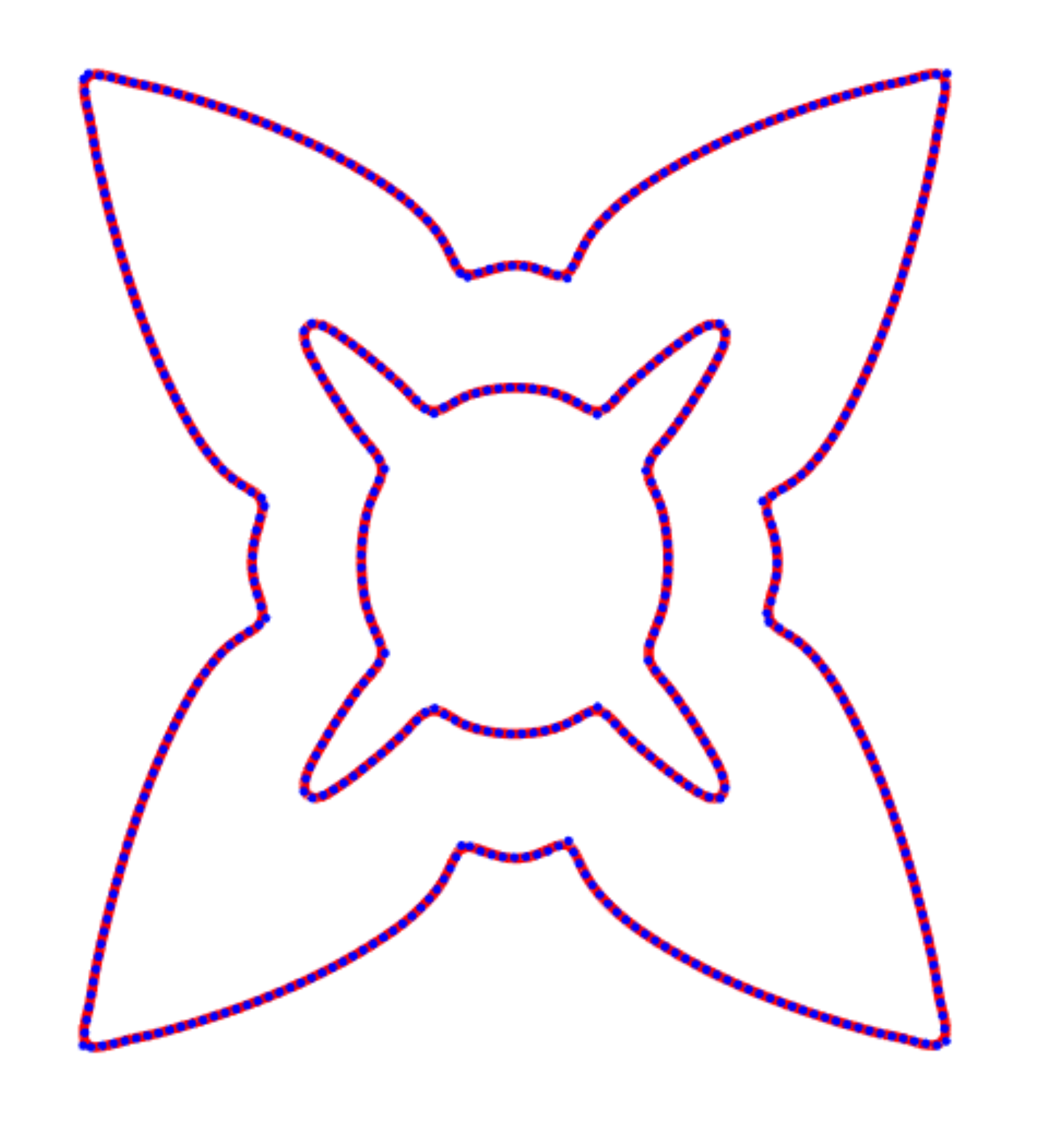}
 }
 \subfigure[Iteration 1.]
 {
    \includegraphics[width=0.22\linewidth]{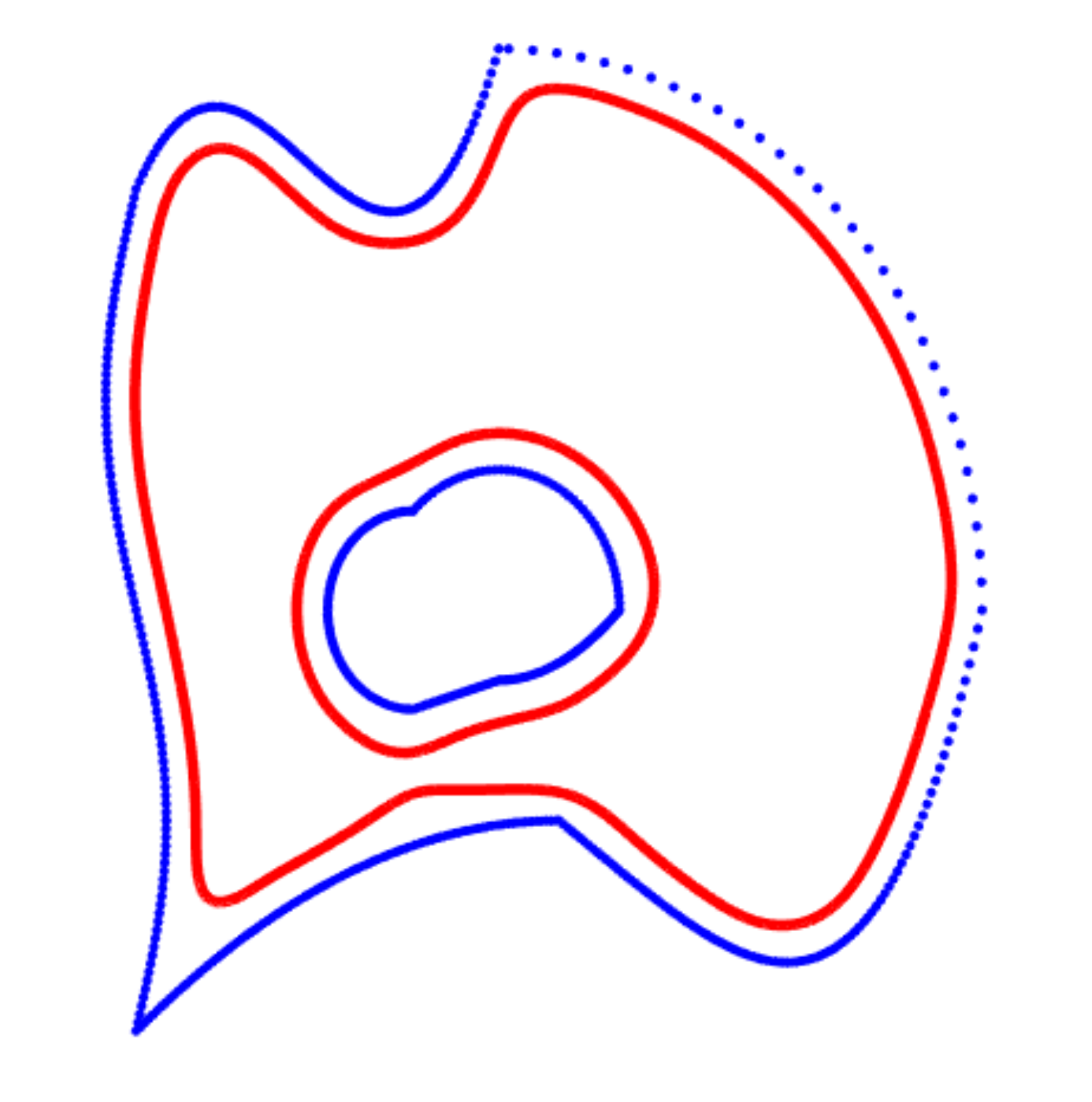}
 }
 \subfigure[Iteration 5.]
 {
    \includegraphics[width=0.22\linewidth]{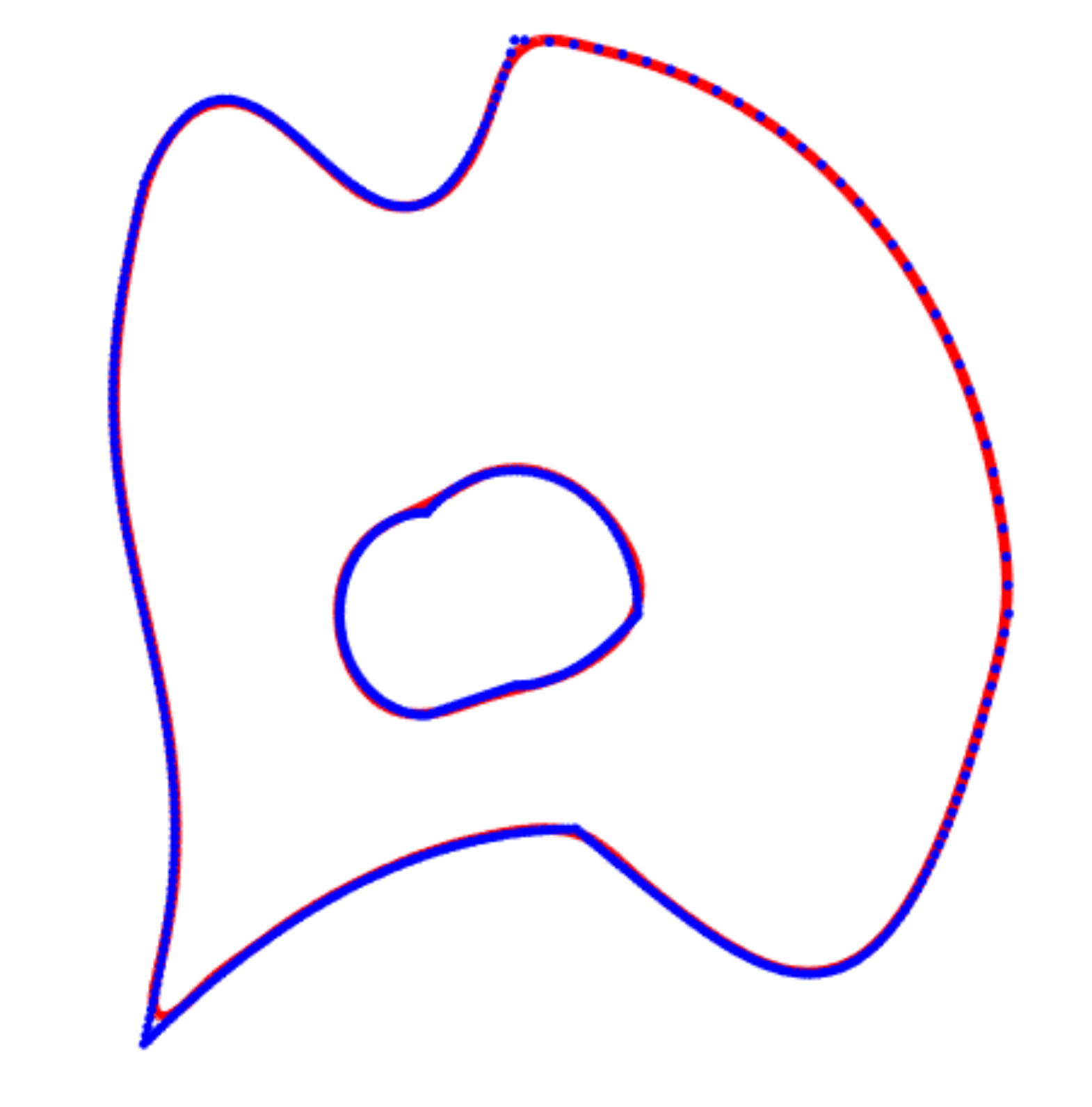}
 }
 \subfigure[Iteration 10.]
 {
    \includegraphics[width=0.22\linewidth]{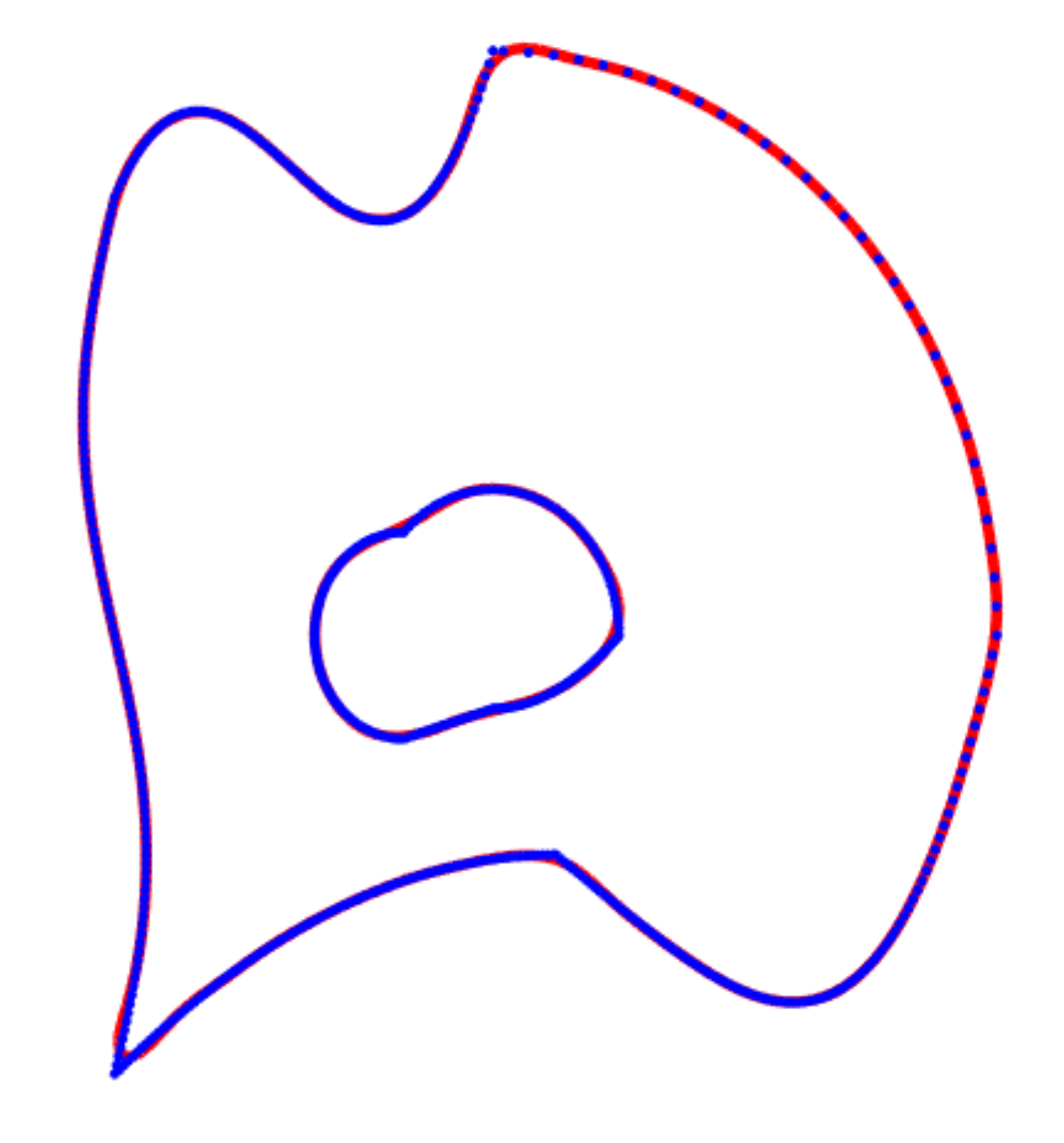}
 }
 \subfigure[Iteration 15.]
 {
    \includegraphics[width=0.22\linewidth]{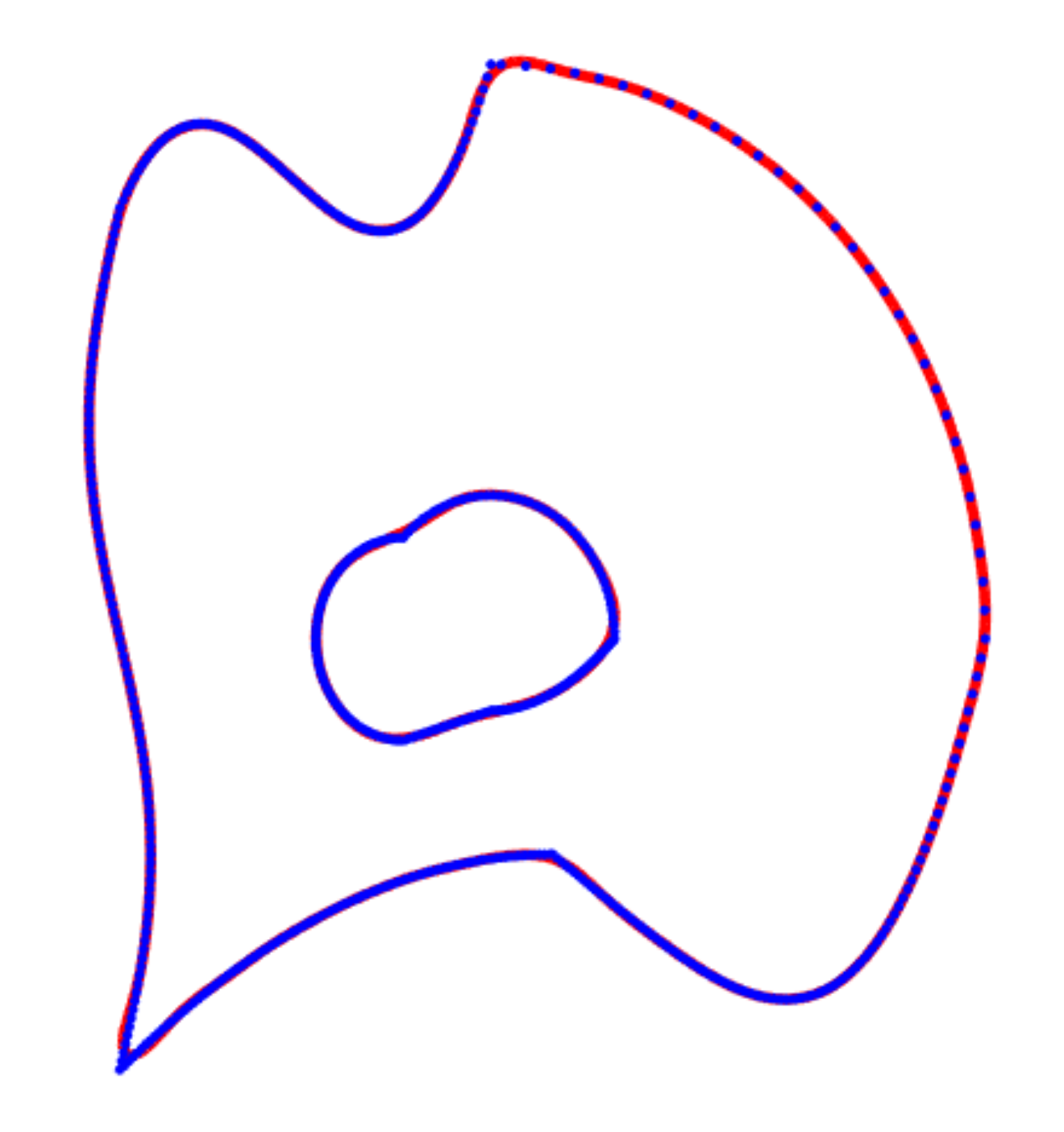}
 }
\caption
{
    Iterations in the reconstruction of $2D$ data sets: First row, \emph{flower} model, and second row, \emph{Coons curve} model.
    Blue points are the given data sets, and the red line is the reconstructed curve.
		From left to right: the 1st, 5th, 10th, 15th iteration steps.
		%In the iteration 15th, the implicit curve almost coincides with the given data sets.
 }
\label{effective-curve}
\end{figure*}

\begin{figure*}[!htb]
\centering
 \subfigure[Iteration 1.]
 {
    \includegraphics[width=0.21\linewidth]{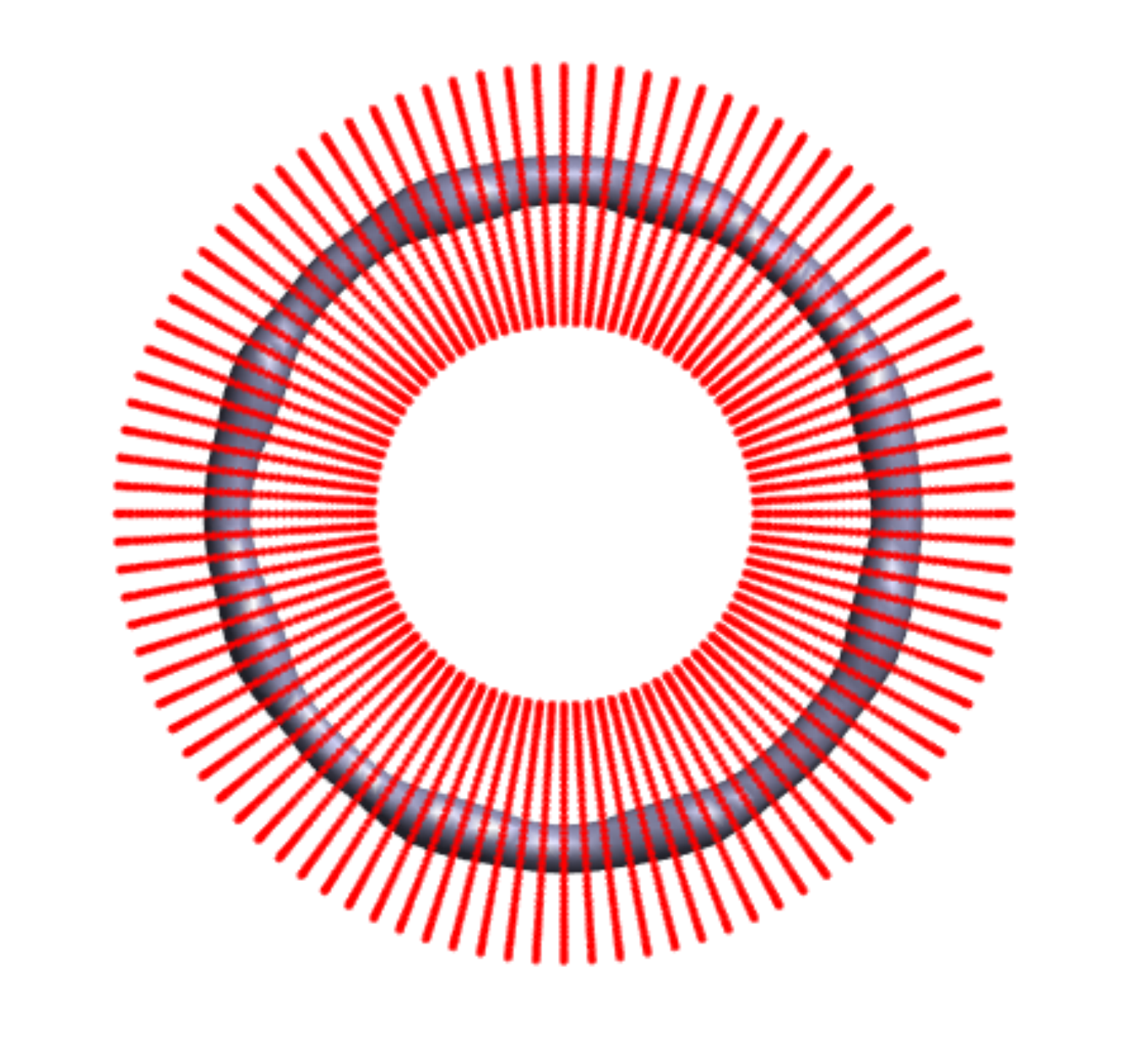}
 }
 \subfigure[Iteration 5.]
 {
    \includegraphics[width=0.21\linewidth]{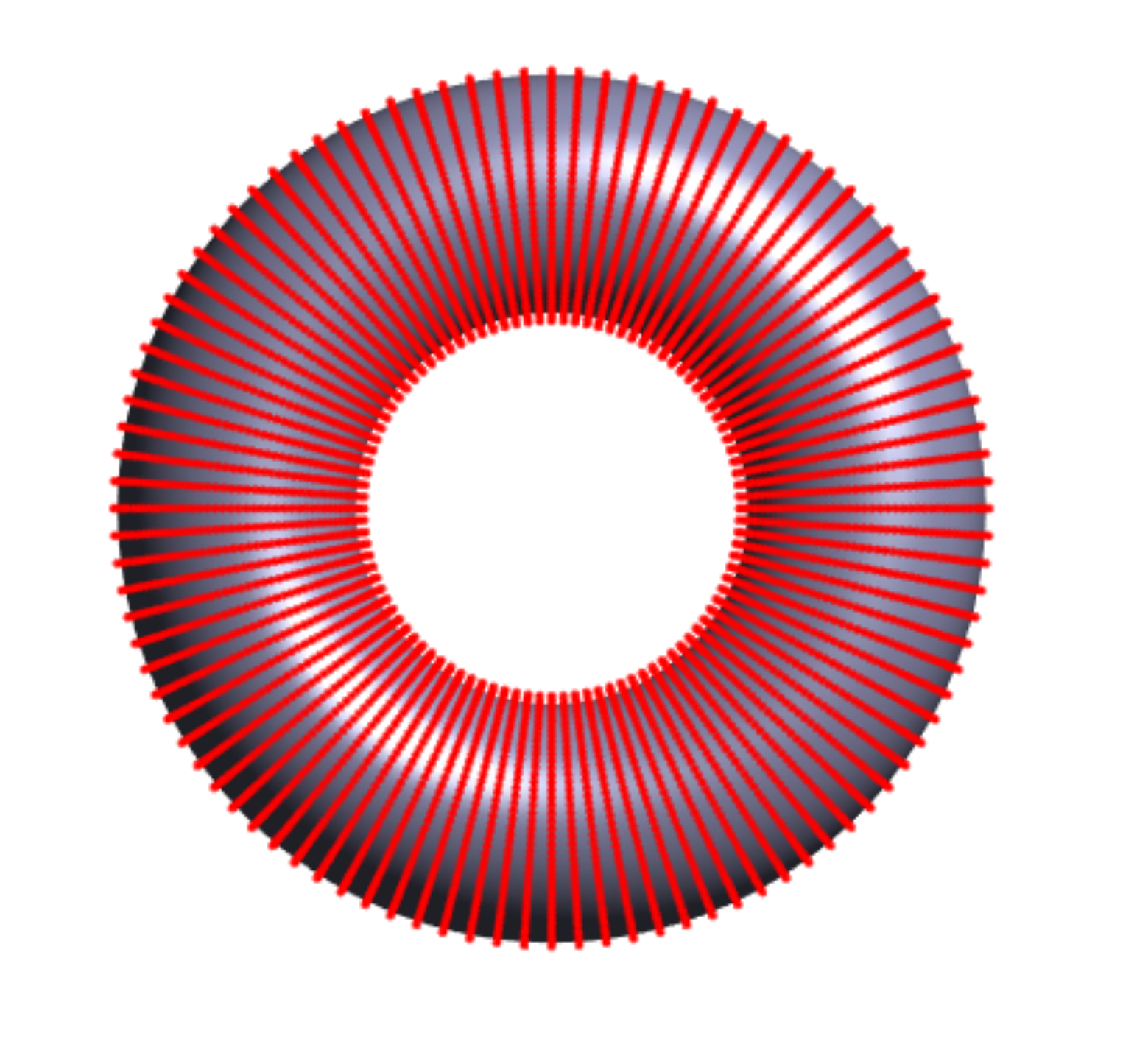}
 }
 \subfigure[Iteration 10.]
 {
    \includegraphics[width=0.20\linewidth]{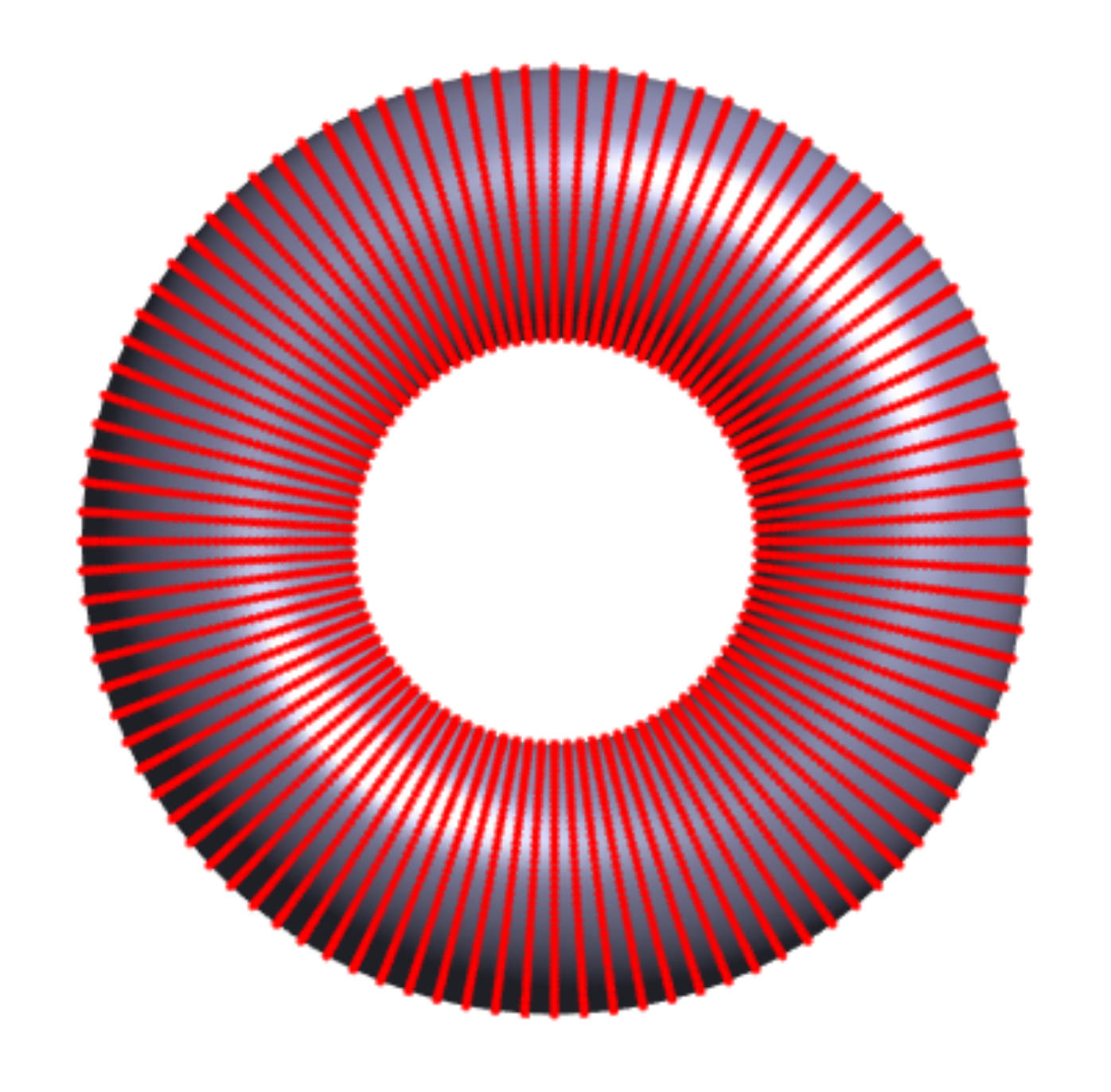}
 }
 \subfigure[Iteration 15.]
 {
    \includegraphics[width=0.21\linewidth]{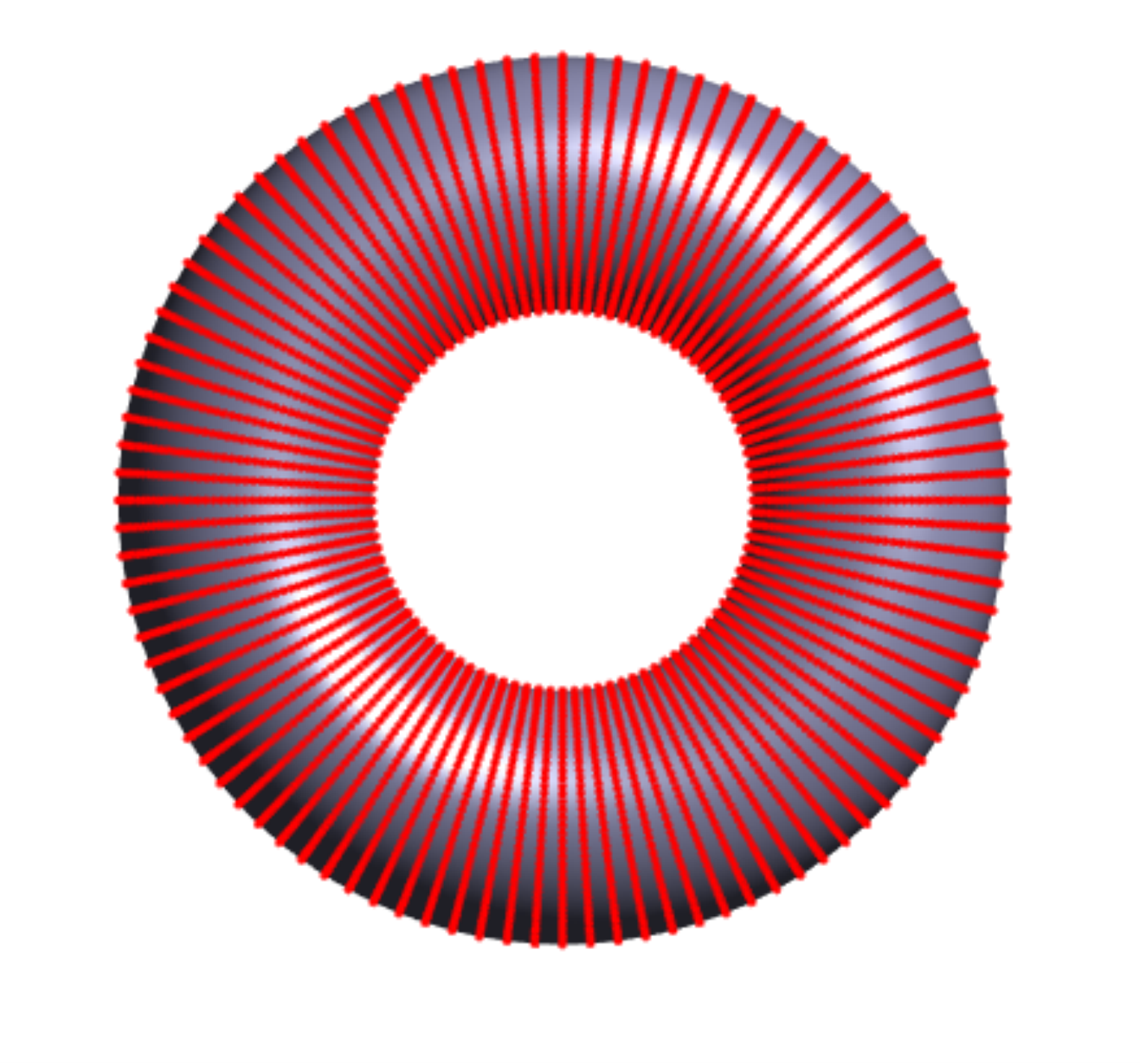}
 }
 \subfigure[Iteration 1.]
 {
    \includegraphics[width=0.21\linewidth]{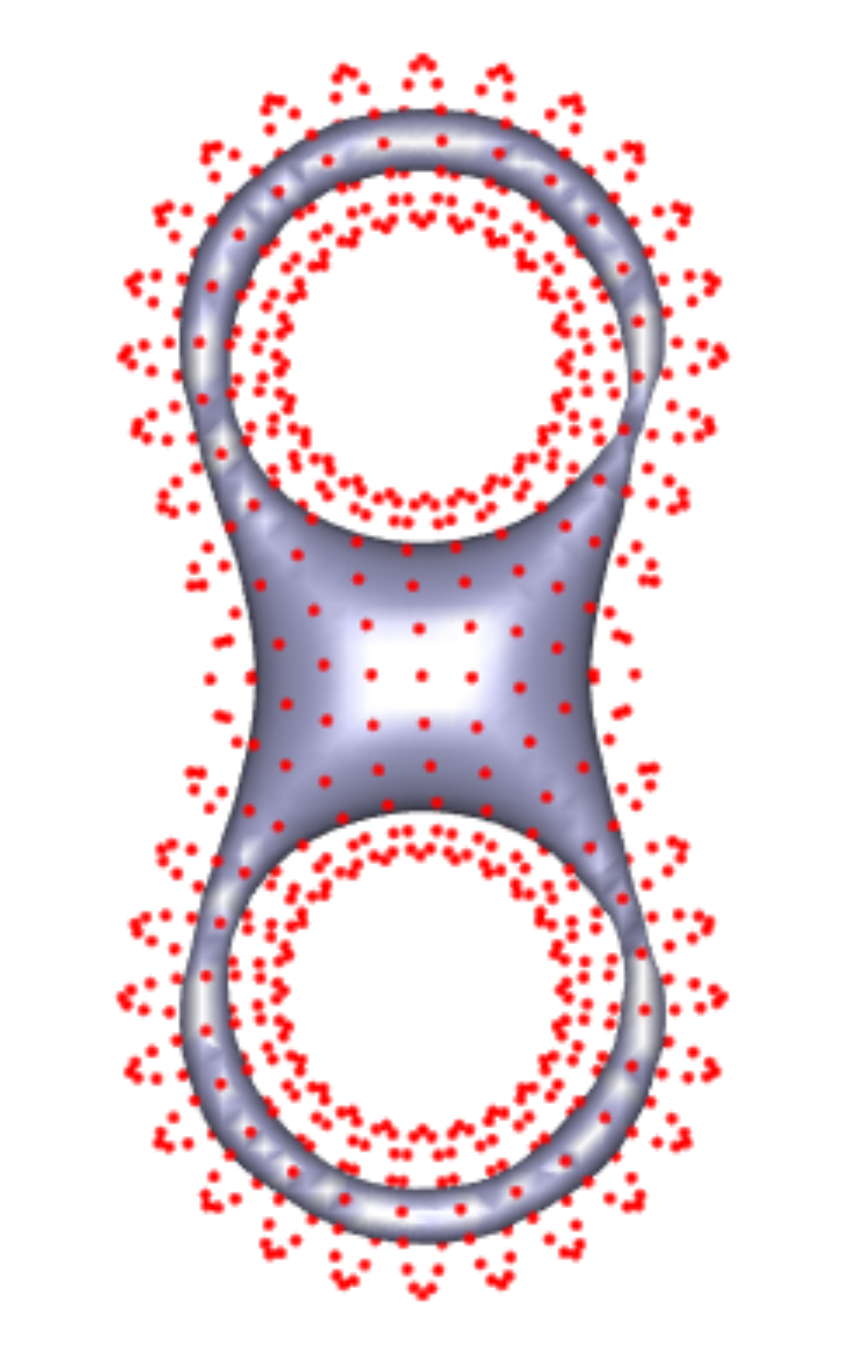}
 }
 \subfigure[Iteration 5.]
 {
    \includegraphics[width=0.21\linewidth]{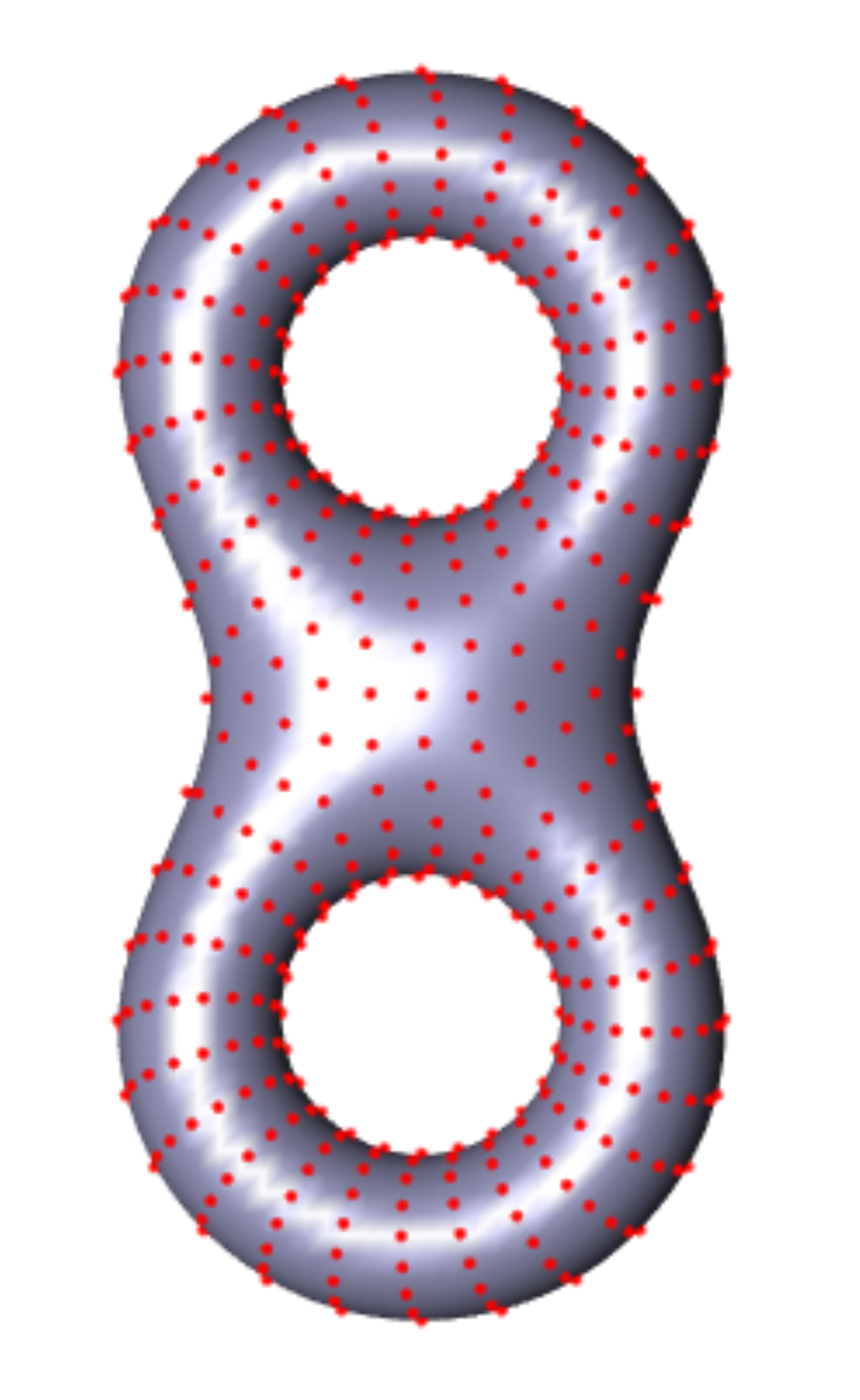}
 }
 \subfigure[Iteration 10.]
 {
    \includegraphics[width=0.21\linewidth]{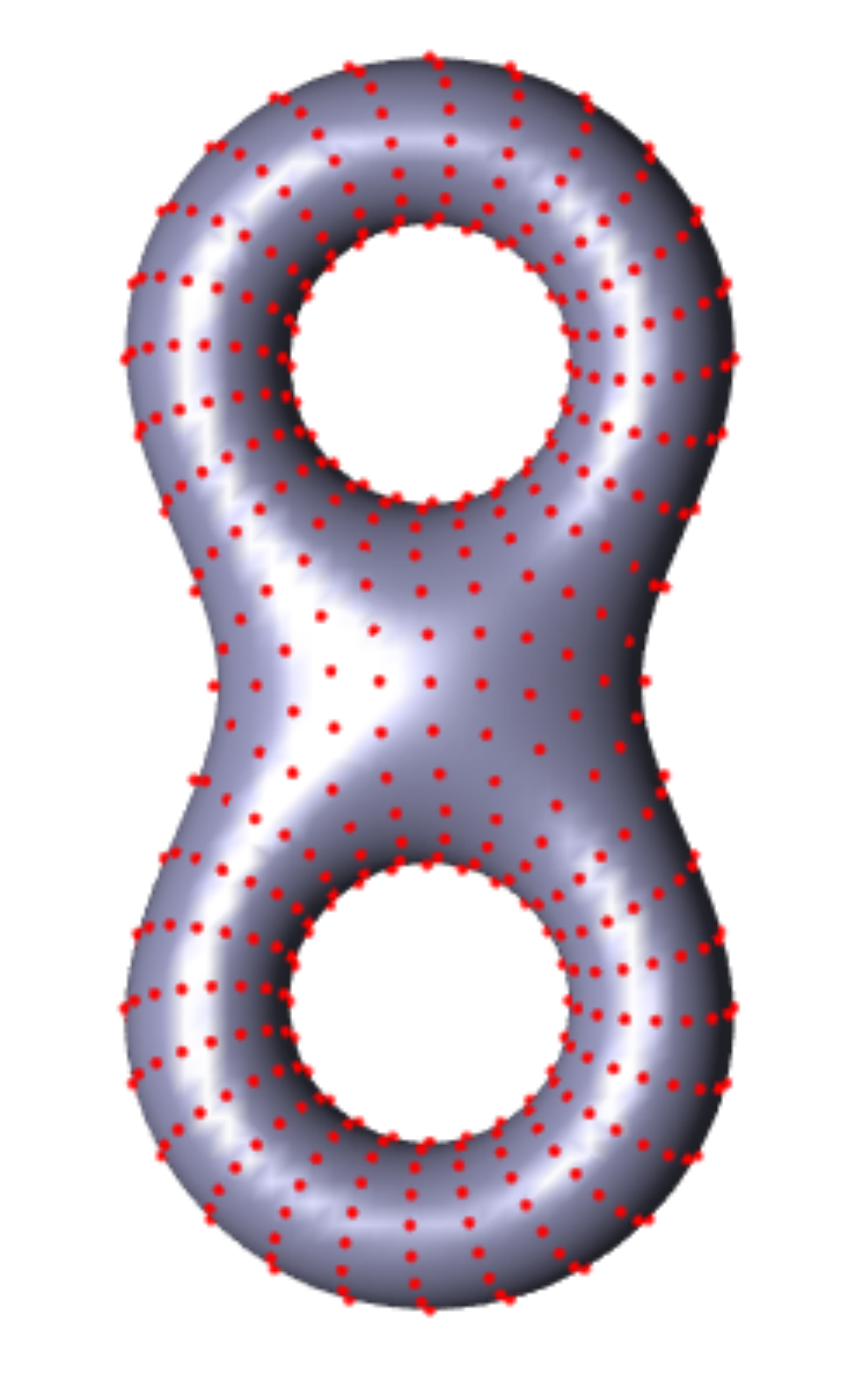}
 }
 \subfigure[Iteration 15.]
 {
    \includegraphics[width=0.23\linewidth]{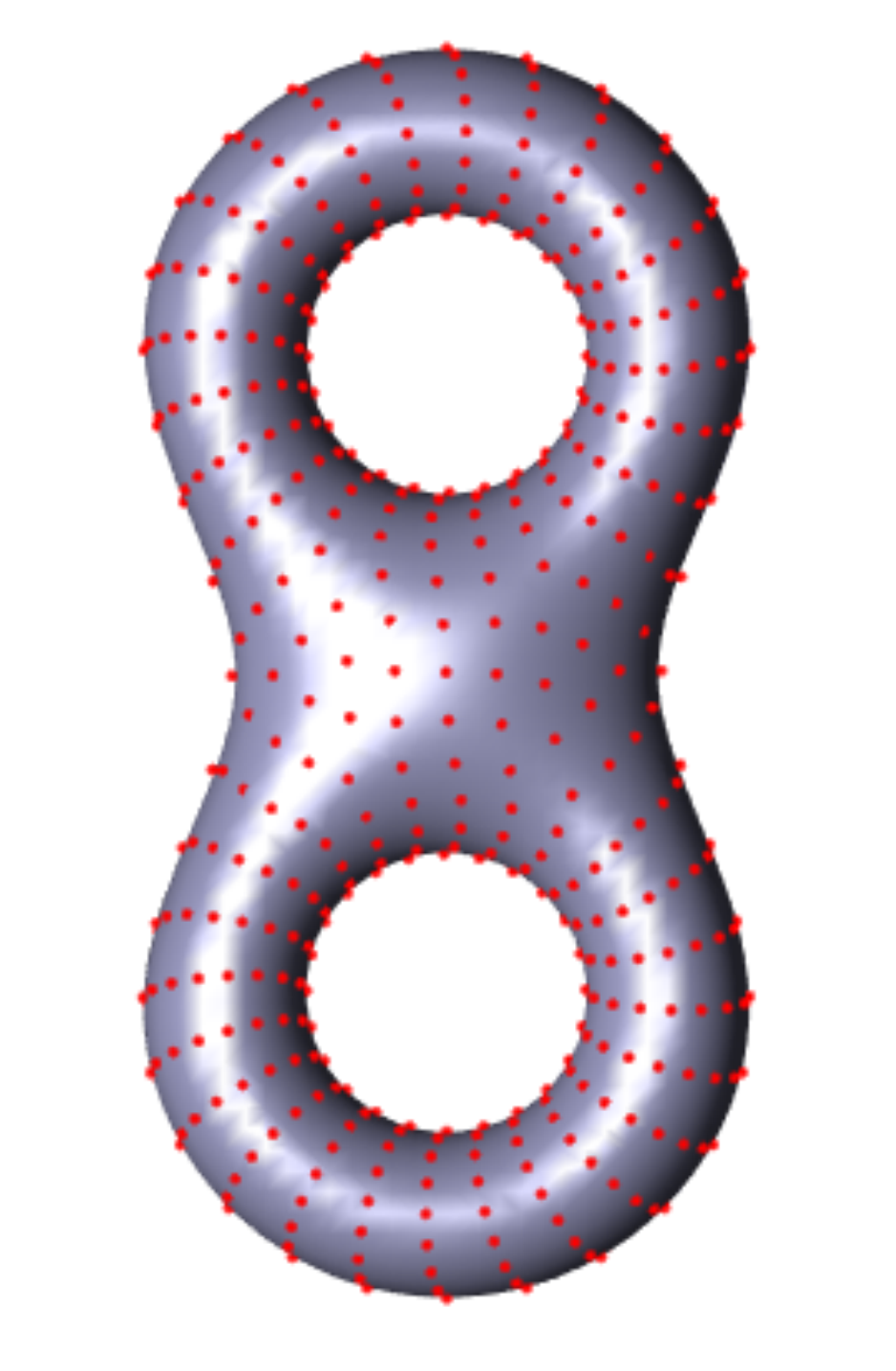}
 }
 \caption
 {
    Iterations in the reconstruction of $3D$ data sets,
    \emph{Torus} and \emph{double-Torus}.
    The red points are the given data sets.
		From left to right: the 1st, 5th, 10th, 15th iteration steps.
		%In the iteration 15th, the implicit surface almost coincides with the
		 %given data sets.
 }
\label{effective-surface}
\end{figure*}

%------------------------------------------------------------------
%%%%%%%%%%%%%%%%%%%%%%%%%%%%%%%%%%%%%%%%%%%%%%%%%%%%%%%%%%%%%%%%%%-
\section{Experiments and discussion}
\label{sec:Numerical Examples}

 Several experiments have been carried out to evaluate the
    performance of I-PIA.
 All the experiments are performed in MATLAB on a PC with an Intel-core i7 @
     3.6 GHz processor and 16 GB of RAM.
 The results and discussions specifically focus on the following areas:
    effectiveness, robustness to inaccurate distance field,
		hole filling, non-uniform sampling and noisy data, porous
    surface (open surface), and fine details.
 The statistical of I-PIA,
    as well as the comparison with the state-of-the-art method~\cite{Liu2017Implicit},
    is listed in Tables~\ref{tab:2d_data} and~\ref{tab:3d_data}.

\subsection{Effectiveness}
\label{sec:Effectiveness}

 As stated above,
    the extra zero-level sets generated in the implicit curve and
    surface reconstruction procedure make the reconstruction results challenging to be interpreted,
    and the elimination of extra zero-level sets is the main problem in designing implicit curve and surface reconstruction method.
 With the I-PIA developed in this paper,
    no extra zero-level set appears in the reconstruction procedure.
 To show the effectiveness of I-PIA in the reconstruction of implicit
    curves and surfaces without the appearance of extra level sets,
    we test our algorithm on planar curves and 3D surfaces.
The initial control coefficients are taken as zero.
After every iteration, the  control coefficients are updated
    and refined by the difference vectors for the control coefficient,
    and the resulting implicit curve/surface will approach the given
 data sets closer than the previous implicit curve/surface.
%After some iterations,
%    the implicit curve/surface will approach the given data sets very close.

 Fig.~\ref{effective-curve}  shows the reconstruction process of
    two $2D$ data sets.
 Similarly, Fig.~\ref{effective-surface} illustrates the reconstruction
     process of two $3D$ data sets.
 From the results presented in Fig.~\ref{effective-curve} and
     Fig.~\ref{effective-surface},
 we can see that no extra zero-level sets exist in the reconstruction
 process of the $2D$ and $3D$ data sets.

\begin{figure*}[!htb]
\centering
 \subfigure[input]
 {
    \includegraphics[width=0.18\linewidth, height = 0.25\linewidth]
        {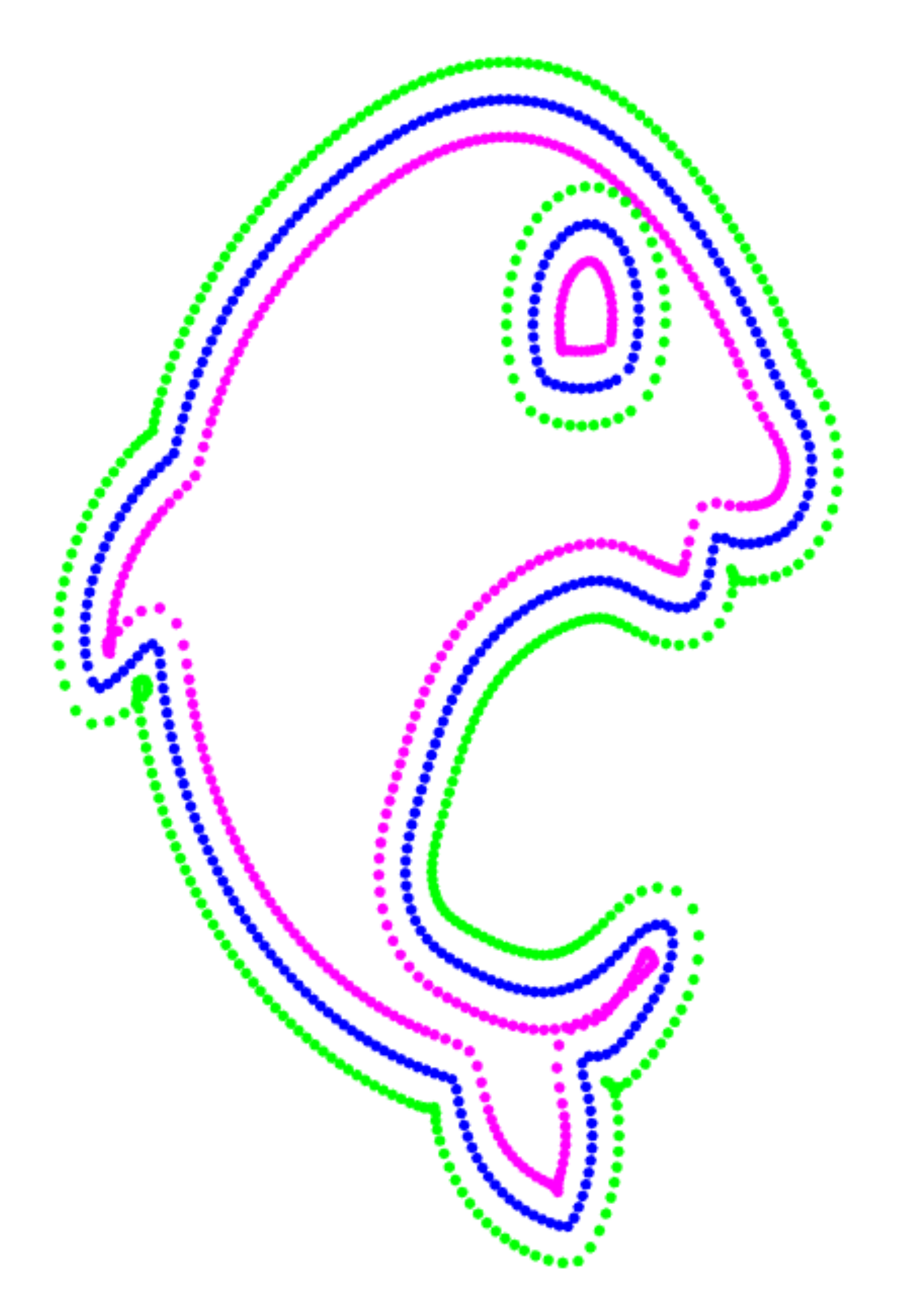}
 }
 \subfigure[$\sigma=0$]
 {
    \label{subfig:robust_0}
    \includegraphics[width=0.175\linewidth, height = 0.25\linewidth]
        {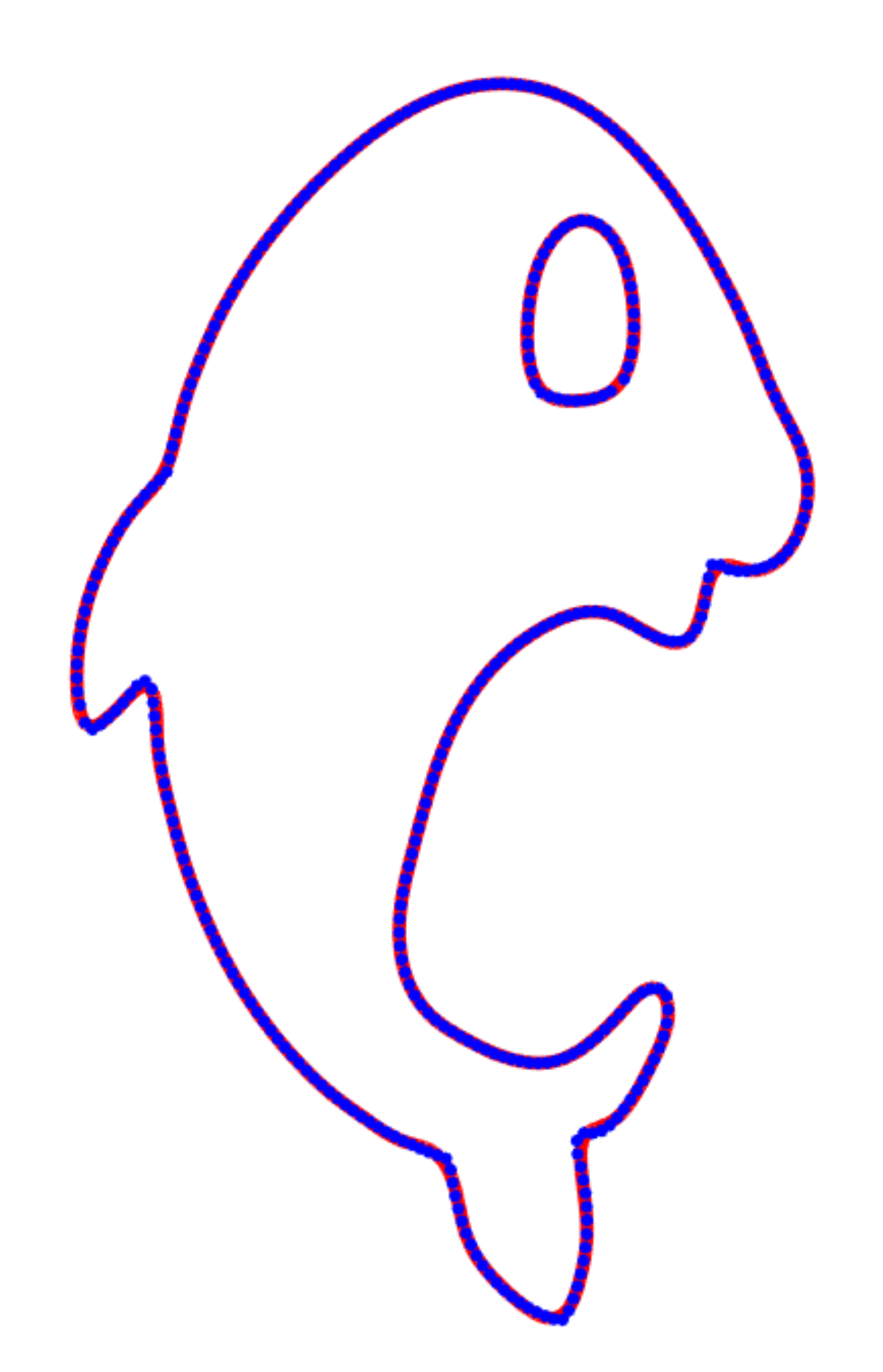}
 }
 \subfigure[$\sigma=0.02$]
 {
    \label{subfig:robust_02}
    \includegraphics[width=0.175\linewidth, height = 0.25\linewidth]
        {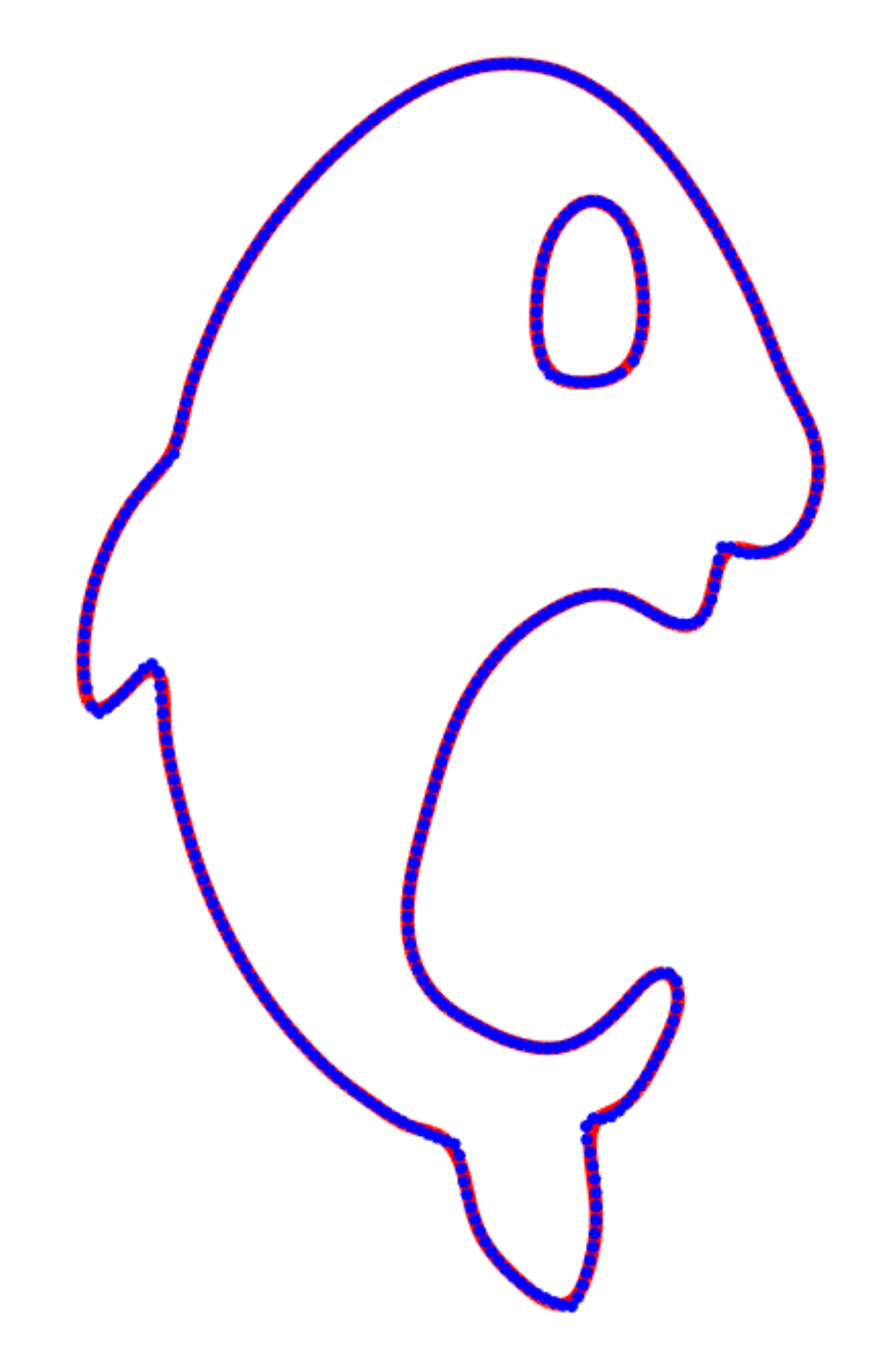}
 }
 \subfigure[$\sigma=0.05$]
 {
    \label{subfig:robust_05}
    \includegraphics[width=0.175\linewidth, height = 0.25\linewidth]
        {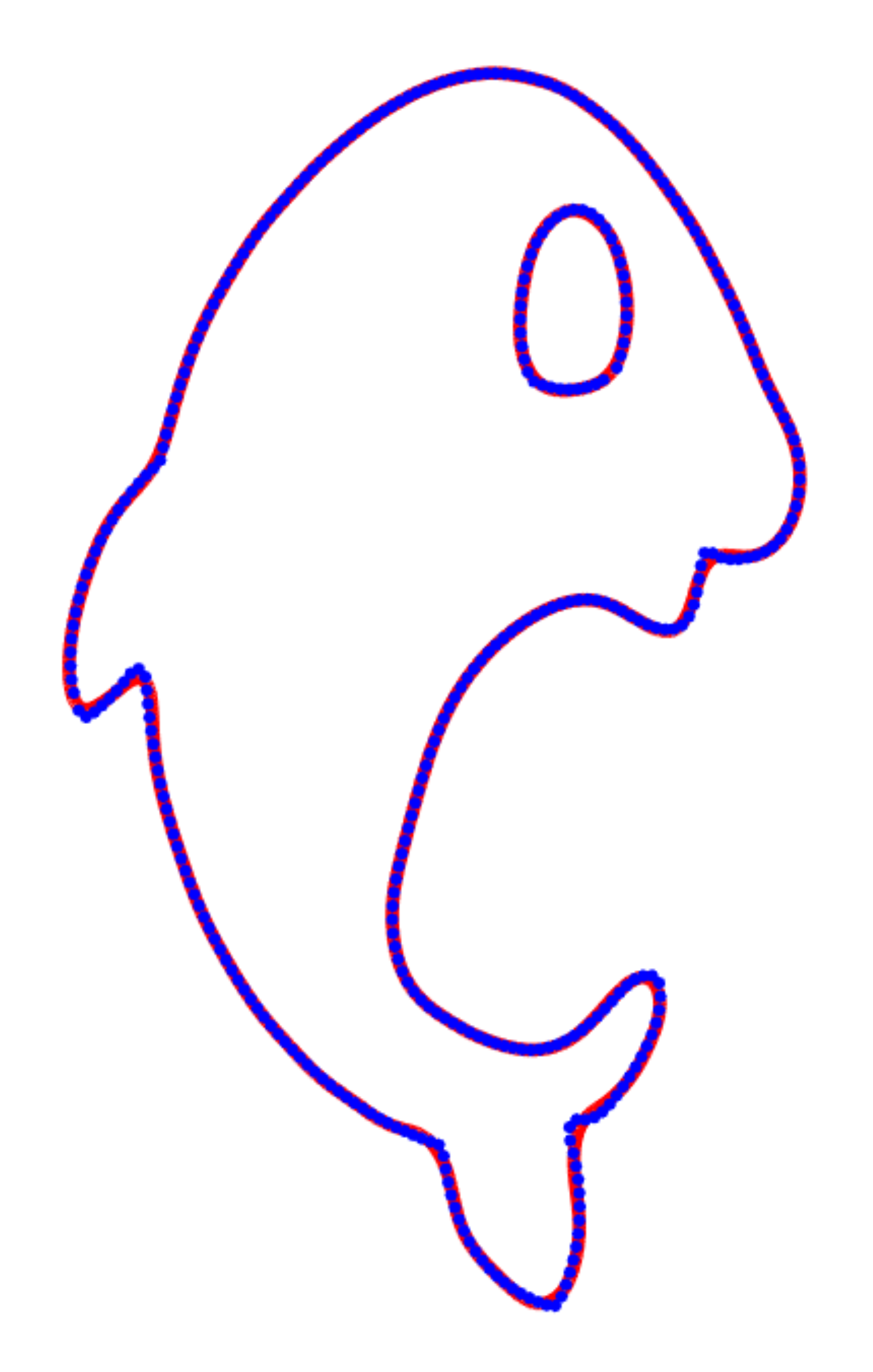}
 }
 \subfigure[$\sigma=0.1$]
 {
    \label{subfig:robust_10}
    \includegraphics[width=0.175\linewidth, height = 0.25\linewidth]
        {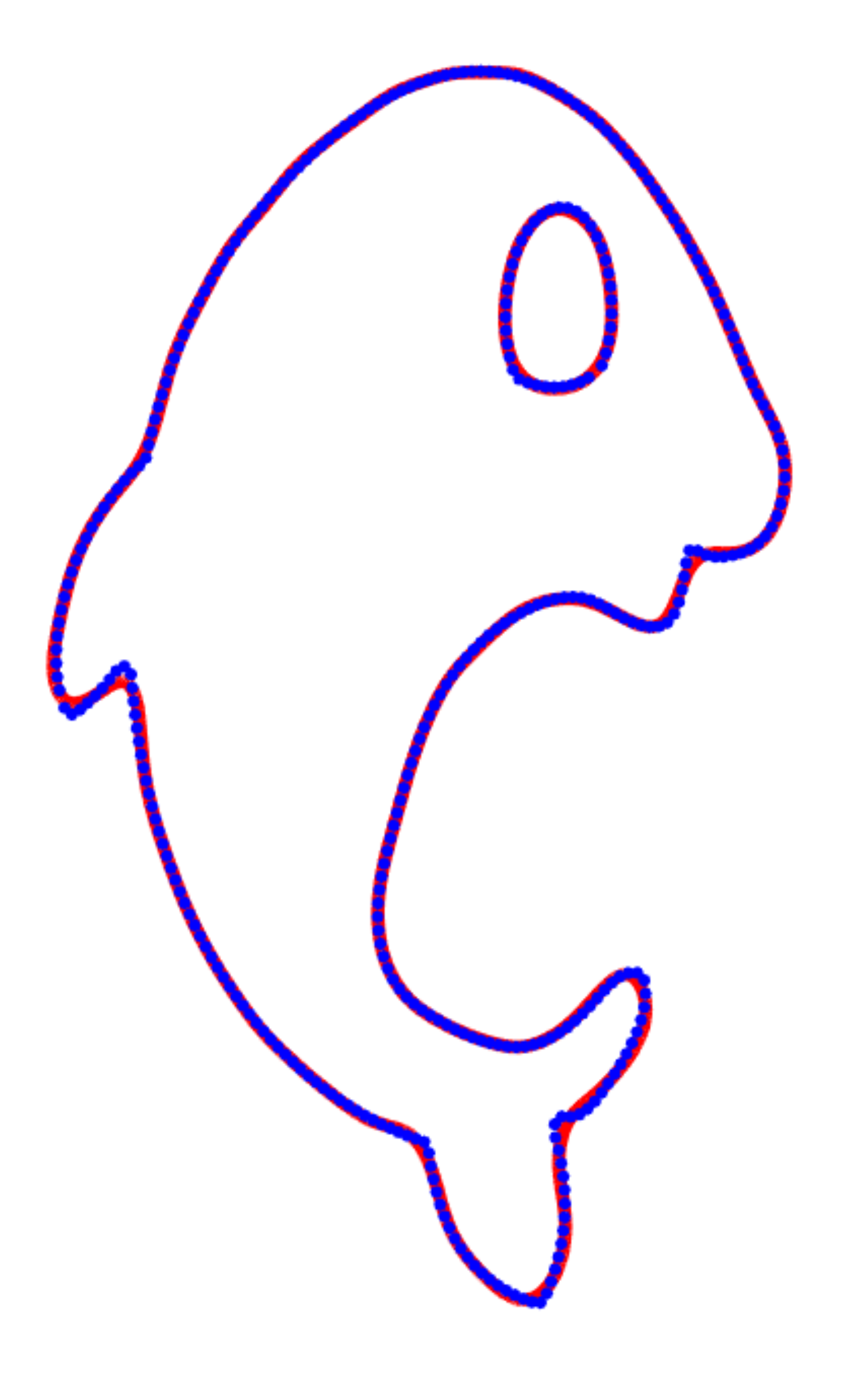}
 }\\
\subfigure[input]
 {
    \label{subfig:robust_butterfly_input}
    \includegraphics[width=0.25\linewidth, height = 0.25\linewidth]
        {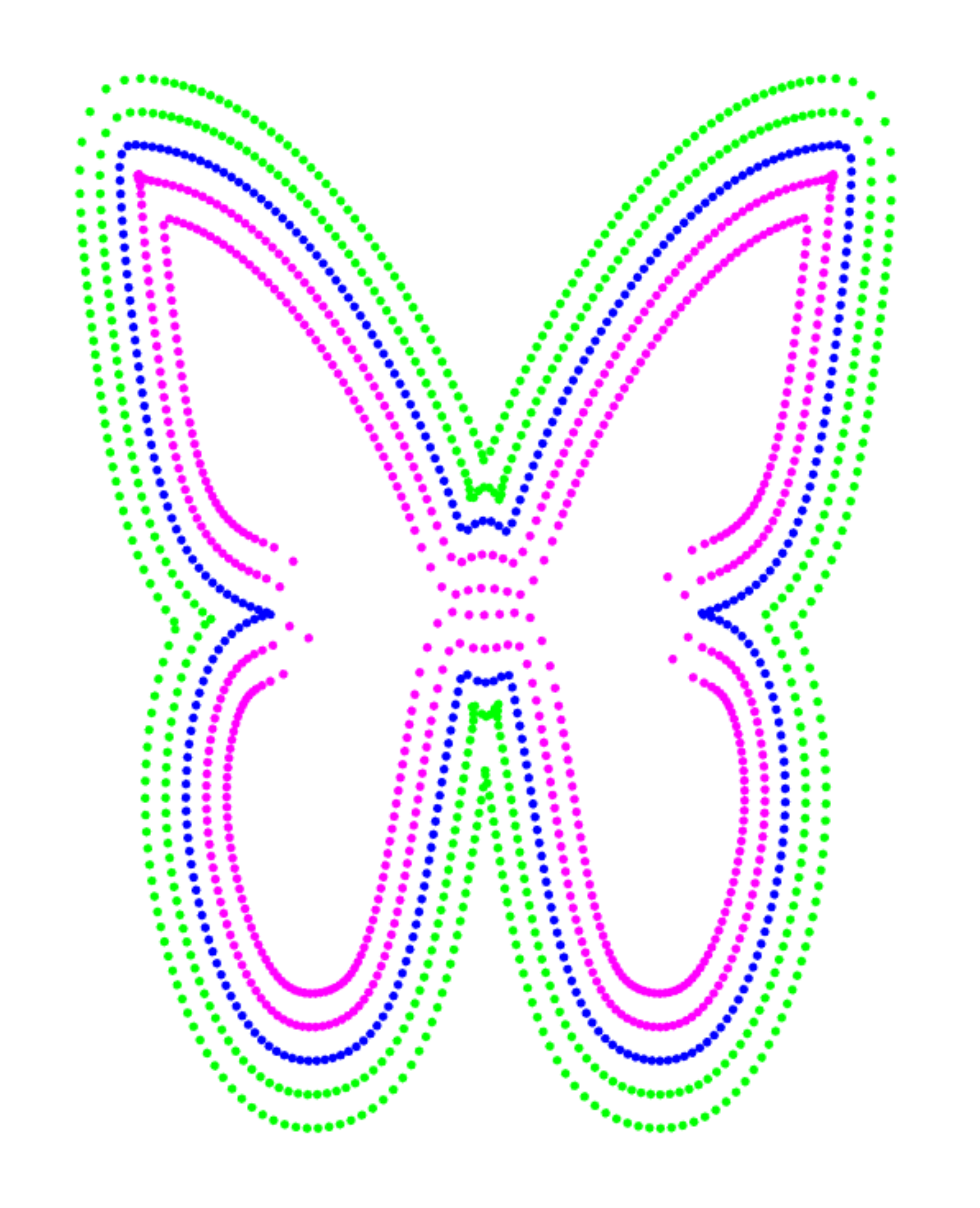}
 }
 \subfigure[reconstructed curve]
 {
    \label{subfig:robust_butterfly_curve}
    \includegraphics[width=0.27\linewidth, height = 0.25\linewidth]
        {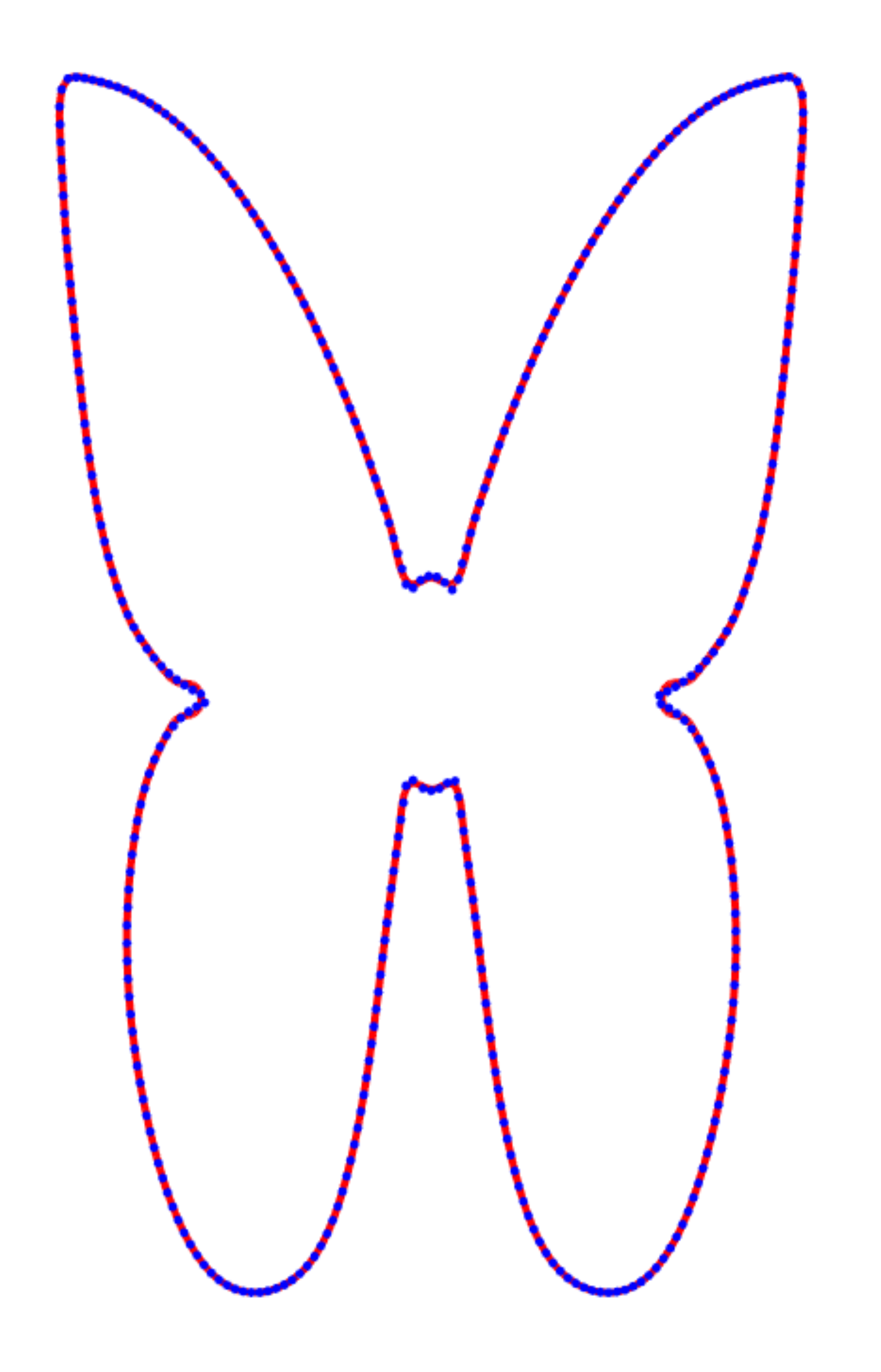}
 }
 \subfigure[input]
 {
    \label{subfig:robust_dolphin_input}
    \includegraphics[width=0.18\linewidth, height = 0.25\linewidth]
        {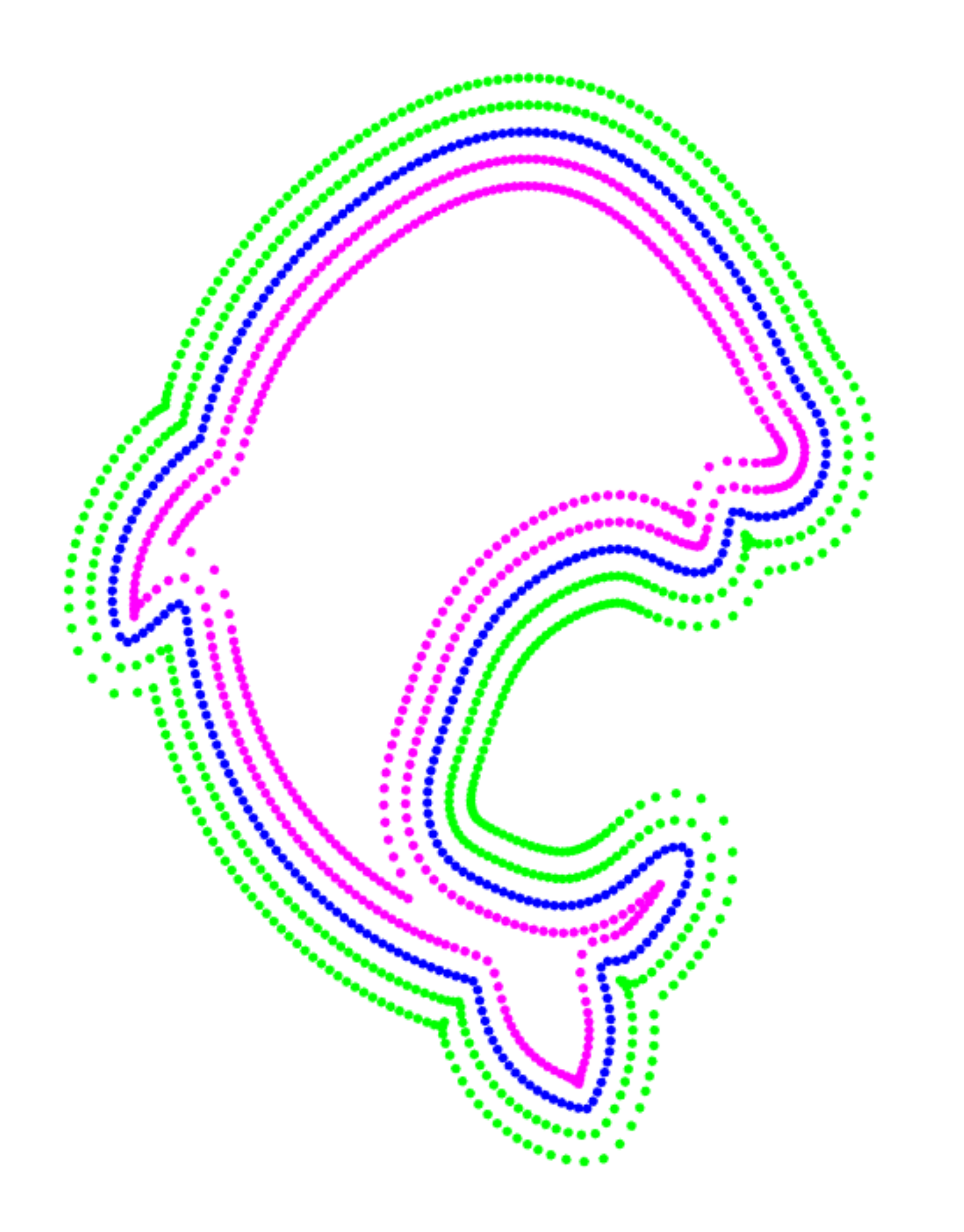}
 }
  \subfigure[reconstructed curve]
 {
    \label{subfig:robust_dolphin_curve}
    \includegraphics[width=0.175\linewidth, height = 0.25\linewidth]
        {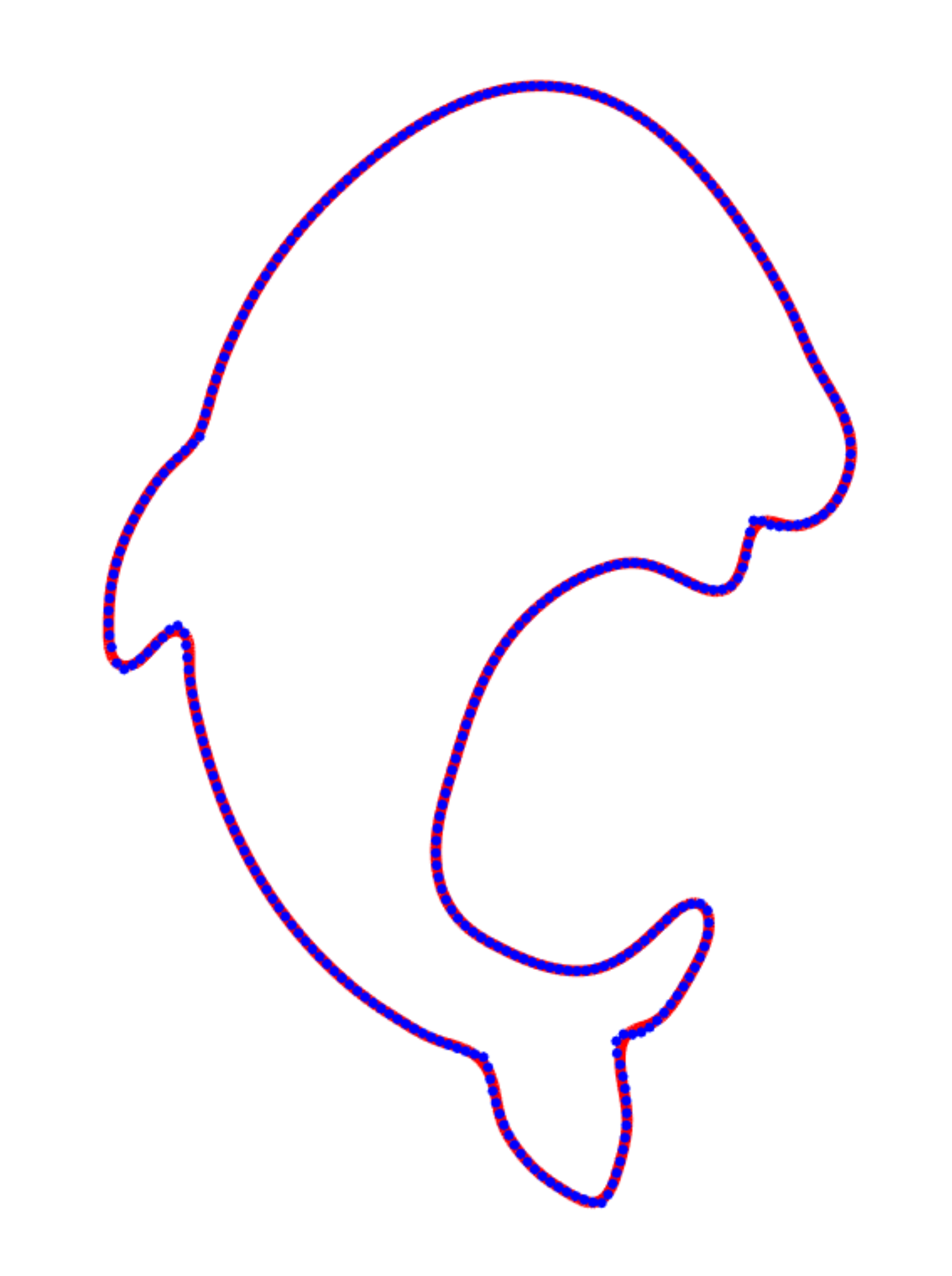}
 }
 \caption
 {
    Reconstruction of $2D$ data sets with inaccurate distance fields.
    First row: (a) the input data of \emph{dolphin} (blue) and offset points
        (outside offset in green and inside offset in magenta).
    The preassigned function values selected as uniformly distributed random number in $[\pm 0.5-\sigma, \pm 0.5+\sigma]$ for the outside and inside offsets, respectively, (b)--(e): $\sigma=0, 0.02, 0.05, 0.1$.
	Second row: (f) and (h) the input data of \emph{butterfly} and
    \emph{dolphin} (blue), respectively,
	and offset points (outside offset in green and inside offset in magenta).
	Larger value are assign to closer offset point and vise versa.
	From inside to outside, the distance values of the four sets of
        offset points are $-0.1,-0.2,0.2,0.1$.
 (g) and (i) the reconstructed curves from (f) and (h),
			respectively.
	%The reconstructed curves (red lines), and the given data set (blue points).
 }
\label{inaccurate-distance}
\end{figure*}

\subsection{Robustness to inaccurate distance field}
\label{sec:Robustness}

 The auxiliary offset points appended to the data points helps to orient
	   the surface and avoid the appearance of artifacts in curve and surface
	   reconstruction~\cite{Carr2001Reconstruction,Liu2017Implicit}.
 The I-PIA algorithm is insensitive to the distance values assigned at the
    offset points,
    thus robust to the inaccurate distance field.
 The first row of Fig.~\ref{inaccurate-distance} illustrates the
    reconstruction of a $2D$ data sets with synthesis offset points.
	%The distance value to the curve is inaccurately computed and
 The  function values on the offset points are assigned in a random strategy.
 %That is, the function values are initially set to $0.5$ for the
%        outside offsets and $-0.5$ for the inside offsets.
 That is, the preassigned values are selected as uniformly distributed
    random numbers in $[0.5-\sigma, 0.5+\sigma]$ for the outside offsets
    and $[-0.5-\sigma, -0.5+\sigma]$ for the inside offset,
	respectively,
    where $\sigma=0, 0.02, 0.05, 0.1$,
    as demonstrated in Figs.~\ref{subfig:robust_0}-~\ref{subfig:robust_10}.
 The results show that I-PIA can still reconstruct the given
	  data sets in a corrupted distance field.
		
 Moreover, given a data points set,
    we generate two sets of offset points inside and
	outside the data points set respectively.
 However, the distance values are assigned in reverse order, i.e.,
		larger distance values are given to nearer offset points,
	  and smaller distance values are given to far offset points.
 As illustrated in Figs.~\ref{subfig:robust_butterfly_input}
    and~\ref{subfig:robust_dolphin_input},
	from inside to outside, the distance values of the four sets of
        offset points are $-0.1,-0.2,0.2,0.1$.
 The curves is robustly reconstructed using I-PIA,
    and demonstrated in Figs.~\ref{subfig:robust_butterfly_curve}
    and~\ref{subfig:robust_dolphin_curve}.

\subsection{Holes filling}
\label{sec:HolesFilling}

\begin{figure*}[!htb]
\centering
 \subfigure[]
 {
    \includegraphics[width=0.24\linewidth]{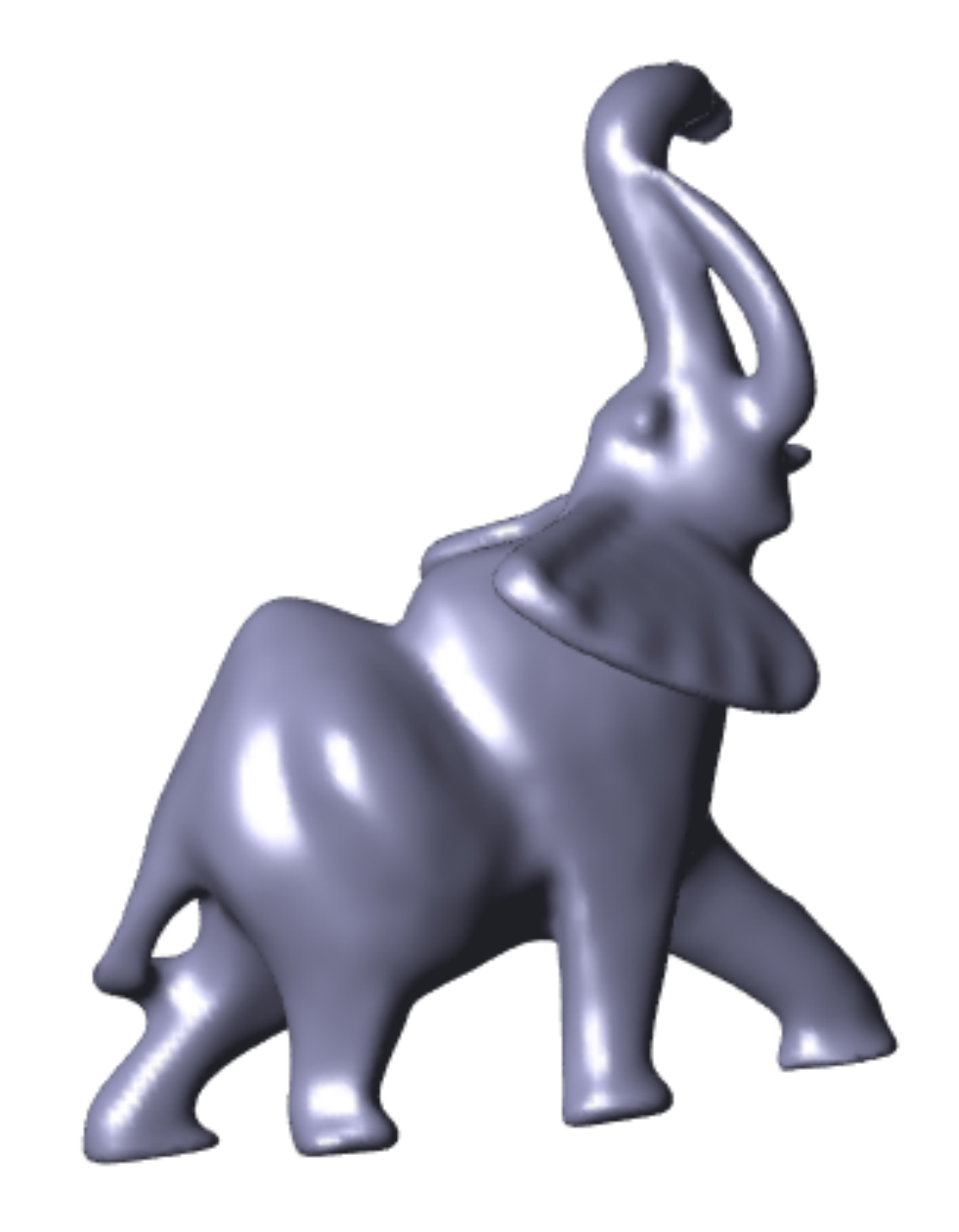}
 }
 \subfigure[]
 {
    \includegraphics[width=0.24\linewidth]{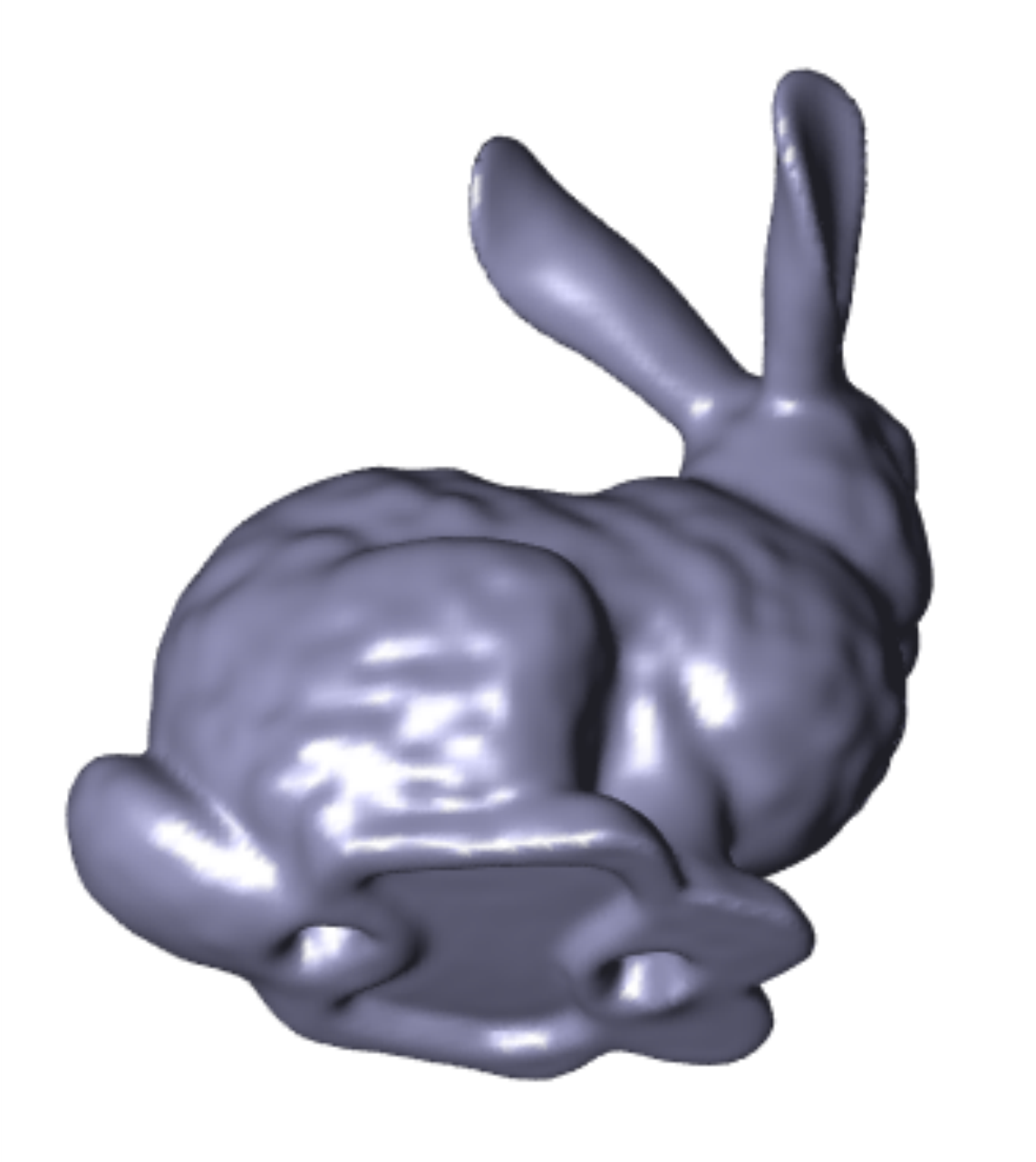}
 }
 \subfigure[]
 {
    \includegraphics[width=0.27\linewidth]{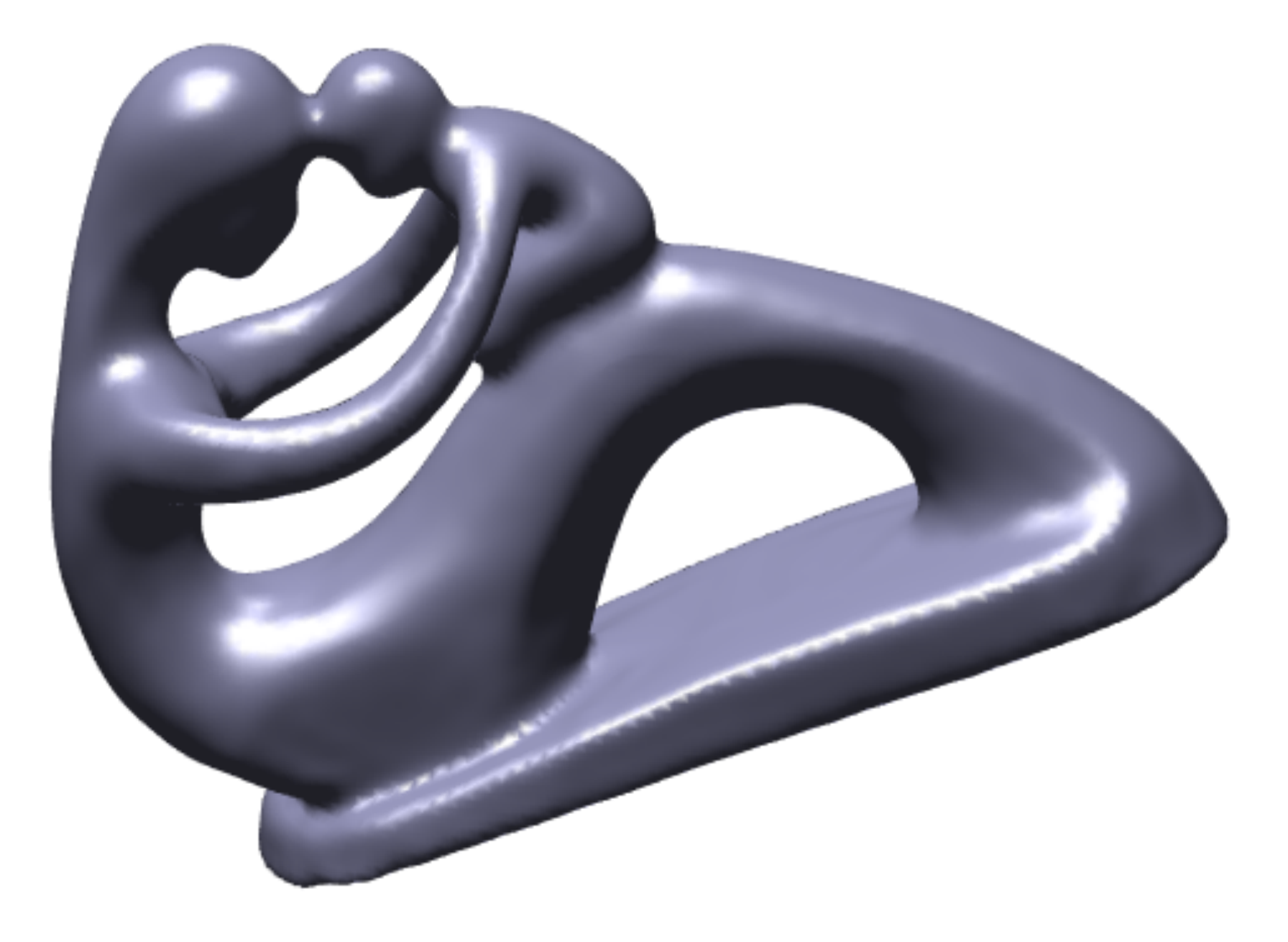}
 }\\
 \subfigure[]
 {
    \includegraphics[width=0.24\linewidth]{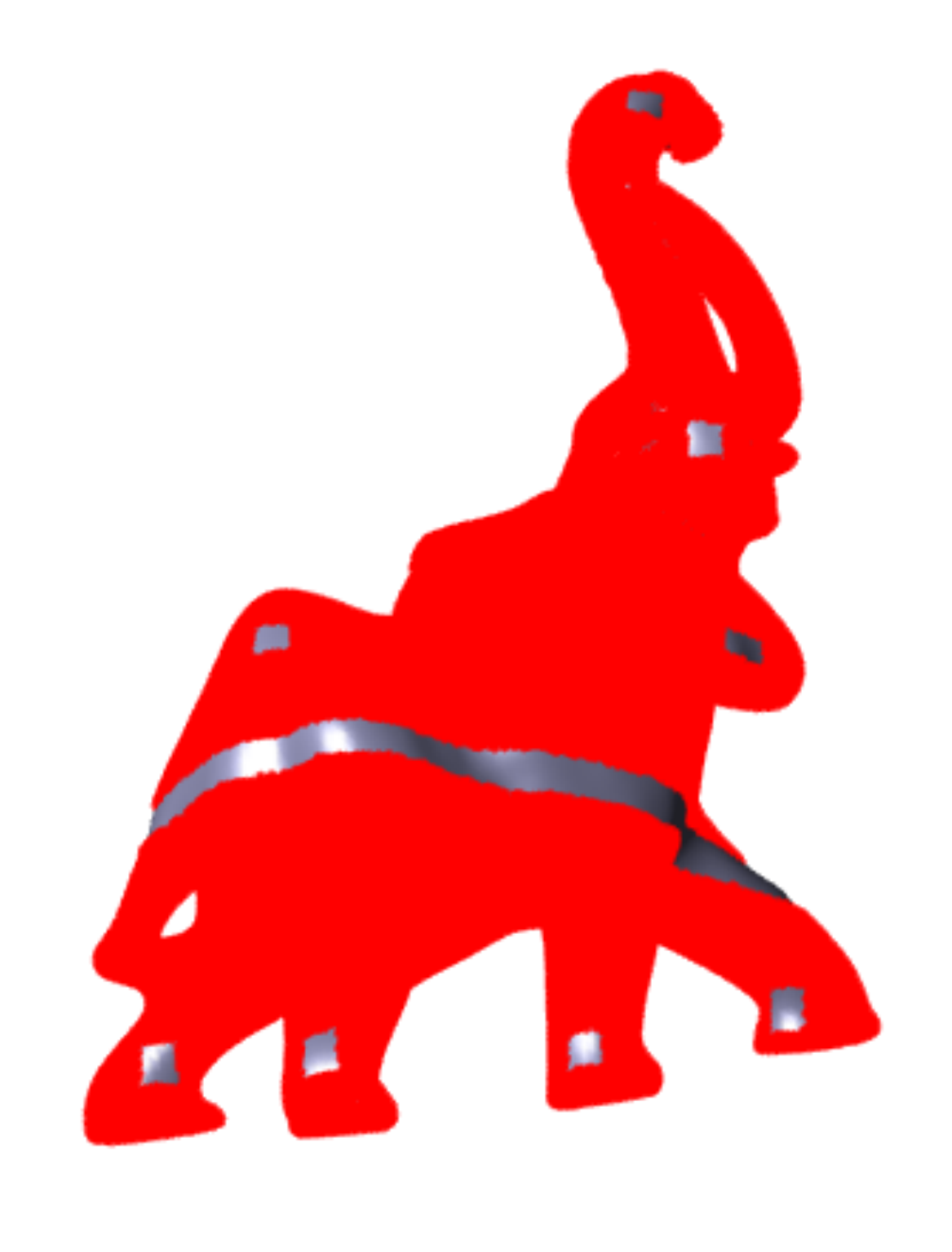}
 }
\subfigure[]
 {
    \includegraphics[width=0.24\linewidth]{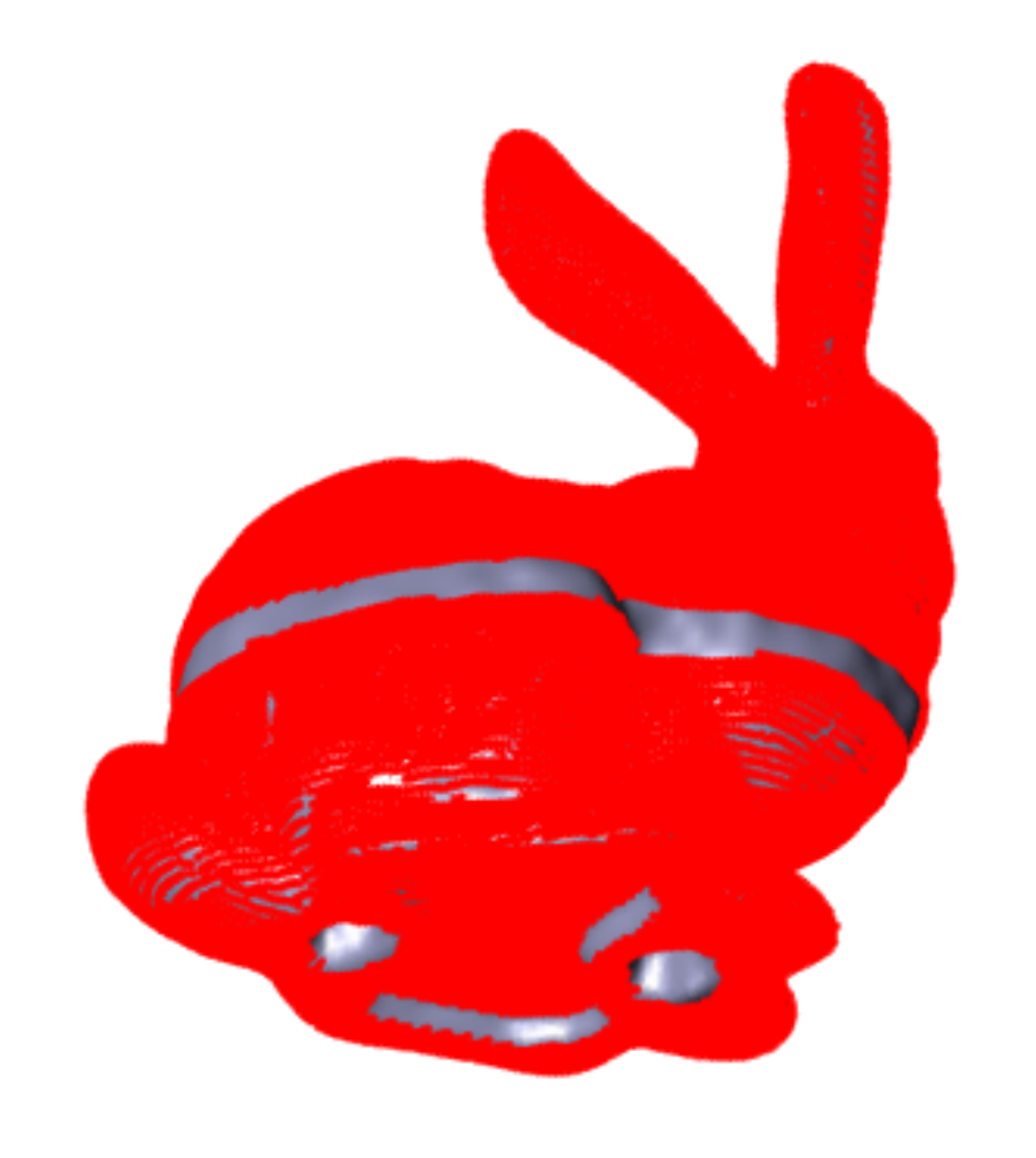}
 }
 \subfigure[]
 {
    \includegraphics[width=0.27\linewidth]{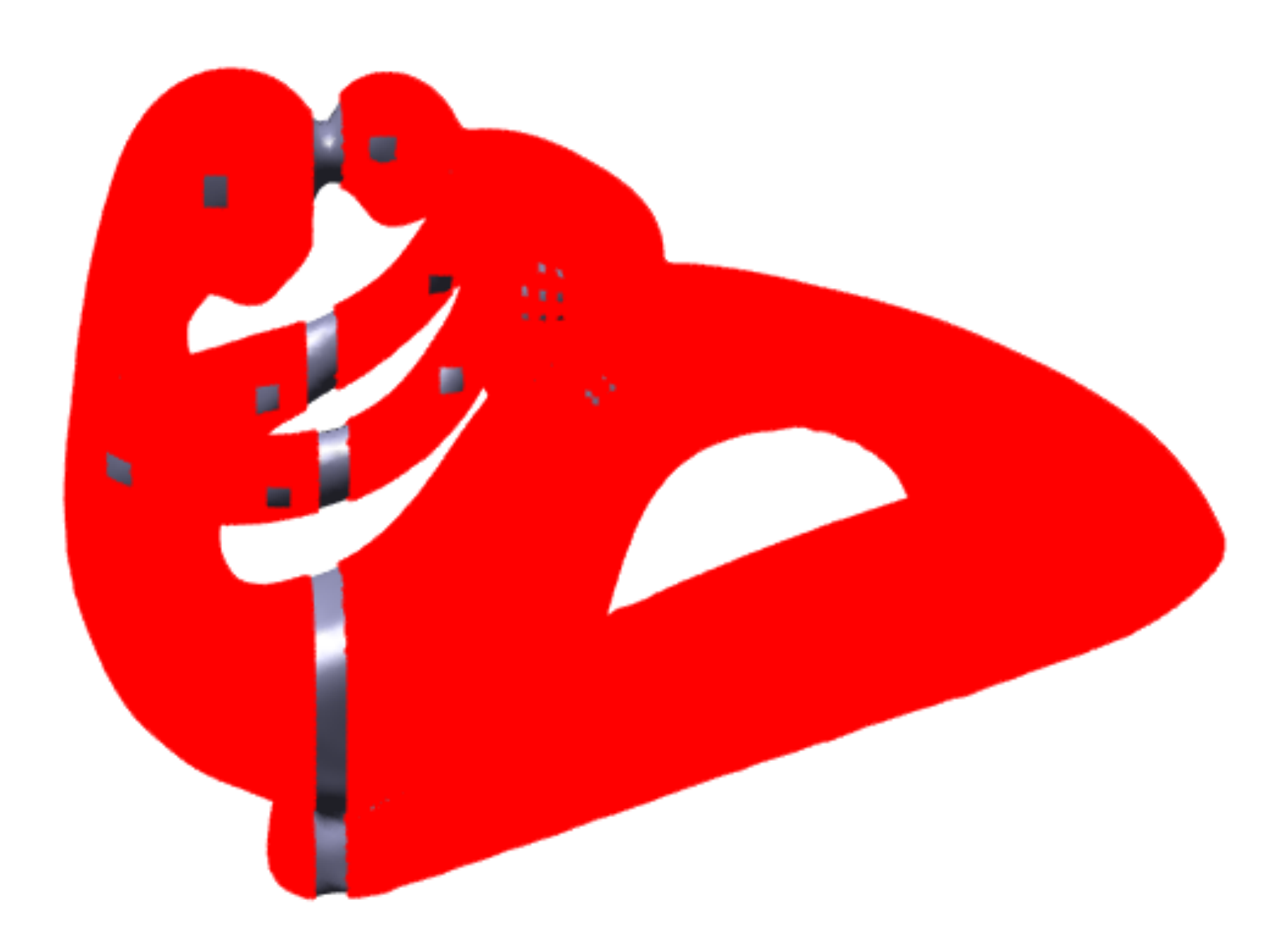}
 }
 \caption
 { Holes and gaps filling. (a,d) \emph{elephant}; (b,e) \emph{bunny}; (c,f) \emph{fertility}. First row: reconstructed surfaces. Second row: the reconstructed surfaces with the data points superimposed.
 }
\label{holes}
\end{figure*}

 Holes and missing data often arise in the point cloud generated by
    scanning devices.
 Fig.~\ref{holes} (a)--(c) depicts the reconstructed surfaces of elephant,
  fertility, and bunny respectively,
  in which some portions of the data points are removed to creates holes and gaps.
 The reconstructed surfaces are illustrated in the first row, and the second row shows the reconstructed surfaces with the data points superimposed.
 As can be seen, I-PIA successfully filled the holes and gaps.
 Thus, I-PIA performs well in filling holes and missing data.

\subsection{Non-uniform sampling and noisy data}
\label{sec:DownSampledAndNoisyData}

\begin{figure*}[!htb]
\centering
 \subfigure[]
 {
    \label{subfig:bunny_original}
    \includegraphics[width=0.27\linewidth]{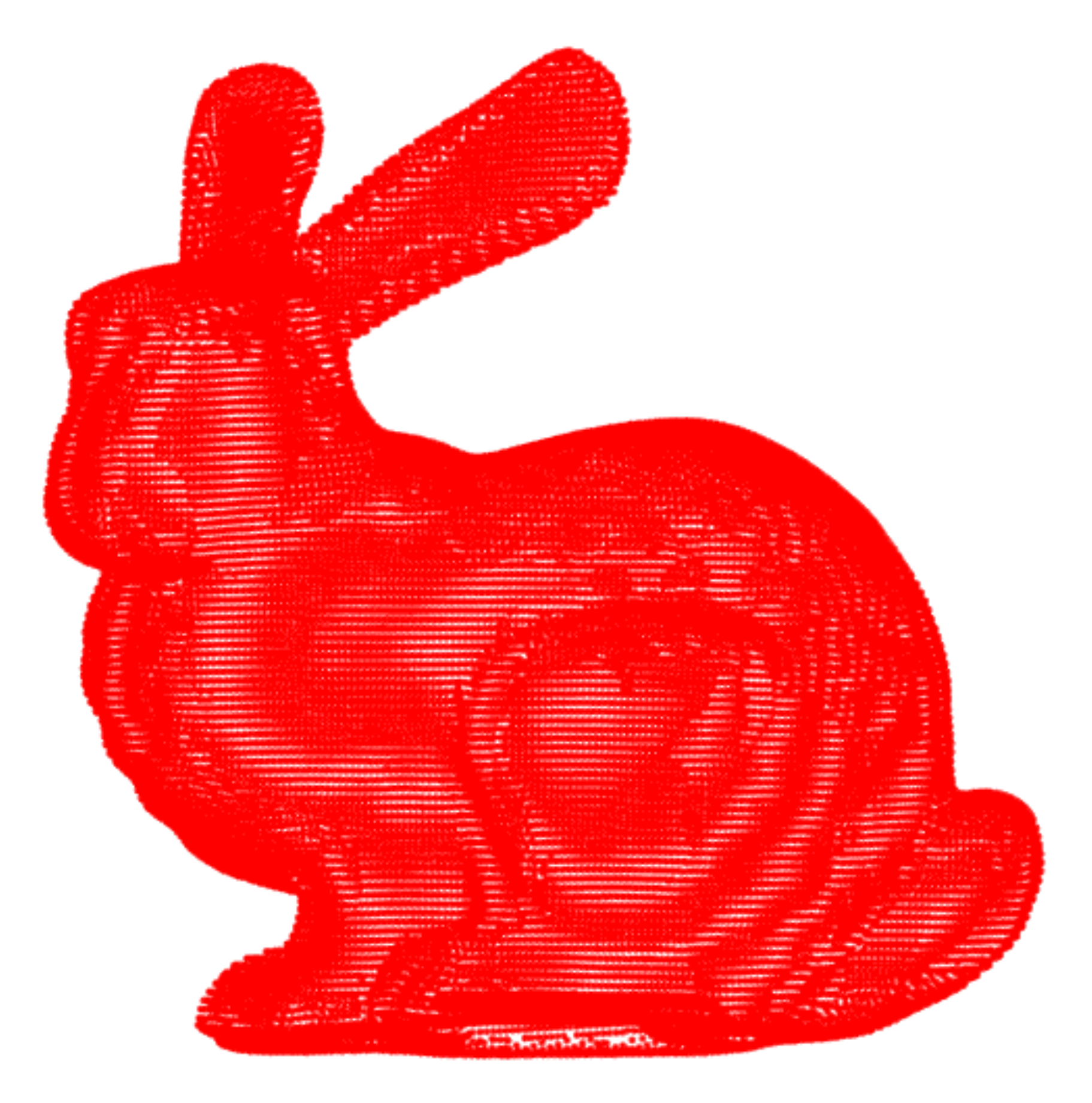}
 }
 \subfigure[]
 {
    \label{subfig:bunny_downsample}
    \includegraphics[width=0.25\linewidth]{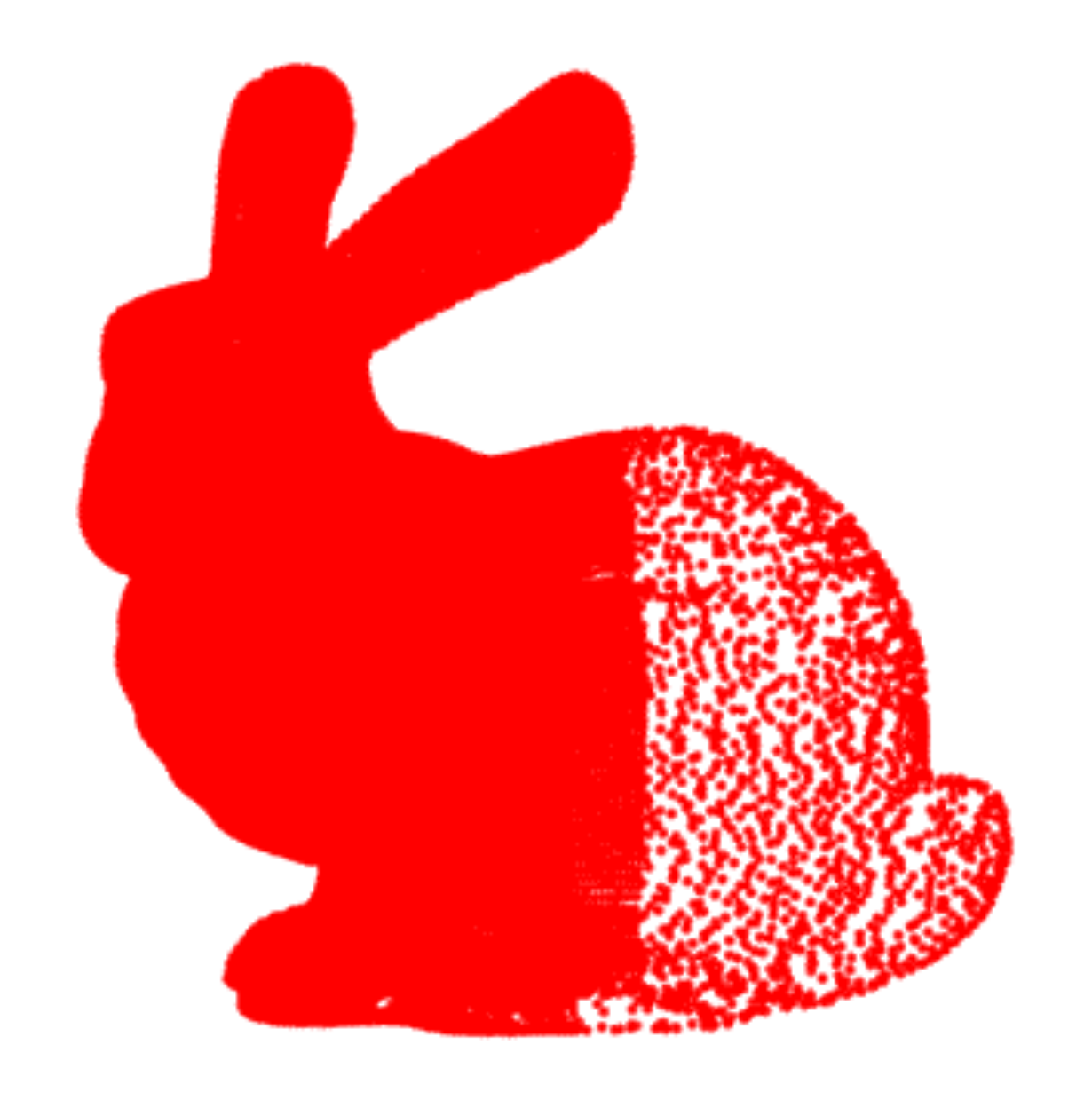}
 }
 \subfigure[]
 {
    \label{subfig:bunny_noise}
    \includegraphics[width=0.25\linewidth]{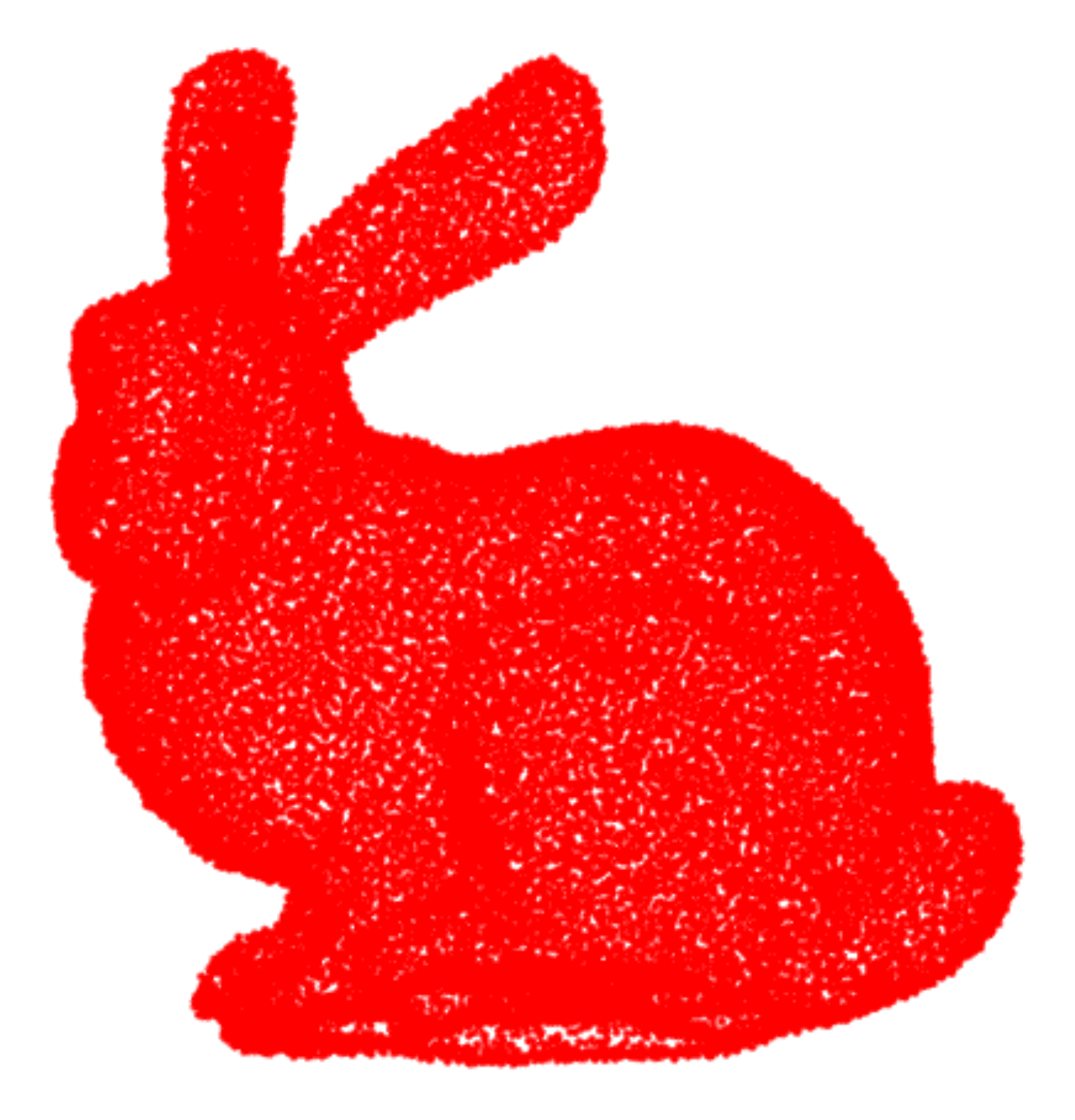}
 }
 \subfigure[]
 {
    \label{subfig:bunny_original_recons}
    \includegraphics[width=0.23\linewidth]{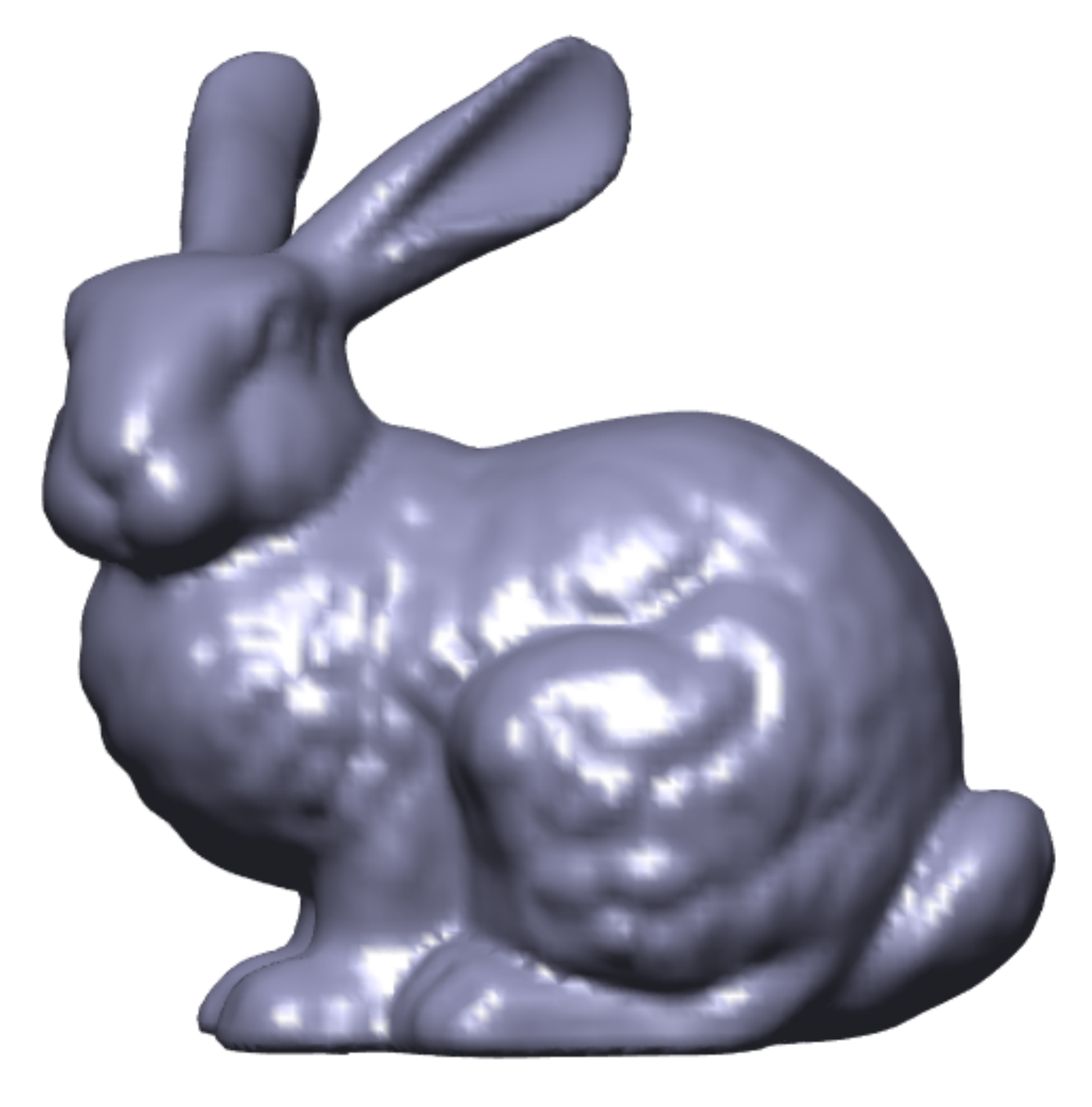}
 }
 \subfigure[]
 {
    \label{subfig:bunny_downsample_recons}
    \includegraphics[width=0.25\linewidth]{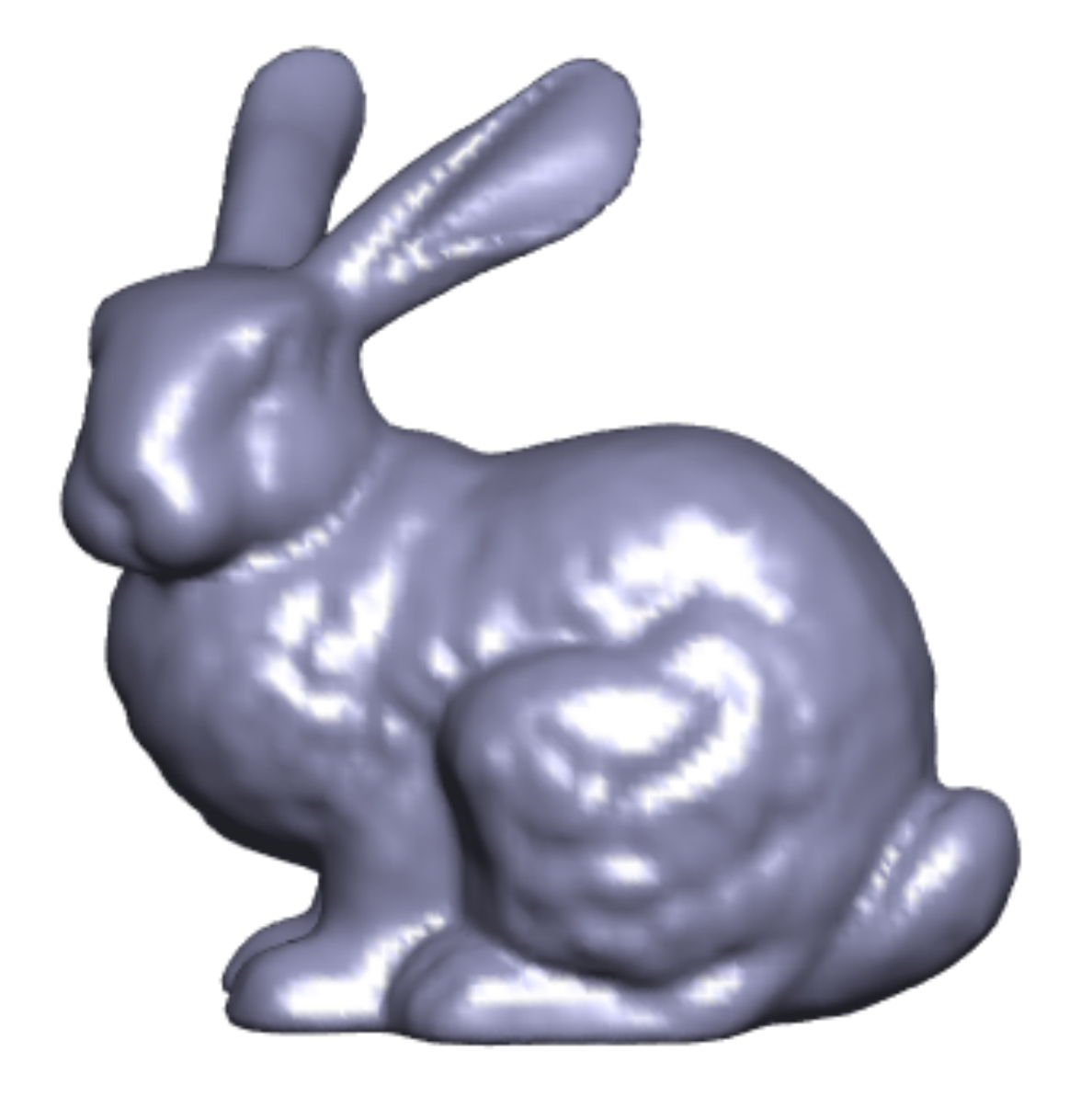}
 }
\subfigure[]
 {
    \label{subfig:bunny_noise_recons}
    \includegraphics[width=0.25\linewidth]{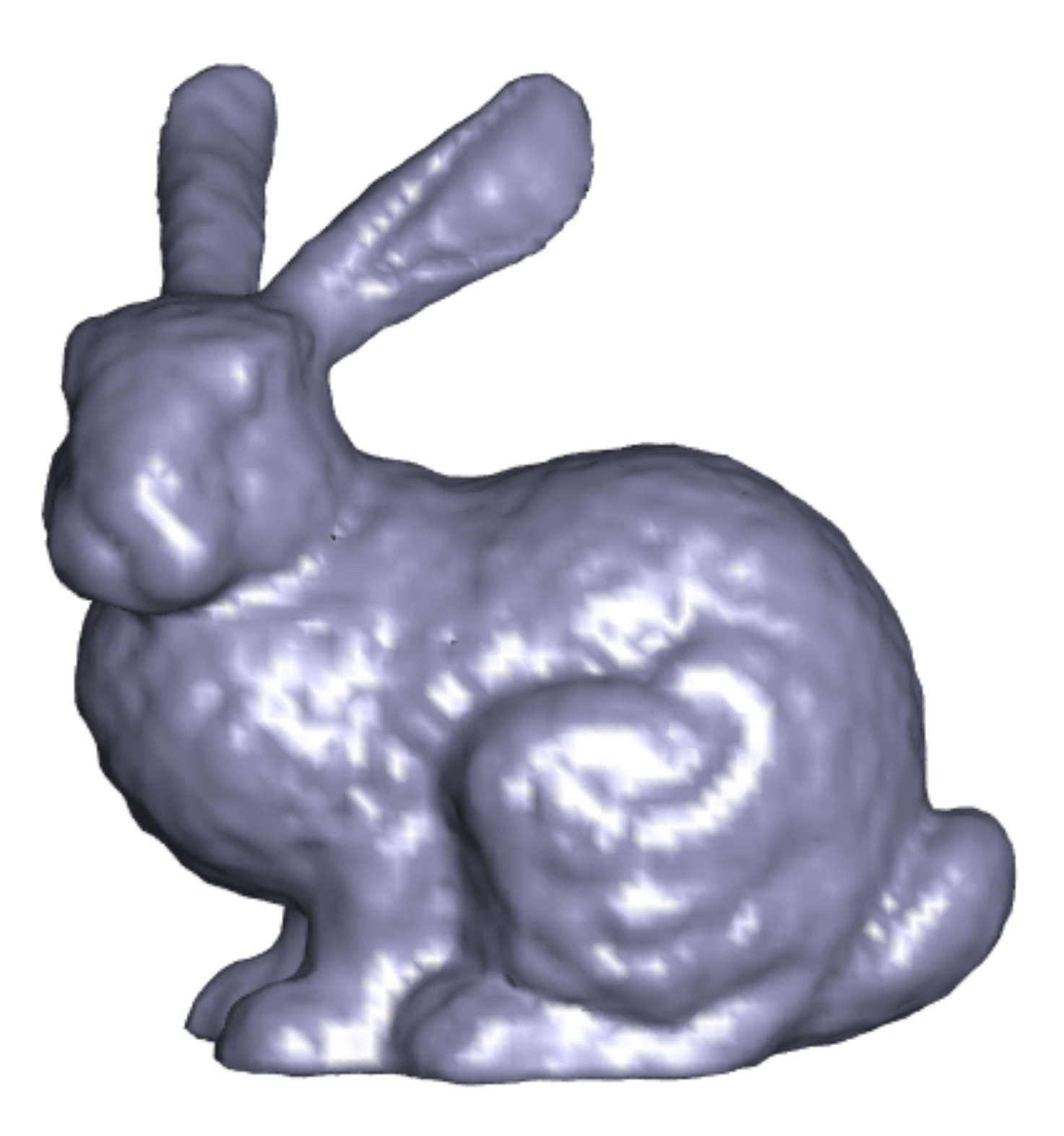}
 }
 \caption
 {
    Surface reconstructions of the \emph{bunny} model from different point clouds.
    (a) The original point cloud of the \emph{bunny} model.
    (b) 90\% of the right part of the \emph{bunny} model are removed.
    (c) The original point cloud is disturbed by noise.
	(d)--(f) are the reconstructed surfaces from (a)--(c) by I-PIA,
    respectively.}
\label{noise}
\end{figure*}

 In this section, we demonstrates the capability of I-PIA in handling
    non-uniform sampling and noisy data.
 In Fig.~\ref{subfig:bunny_original},
    the original point cloud of the bunny model is illustrated.
 In Fig.~\ref{subfig:bunny_downsample},
    the right part of the bunny model is down-sampled,
    and $90\%$ of the original data points are removed.
 Moreover, in Fig.~\ref{subfig:bunny_noise},
    the data points are disturbed by noise.
 However, by I-PIA,
    the input point clouds in  Fig.~\ref{subfig:bunny_downsample} and Fig.~\ref{subfig:bunny_noise} do not cause artifacts
   in the reconstructed surfaces in Fig.~\ref{subfig:bunny_downsample_recons} and Fig.~\ref{subfig:bunny_noise_recons}, respectively.
	
\begin{figure*}[!htb]
\centering
 \subfigure[]
 {
    \label{subfig:arm_original}
    \includegraphics[width=0.20\linewidth]{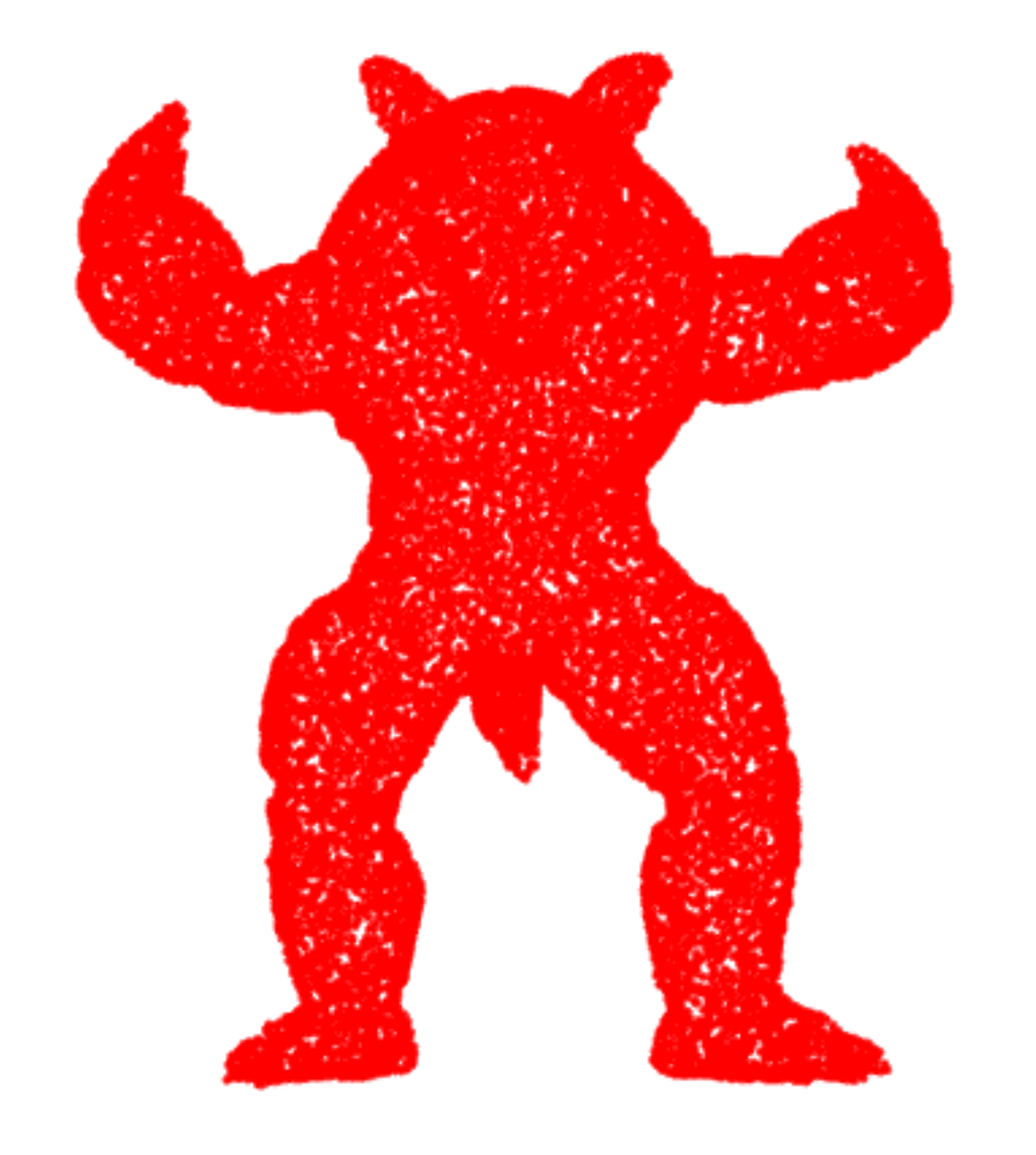}
 }
 \subfigure[]
 {
    \label{subfig:arm_reduced}
    \includegraphics[width=0.20\linewidth]{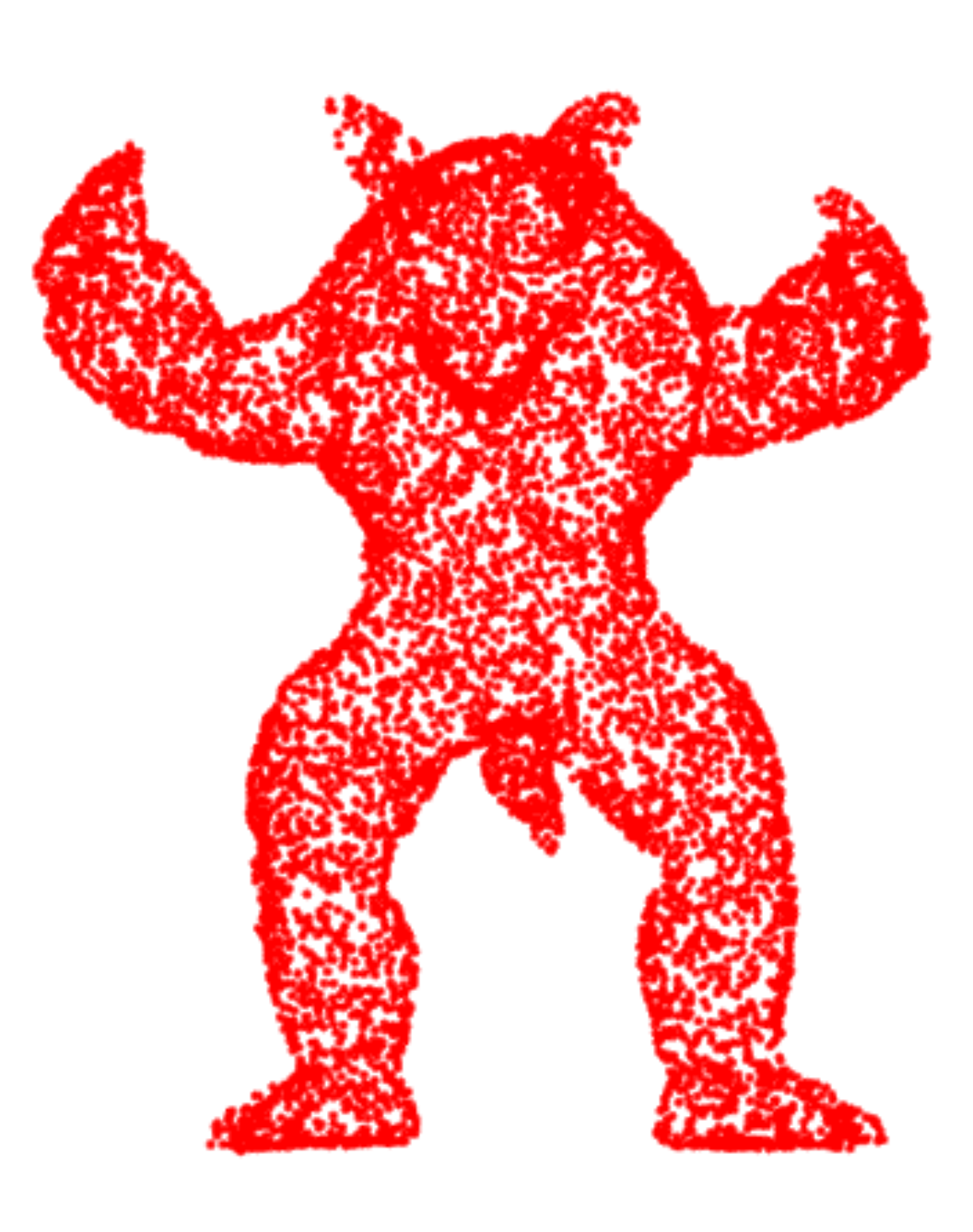}
 }
 \subfigure[]
 {
    \label{subfig:arm_recons_100}
    \includegraphics[width=0.22\linewidth]{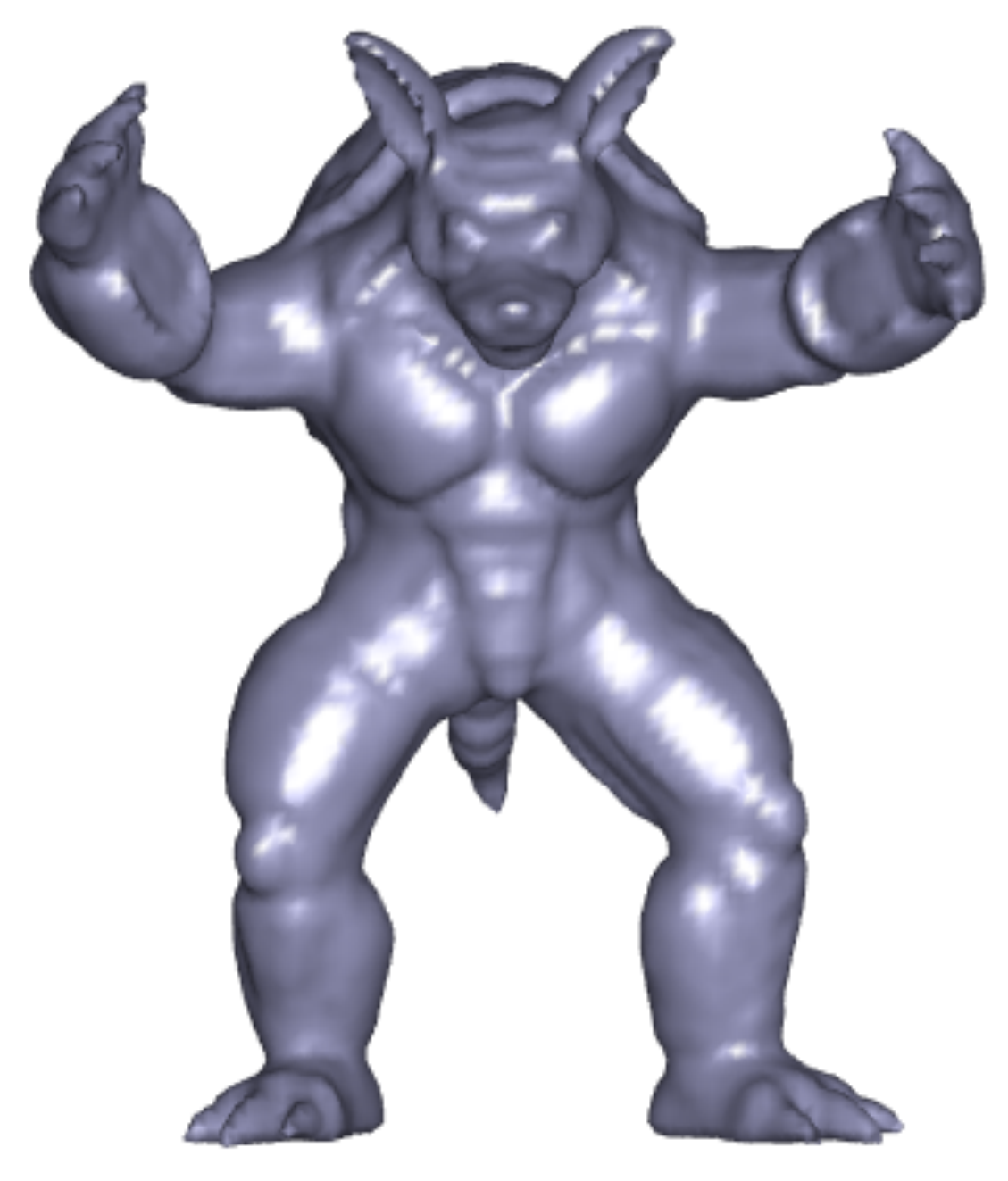}
 }
 \subfigure[]
 {
    \label{subfig:arm_recons_190}
    \includegraphics[width=0.22\linewidth]{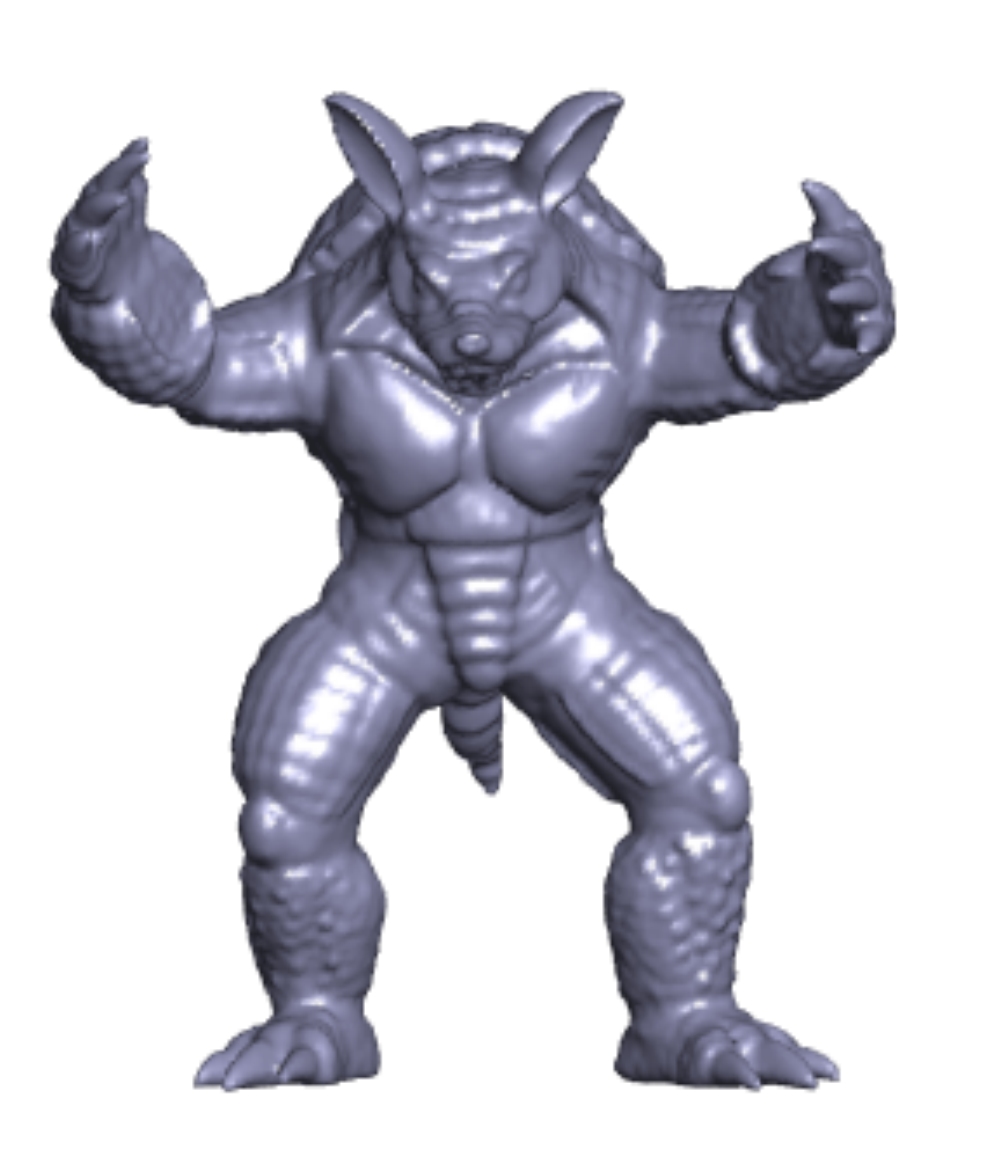}
 }

 \caption
 {
    Surface reconstructed from point clouds with different sampling densities of the \emph{armadillo} model.
    (a) and (b) Input point clouds with different sampling densities.
    (c) Reconstructed surface from the point cloud in (a)  based on a coarse grid ($100\times100\times100$).
    (d) Reconstructed surface from the point cloud in (b) based on a fine grid ($190\times190\times190$).
		%For better visualization, the points are displayed more sparsely than the real density.
		}
\label{lowdensity}
\end{figure*}

 The features and details of the reconstructed surface by I-PIA depends
    on the mesh grid used.
 If the density of the input point cloud is reduced,
    I-PIA can still reconstruct the surface with fine details and features if the mesh grid is refined.
 Fig.~\ref{lowdensity} shows the reconstructed surfaces from
    the point clouds of the model \emph{armadillo} with
    different density,
    where Fig.~\ref{subfig:arm_original} is the original point cloud,
    and Fig.~\ref{subfig:arm_reduced} is the law density point cloud with $30\%$ points of the original point cloud.
 While the reconstructed surface from the original point cloud,
    which is based on a coarse grid $100\times100\times100$,
    is illustrated in Fig.~\ref{subfig:arm_recons_100},
    the reconstructed surface from the low density point cloud,
    based on a fine grid $190\times190\times190$,
    is demonstrated in Fig.~\ref{subfig:arm_recons_190}.
 It can be seen that,
    though the density of the point cloud is lower than that of the original one,
    the reconstructed surface with finer grid $190\times190\times190$ captures more details and features.

\begin{figure*}[!htb]
\centering
 \subfigure[]
 {
    \includegraphics[width=0.26\linewidth]{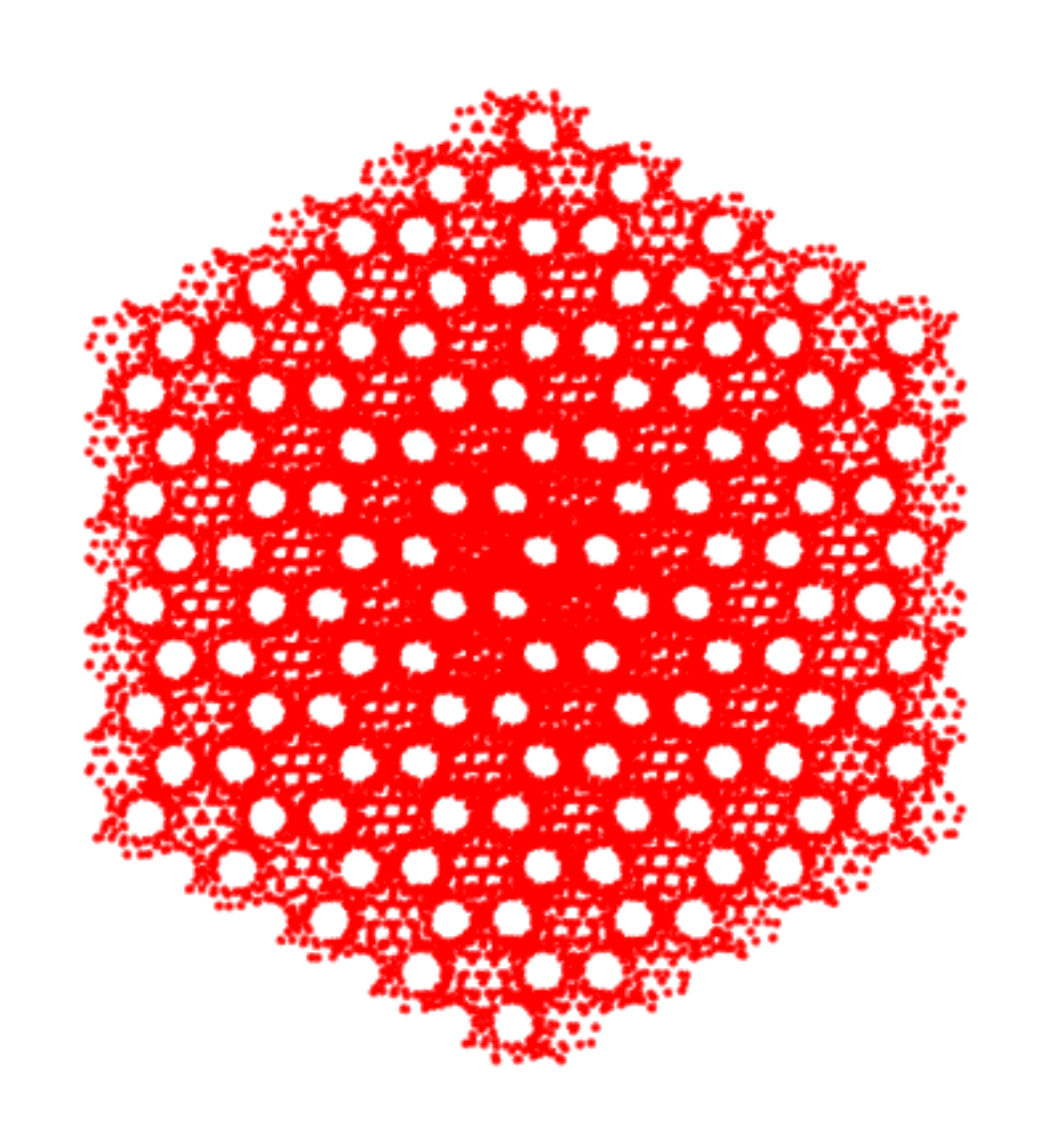}
 }
 \subfigure[]
 {
    \includegraphics[width=0.26\linewidth]{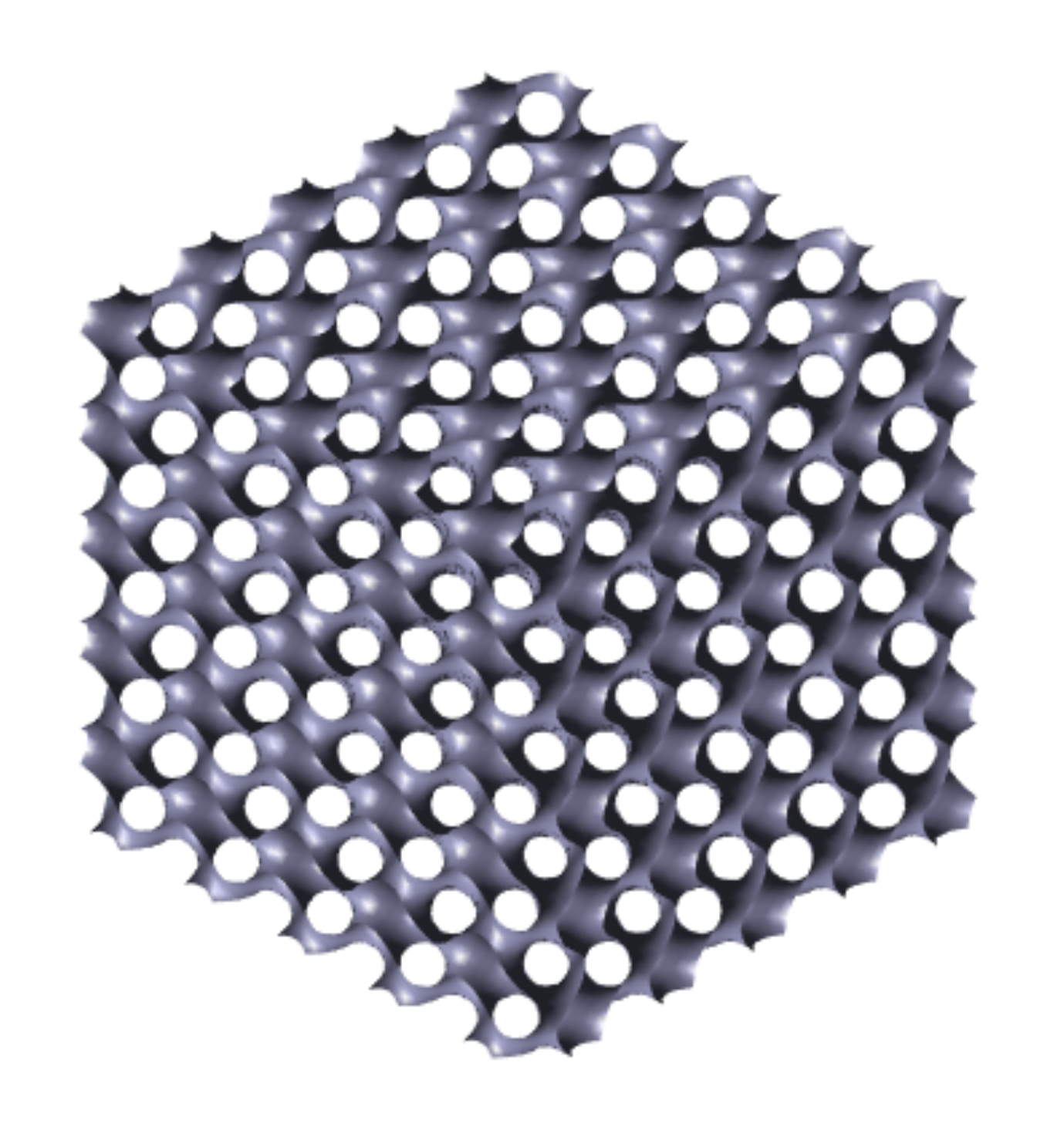}
 }
 \subfigure[]
 {
    \includegraphics[width=0.26\linewidth]{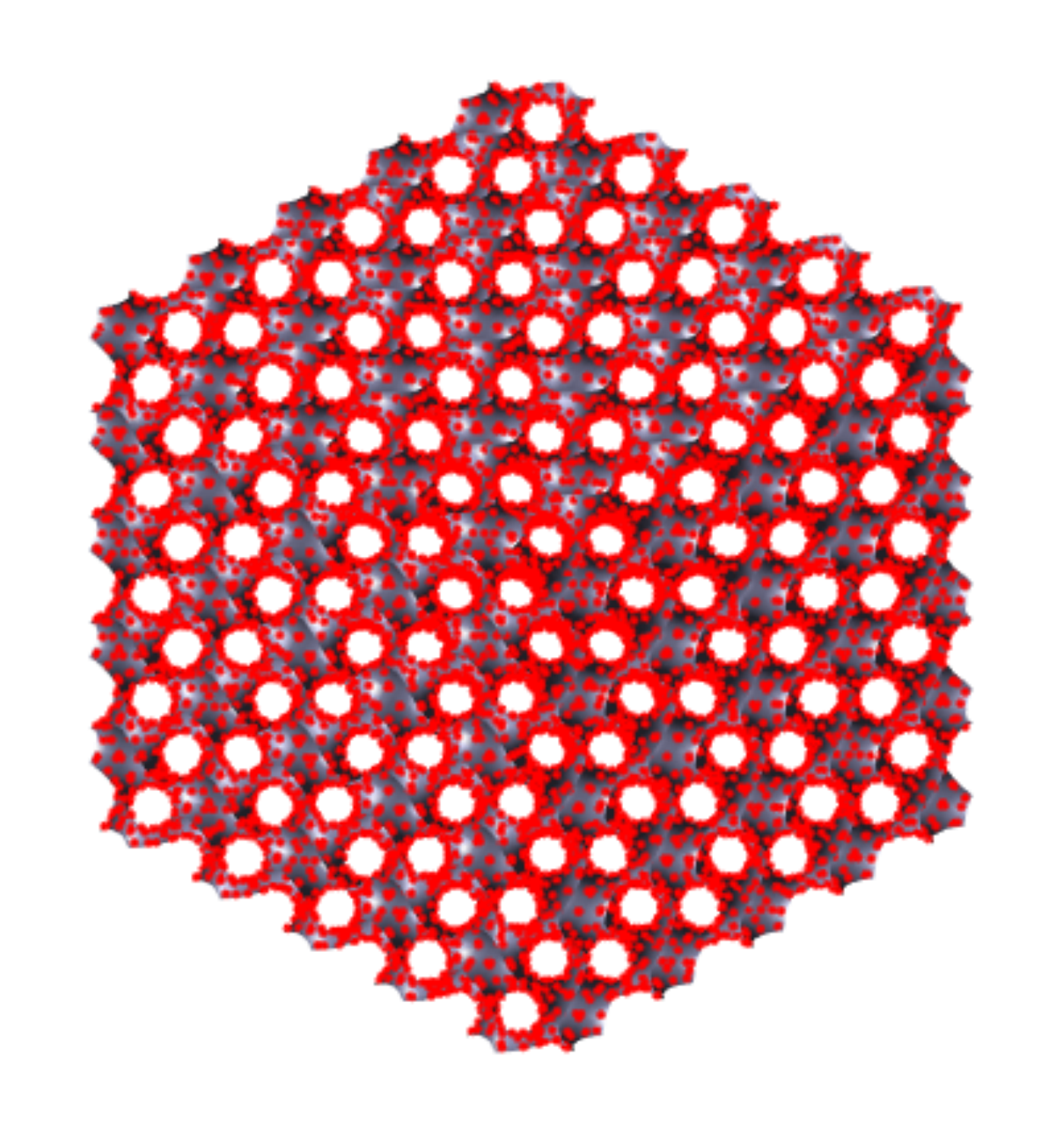}
 }
 \caption
 {
   Surface reconstruction of open surface (porous surfac).
   (a) The point cloud sampled from the \emph{gyroid} surface.
   (b) Reconstructed surface by I-PIA.
   (c) Reconstructed surface with data points superimposed.
 }
\label{Porous}
\end{figure*}

\subsection{Open surface (Porous surface)}
\label{sec:PorousSurfaceOpenSurface}

 I-PIA can perform well for the reconstruction of the open surfaces.
 Fig.~\ref{Porous} depict the reconstructed surface of \emph{gyroid},
 a kind of porous surface.
 From left to right of Fig.~\ref{Porous},
    there are the given point cloud, reconstructed surface, and reconstructed
 surface with the given data superimposed, respectively.
 It can be seen that the implicit B-spline surface reconstructed by I-PIA
    conforms to the given point cloud well.

\begin{figure*}[!htb]
\centering
 \subfigure
 {
    \includegraphics[width=0.20\linewidth]{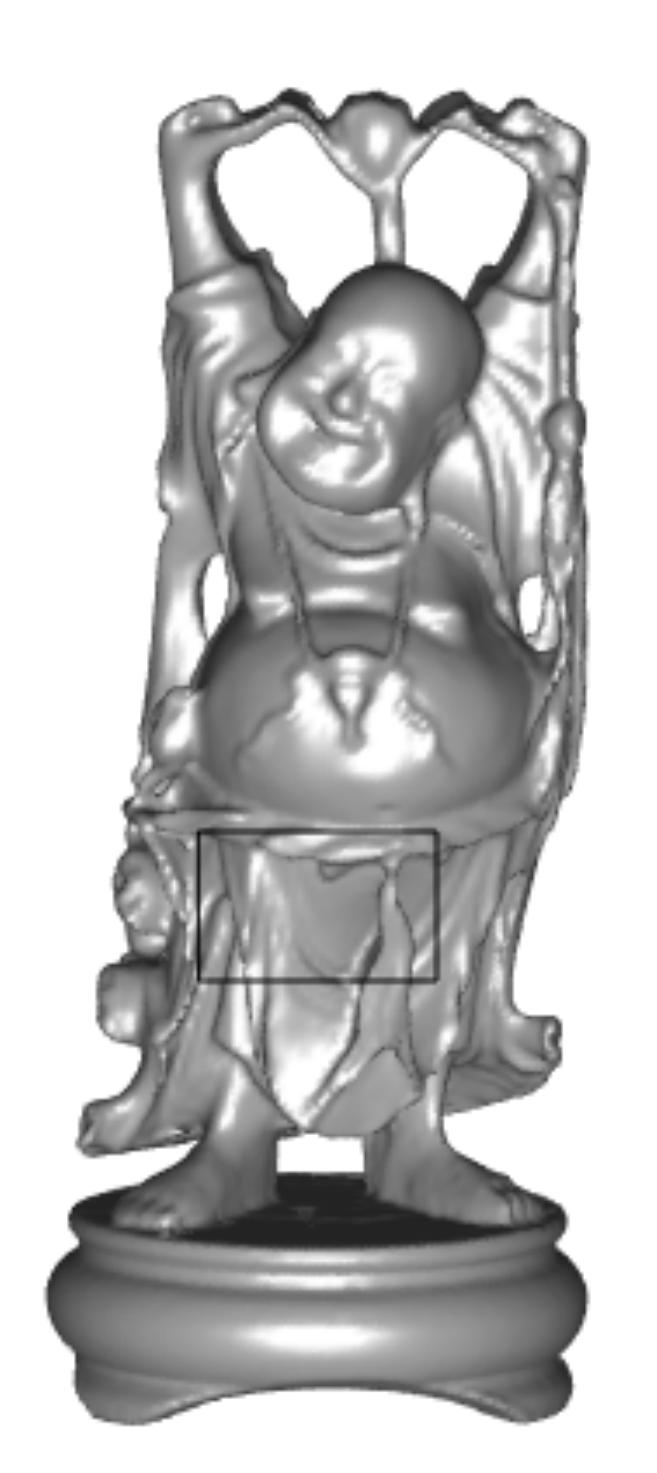}
 }
 \subfigure
 {
    \includegraphics[width=0.20\linewidth]{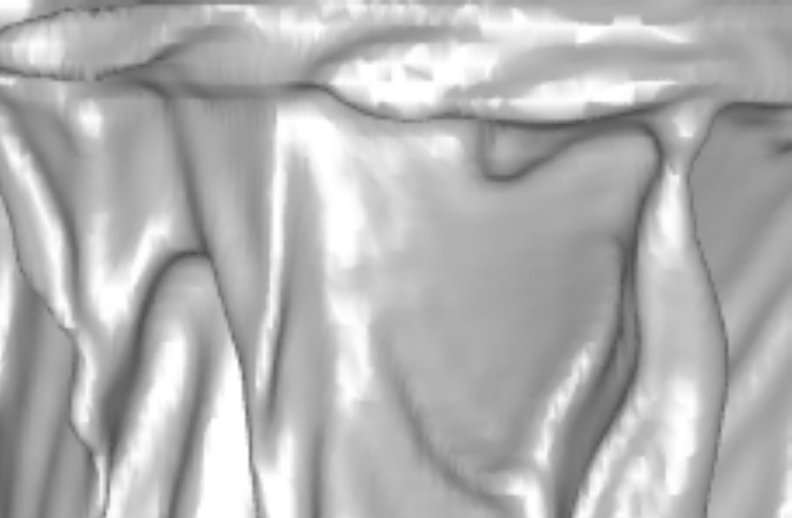}
 }
 \subfigure
 {
    \includegraphics[width=0.23\linewidth]{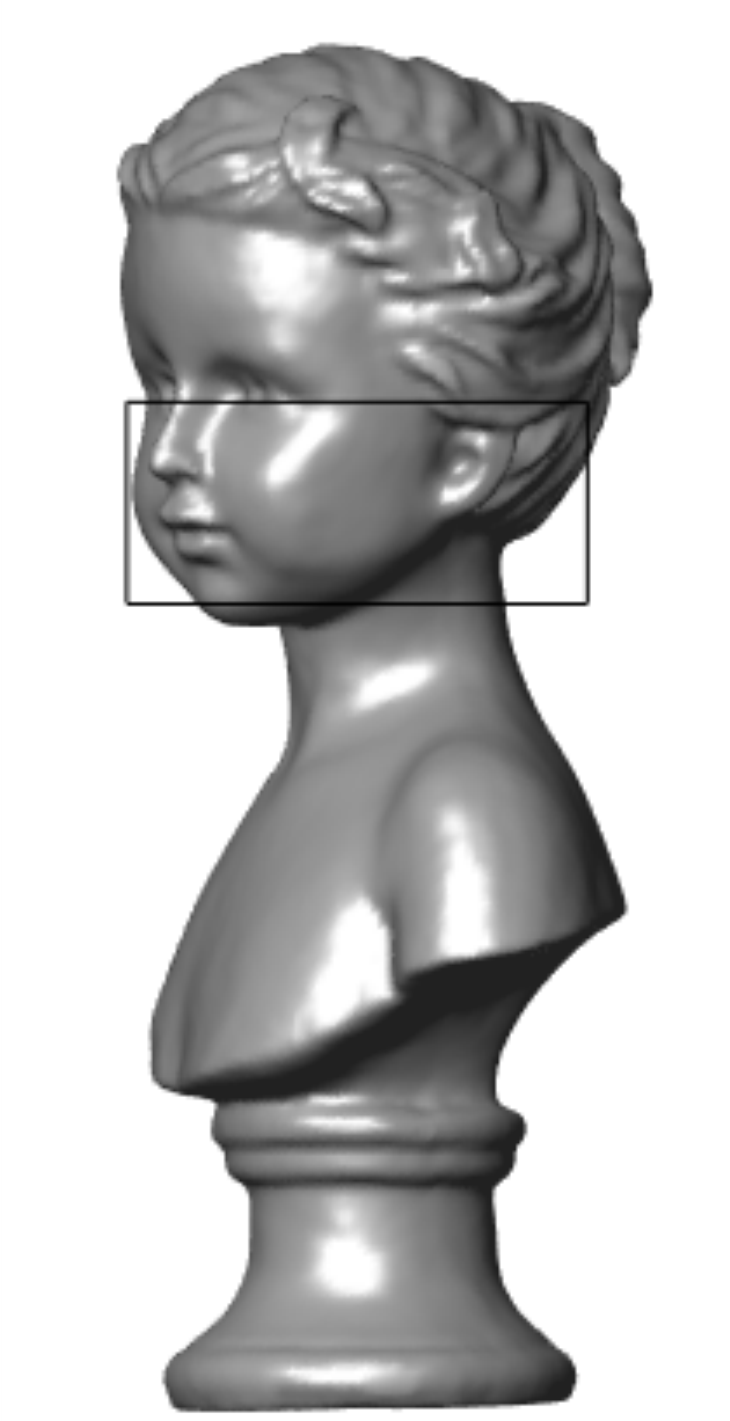}
 }
 \subfigure
 {
    \includegraphics[width=0.23\linewidth]{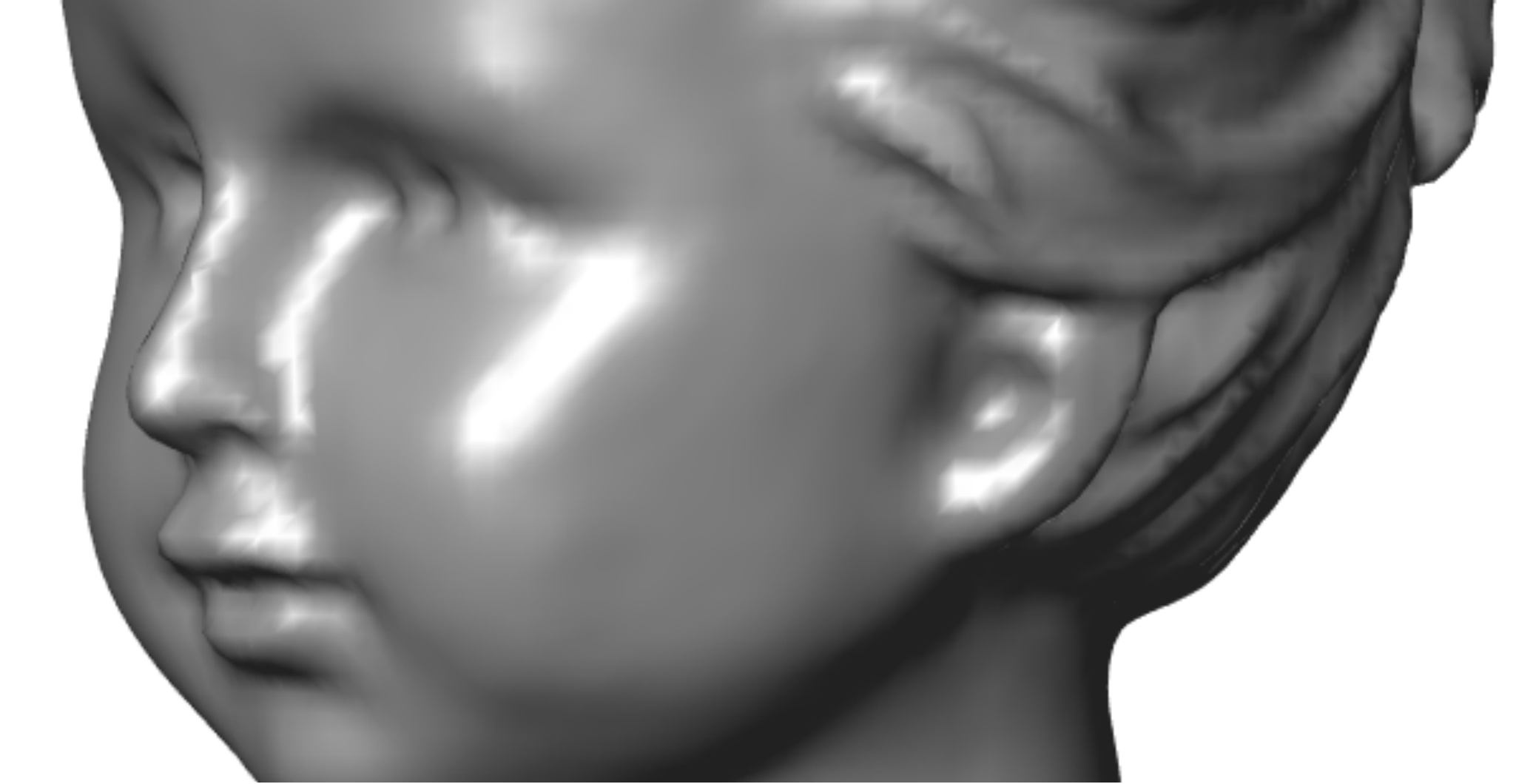}
 }

 \caption
 {Surface reconstruction of \emph{Buddha} and \emph{BU} with sectional enlargements.
 }
\label{enlargement}
\end{figure*}

\subsection{Fine details}
\label{sec:details}
To show the level of quality of the surfaces reconstructed by our
 algorithm, we enlarged some parts of Buddha and BU in Fig.~\ref{enlargement}.
The features of the models are correctly recovered using I-PIA.

%-------------------------------------------------------------------------------
\begin{table}[!htb]
\begin{center}
\caption
{
    \label{tab:2d_data}
    Performance of our method on different planer curves
		  (all timings are measured in seconds)
}
{\small
    \begin{center}
    %\centerline{\small {\bf Table 1}~~Performance of our method on different planer curves
%		  (all timings are measured in seconds).  }\vskip 1mm

    {
    \begin{tabular*}{13cm}{cccccc}
        \toprule

           &Model & \quad  Number of data point &\quad  Grid size    &\quad    Maximum Error &\quad  Time   \\
        \midrule

            &\emph{Flower}      &509    &$30\times30$    &2.924$e$--3    &1.88   \\
            &\emph{Coons curve}        &780    &$30\times30$    &1.2524$e$--3    &2.05   \\
            &\emph{Dolphin}     &371    &$30\times30$    &3.2457$e$--3    &1.67   \\
			&\emph{Butterfly}   &508    &$28\times28$    &2.8881$e$--3    &1.79   \\

        \bottomrule
    \end{tabular*}}
\end{center}}
\end{center}
\end{table}
%------------------------------------------------------------------------------

%------------------------------------------------------------------------------
\begin{table*}
\begin{center}
 \caption
 {
    \label{tab:3d_data}
    Performance of I-PIA on 3D point clouds (all timings are measured in seconds)
 }
{\small
    \begin{center}
    %\centerline{\small {\bf Table 2}~~Performance of our method on different 3D data sets (all timings are measured in seconds).  }\vskip 1mm
%    \label{tab:3d}
    {
    \begin{tabular*}{14.0cm}{ccccccccc}
        \toprule
           &Model &  Data point &  Grid size& Maximum Error      & Method in \cite{Liu2017Implicit}& I-PIA&     \\
        \cmidrule(l){6-7}

              &&&   && Time   & Time \\
        \midrule
				    &\emph{Double-torus}  &766       &$30\times30\times30$     &1.0574$e$--3    &90.07       &3.61   \\[0.8ex]
            &\emph{Torus}         &10100       &$30\times30\times30$   &3.0412$e$--6     &189.13      &8.22      \\[0.8ex]
            &\emph{Gyroid}        &36627    &$50\times50\times50$      &4.8708$e$--6     &11149.28    &34.26     \\[0.8ex]
            &\emph{Elephant}      &52100    &$50\times50\times50$      &9.6839$e$--4     &2773.65     &44.55      \\[0.8ex]
						&\emph{Bunny}         &35947    &$70\times70\times70$      &9.5983$e$--2     &26428.69    &88.24    \\[0.8ex]
						&\emph{Fertility}     &241607    &$70\times70\times70$     &9.5989$e$--5     &25298.49    &902.35      \\[0.8ex]
						
						&\emph{Armadillo}     &51892    &$190\times190\times190$   &1.9831$e$--2     & $>$ 24hr       &2732.32  \\[0.8ex]
						
						&\emph{Buddha}        &80000    &$150\times150\times150$   &7.4594$e$--2     &$>$24hr       &1835.73     \\[0.8ex]
						&\emph{Bu}            &70000    &$160\times160\times160$   &1.1899$e$--3     &$>$24hr       &2028.01     \\[0.8ex]				
        \bottomrule
    \end{tabular*}}
\end{center}}
\end{center}
\end{table*}
%------------------------------------------------------------------------------

\subsection{Performance}
\label{sec:Performance}

 Table~\ref{tab:2d_data} summarizes the running time performance of I-PIA
    on different planner curves.
 Table~\ref{tab:3d_data} shows the comparison of running time performance
    of I-PIA and the state-of-the-art method in \cite{Liu2017Implicit}.
 The computational cost of \cite{Liu2017Implicit} is higher than
    that of I-PIA,
    because in each iteration of the method in  \cite{Liu2017Implicit},
    a linear system needs to be solved,
    while I-PIA avoids solving the linear system.
 The surface reconstructions of \emph{Armadillo}, \emph{Buddha},
    and \emph{BU} are not completed in 24 hours using the method in \cite{Liu2017Implicit}.
 In summary, the implicit curve and surface reconstruction time by I-PIA
        is improved at least one to three orders of magnitude,
        compared with the state-of-the-art method in~\cite{Liu2017Implicit}.

 With I-PIA, the construction of the iterative matrix
		$B$ (Eq.~\pref{eq:b_2d_collocation}~\pref{eq:b_3d_collocation})
    consumes more than 90\% of the running time.
 For example,
	 $88.24$ seconds was spent in the reconstruction of the \emph{bunny}
    model,
    but the construction of the iterative matrix $B$ cost $81.92$ seconds.

\section{Conclusions}
\label{sec:conclude}
In this paper,
  we proposed a novel approach for implicit curve and surface reconstruction based on the progressive-iterative approximation method,
  named I-PIA.
 I-PIA solves the minimization problem with regularization terms naturally
    without any extra computation effort.
Thus, it not only avoids the spurious sheets and artifacts but also reduces
    the computational cost effectively.
Several kinds of  experiments  presented  demonstrate that I-PIA is robust
    to inaccurate distance field, data holes, non-uniform sampling and point noise,
    and thus produces high-quality reconstruction results.

%------------------------------------------------------------------------------
% Acknowledgements:
%------------------------------------------------------------------------------
\section*{Acknowledgments}
This work is supported by the National Natural Science Foundation of China
    under Grant No.61872316,
    and the National Key R\&D Plan of China under Grant No.2016YFB1001501.

% \appendix

% \section{}
% \label{}

%\bibliographystyle{cag-num-names}
%\bibliography{Shape_Retrieval_Neural_Network}

%\chapter{References}
%\label{sec:conclude}

%\section*{References}
\bibliographystyle{elsart-num}
\bibliography{IPIAbibfile}

\begin{thebibliography}{10}
\expandafter\ifx\csname url\endcsname\relax
  \def\url#1{\texttt{#1}}\fi
\expandafter\ifx\csname urlprefix\endcsname\relax\def\urlprefix{URL }\fi

\bibitem{Liu2017Implicit}
Y.~Liu, Y.~Song, Z.~Yang, J.~Deng, Implicit surface reconstruction with total
  variation regularization, Computer Aided Geometric Design 52 (2017) 135--153.

\bibitem{Juttler2002Least}
B.~J{\"u}ttler, A.~Felis, Least-squares fitting of algebraic spline surfaces,
  Advances in Computational Mathematics 17~(1-2) (2002) 135--152.

\bibitem{Rouhani2011Implicit}
M.~Rouhani, A.~D. Sappa, Implicit {B}-spline fitting using the 3{L} algorithm,
  in: 2011 18th IEEE International Conference on Image Processing, IEEE, 2011,
  pp. 893--896.

\bibitem{Rouhani2015Implicit}
M.~Rouhani, A.~D. Sappa, E.~Boyer, Implicit {B}-spline surface reconstruction,
  IEEE Transactions on Image Processing 24~(1) (2014) 22--32.

\bibitem{Yang2005Fitting}
Z.~Yang, J.~Deng, F.~Chen, Fitting unorganized point clouds with active
  implicit {B}-spline curves, The Visual Computer 21~(8-10) (2005) 831--839.

\bibitem{Chen2008Subdivision}
Z.~Chen, X.~Luo, L.~Tan, B.~Ye, J.~Chen, Progressive interpolation based on
  {C}atmull-{C}lark subdivision surfaces, Computer Graphics Forum 27~(7) (2008)
  1823--1827.

\bibitem{Cheng2009Loop}
F.-H.~F. Cheng, F.-T. Fan, S.-H. Lai, C.-L. Huang, J.-X. Wang, J.-H. Yong, Loop
  subdivision surface based progressive interpolation, Journal of Computer
  Science and Technology 24~(1) (2009) 39--46.

\bibitem{Deng2014Progressive}
C.~Deng, H.~Lin, Progressive and iterative approximation for least squares
  {B}-spline curve and surface fitting, Computer-Aided Design 47 (2014) 32--44.

\bibitem{Deng2014Weighted}
C.~Deng, W.~Ma, Weighted progressive interpolation of {L}oop subdivision
  surfaces, Computer-Aided Design 44~(5) (2012) 424--431.

\bibitem{Lin2010Local}
H.~Lin, Local progressive-iterative approximation format for blending curves
  and patches, Computer Aided Geometric Design 27~(4) (2010) 322--339.

\bibitem{Lin1707convergence}
H.~Lin, Q.~Cao, X.~Zhang, The convergence of least-squares progressive
  iterative approximation with singular iterative matrix, arXiv preprint
  arXiv:1707.09109.

\bibitem{Lin2004non-uniform}
H.~Lin, G.~Wang, C.~Dong, Constructing iterative non-uniform {B}-spline curve
  and surface to fit data points, Science in China Series: Information Sciences
  47~(3) (2004) 315--331.

\bibitem{Lin11extended}
H.~Lin, Z.~Zhang, An extended iterative format for the progressive-iteration
  approximation, Computers \& Graphics 35~(5) (2011) 967--975.

\bibitem{Lin2005Totally}
H.-W. Lin, H.-J. Bao, G.-J. Wang, Totally positive bases and progressive
  iteration approximation, Computers \& Mathematics with Applications 50~(3-4)
  (2005) 575--586.

\bibitem{Lu2010Weighted}
L.~Lu, Weighted progressive iteration approximation and convergence analysis,
  Computer Aided Geometric Design 27~(2) (2010) 129--137.

\bibitem{Blinn82generalization}
J.~F. Blinn, A generalization of algebraic surface drawing, ACM Transactions on
  Graphics (TOG) 1~(3) (1982) 235--256.

\bibitem{Muraki1991Volumetric}
S.~Muraki, Volumetric shape description of range data using “blobby model”,
  ACM SIGGRAPH Computer Graphics 25~(4) (1991) 227--235.

\bibitem{Hoppe1992Surface}
H.~Hoppe, T.~DeRose, T.~Duchamp, J.~McDonald, W.~Stuetzle, Surface
  reconstruction from unorganized points, ACM SIGGRAPH Computer Graphics 26~(2)
  (1992) 71--78.

\bibitem{Curless96volumetric}
B.~Curless, M.~Levoy, A volumetric method for building complex models from
  range images, in: Proceedings of the 23rd Annual Conference on Computer
  Graphics and Interactive Techniques, ACM, 1996, pp. 303--312.

\bibitem{Carr2001Reconstruction}
J.~C. Carr, R.~K. Beatson, J.~B. Cherrie, T.~J. Mitchell, W.~R. Fright, B.~C.
  McCallum, T.~R. Evans, Reconstruction and representation of 3{D }objects with
  radial basis functions, in: Proceedings of the 28th Annual Conference on
  Computer Graphics and Interactive Techniques, SIGGRAPH '01, ACM, 2001, pp.
  67--76.

\bibitem{Morse2005Interpolating}
B.~S. Morse, T.~S. Yoo, P.~Rheingans, D.~T. Chen, K.~R. Subramanian,
  Interpolating implicit surfaces from scattered surface data using compactly
  supported radial basis functions, in: ACM SIGGRAPH 2005 Courses, ACM, 2005,
  p.~78.

\bibitem{Turk1998Variational}
G.~Turk, J.~F. O'brien, Variational implicit surfaces, Technical Report
  GIT-GUV-99-15, Georgia Institute of Technology.

\bibitem{Kojekine2003Software}
N.~Kojekine, I.~Hagiwara, V.~Savchenko, Software tools using {CSRBF}s for
  processing scattered data, Computers \& Graphics 27~(2) (2003) 311--319.

\bibitem{ohtake2005multi}
Y.~Ohtake, A.~Belyaev, H.-P. Seidel, Multi-scale and adaptive {CS-RBFS} for
  shape reconstruction from clouds of points, in: Advances in Multiresolution
  for Geometric Modelling, Springer, 2005, pp. 143--154.

\bibitem{Pan2016Compact}
M.~Pan, W.~Tong, F.~Chen, Compact implicit surface reconstruction via low-rank
  tensor approximation, {C}omputer-Aided Design 78 (2016) 158--167.

\bibitem{Ohtake2003Multi-level}
Y.~Ohtake, A.~Belyaev, M.~Alexa, M.~Alexa, G.~Turk, H.-P. Seidel, Multi-level
  partition of unity implicits, ACM Trans. Graph. 22~(3) (2003) 463--470.
\newline\urlprefix\url{http://doi.acm.org/10.1145/882262.882293}

\bibitem{Wang2011Parallel}
J.~Wang, Z.~Yang, L.~Jin, J.~Deng, F.~Chen, Parallel and adaptive surface
  reconstruction based on implicit {PHT}-splines, Computer Aided Geometric
  Design 28~(8) (2011) 463--474.

\bibitem{Kazhdan2006Poisson}
M.~Kazhdan, M.~Bolitho, H.~Hoppe, Poisson surface reconstruction, in:
  Proceedings of the Fourth Eurographics Symposium on Geometry Processing,
  Vol.~7, 2006, pp. 61--70.

\bibitem{Kazhdan2013Screened}
M.~Kazhdan, H.~Hoppe, Screened poisson surface reconstruction, ACM Transactions
  on Graphics (ToG) 32~(3) (2013) 29.

\bibitem{Qi1975numeric}
D.~Qi, Z.~Tian, Y.~Zhang, J.~Feng, The method of numeric polish in curve
  fitting, Acta Mathematica Sinica 18~(3) (1975) 173--184.

\bibitem{deBoor}
C.~de~Boor, How does {A}gee’s smoothing method work, in: Proceedings of the
  1979 army numerical analysis and computers conference, ARO Report, 1979, pp.
  79--3.

\bibitem{Shi2006terative}
S.~Limin, W.~Renhong, An iterative algorithm of {NURBS} interpolation and
  approximation, Journal of Mathematical Research and Exposition 26~(4) (2006)
  735--743.

\bibitem{Lin2018Survey}
H.~Lin, T.~Maekawa, C.~Deng, Survey on geometric iterative methods and their
  applications, Computer-Aided Design 95 (2018) 40--51.

\bibitem{Fasshauer2007Meshfree}
G.~E. Fasshauer, Meshfree approximation methods with MATLAB, Vol.~6, World
  Scientific, 2007.

\bibitem{james1978generalised}
M.~James, The generalised inverse, The Mathematical Gazette 62~(420) (1978)
  109--114.

\end{thebibliography}

\end{document}